\documentclass{amsart}
\usepackage{amssymb,amsmath,accents}
\newcommand{\RR}{{\mathbb R}}

\newcommand{\AH}{{\mathbb A}}
\newcommand{\A}{{\mathcal A}}
\newcommand{\Diff}{{\rm Diff}}
\newcommand{\D}{{\mathcal D}}
\newcommand{\x}{\langle x\rangle}

\newcommand{\e}{\varepsilon}

\newcommand{\Lap}{\Delta}
\newcommand{\de}{\delta}
\newcommand{\del}{\partial}
\newcommand{\om}{\omega}

\newcommand{{\loc}}{{\ell\mathrm oc}}

\def\meanint{{\diagup\hskip -.42cm\int}}

\newtheorem{theorem}{Theorem}[section]
\newtheorem{lemma}{Lemma}[section]
\newtheorem{proposition}{Proposition}[section]
\newtheorem{corollary}{Corollary}[section]

\newtheorem{definition}{Definition}[section]
\newtheorem{remark}{Remark}[section]
\newtheorem{example}{Example}[section]
\newtheorem*{main}{Theorem}

\begin{document}

\title[Asymptotic Diffeomorphisms]
{Groups of Asymptotic Diffeomorphisms}
\author{Robert McOwen}
\address{Northeastern University}
\email{r.mcowen@neu.edu}


\author{Peter Topalov }
\address{Northeastern University}
\email{p.topalov@neu.edu}

\date{October 27, 2015 \qquad{\it Email}: r.mcowen@neu.edu\ }

\begin{abstract}
We consider classes of diffeomorphisms of Euclidean space with partial asymptotic expansions at infinity; the remainder term lies in a weighted Sobolev space whose properties at infinity fit with the desired application. We show that two such classes of asymptotic diffeomorphisms form topological groups under composition. As such, they can be used in the study of fluid dynamics according to the approach of V.\ Arnold \cite{A}. 

\smallskip\noindent
{\bf Keywords.} Groups of diffeomorphisms, asymptotic expansions, weighted Sobolev spaces, Camassa-Holm equation, Euler equation.

\end{abstract}

\maketitle

\addtocounter{section}{-1}
\section{Introduction}

A modern development in fluid dynamics is to view the motion of an incompressible fluid as a geodesic flow on a group of diffeomorphisms of the underlying physical space. This approach was initiated by V. Arnold \cite{A} and further developed by Ebin \& Marsden \cite{EM} and Bourguignon and Brezis \cite{BB} to obtain well-posedness of initial-value problems associated with the Euler and Navier-Stokes equations. In these papers, the underlying physical space was compact (a compact manifold with or without boundary). Subsequently, Cantor \cite{C} used this approach to study the Euler equations on $\RR^d$ by considering diffeomorphisms $\phi:\RR^d\to\RR^d$ of the form
\begin{equation}\label{Id+f}
\phi=Id + f,
\end{equation}
where $Id$ is the identity map and the function $f$ is in a weighted Sobolev space that requires $f$ to  decay rapidly at infinity.
However, one would like to consider diffeomorphisms of the form (\ref{Id+f}) where $f$ is bounded but not required to decay rapidly.
Moreover, if the initial condition has asymptotics at infinity, one would like to know that the solution has similar asymptotics at infinity (with coefficients depending on $t$). To make these improvements, we require additional structure.

In this paper, we study groups of diffeomorphism on $\RR^d$ of the form
\begin{equation}\label{Id+phi}
\phi=Id + u,
\end{equation}
where $u$ is taken from a function space that we call an {\it asymptotic space}: these consist of  bounded maps on $\RR^d$ having a {\it partial asymptotic expansion at infinity} of the form
\begin{equation}\label{u-asymptotics}
u(x)=a_0(\theta)+\frac{a_1(\theta)}{r}+\cdots + \frac{a_N(\theta)}{r^N}+f_N(x) \quad\hbox{for}\ r=|x|>R,
\end{equation}
where $\theta=x/|x|$, the functions $a_0,\dots,a_N$ lie in certain Sobolev spaces on the unit sphere $S^{d-1}$, and the
 remainder function $f_N$ belongs to a function space ${\mathcal R}_N$ which ensures that 
 \begin{equation}\label{f-remainder}
| f_N(x)|=o(|x|^{-N}) \quad\hbox{ as } |x|\to\infty. 
 \end{equation}
The remainder space ${\mathcal R}_N$ will be a weighted Sobolev space, but there are different possibilities: the choice will depend upon the application, since it must be compatible with the equations being studied. In this paper we shall consider as remainder space two different classes of weighted Sobolev spaces. In one class, that we shall denote by $H_\delta^{m,p}(\RR^d)$,  the functions have derivatives up to order $m$ that are in $L_N^p(\RR^d):=\{f\in L_{\loc}^p(\RR^d): (1+|x|)^\delta f(x)\in L^p(\RR^d)\}$. In the other class of weighted Sobolev spaces, that we shall denote by $W_{\delta}^{m,p}(\RR^d)$,  the derivatives of functions satisfy $D^\alpha f\in L^p_{\delta+|\alpha|}(\RR^d)$ for all $|\alpha|\leq m$. Neither of these classes of weighted Sobolev spaces is new to the literature, but for convenience we shall define them and summarize their properties in Section 1; our exposition is self-contained, with proofs provided in the Appendix.
 
In Section 2 we give the formal definition of our  asymptotic spaces on $\RR^d$. The asymptotic space  $\AH_{N}^{m,p}(\RR^d)$ uses the weighted Sobolev space $H_N^{m,p}(\RR^d)$ as the remainder space; in fact, (\ref{u-asymptotics}) is satisfied by $f\in H_N^{m,p}(\RR^d)$ provided $m>d/p$.
On the other hand, the asymptotic space ${\mathcal A}_{N}^{m,p}(\RR^d)$ uses $W_{\gamma_N}^{m,p}(\RR^d)$ 
for certain values of $\gamma_N$ as the remainder space; in fact, (\ref{u-asymptotics}) is satisfied by $f\in W_{\gamma_N}^{m,p}(\RR^d)$ 
 if $m>d/p$ and $\gamma_N\geq N-d/p$, but we
will need to impose additional constraints on $\gamma_N$ for technical reasons (see Sections 2 and 3).
 At times it is useful to consider functions with $a_0=\cdots=a_{n-1}=0$ for an integer $ n\leq N$; we denote the corresponding spaces by $\AH_{n,N}^{m,p}(\RR^d)$ and
${\mathcal A}_{n,N}^{m,p}(\RR^d)$, and identify $\AH_{0,N}^{m,p}(\RR^d)=\AH_{N}^{m,p}(\RR^d)$ and
${\mathcal A}_{0,N}^{m,p}(\RR^d)={\mathcal A}_{N}^{m,p}(\RR^d)$.
Our primary interest in these asymptotic spaces is to control the behavior of diffeomorphisms at infinity. However, in Section 3, we consider an application of the asymptotic spaces to the Helmholtz decomposition of vector fields on $\RR^d$; this requires an analysis of the inverse of the Laplacian, which is important in many applications, including our study \cite{MT2} of Euler's equation on $\RR^d$.
 
In Section 4, we introduce and study the associated spaces of diffeomorphisms ${\AH\mathcal D}_{n,N}^{m,p}$
and ${\mathcal {AD}}_{n,N}^{m,p}$, i.e.\ diffeomorphisms $\phi:\RR^d\to\RR^d$ that are of the form (\ref{Id+phi}) where the 
components of $u$ are in ${\AH}_{n,N}^{m,p}$ or ${\mathcal A}_{n,N}^{m,p}$ respectively. The main result of this paper is the following:
\begin{main} For integers $m> 2+d/p$ and $N\geq n\geq 0$,  ${\mathcal {AD}}_{n,N}^{m,p}$ is a topological group under composition; in fact, composition is  a $C^1$-map 
${\mathcal {AD}}_{n,N}^{m,p}\times {\mathcal {AD}}_{n,N}^{m-1,p}\to {\mathcal {AD}}_{n,N}^{m-1,p}$ and
the inverse map $\phi \mapsto \phi^{-1}$ is  a $C^1$-map 
${\mathcal {AD}}_{n,N}^{m,p} \to {\mathcal {AD}}_{n,N}^{m-1,p}$. 
The analogous statements also apply to ${\AH\mathcal D}_{n,N}^{m,p}$.
\end{main}

The fact that ${\AH\mathcal D}_{n,N}^{m,p}$ and ${\mathcal {AD}}_{n,N}^{m,p}$ are topological groups allows us to use them to study the asymptotics of various fluid flows on $\RR^d$. With $d=1$, for example, the Camassa-Holm equation \cite{CH} is a completely integrable equation that has attracted considerable attention recently. From the differential geometric point of view, Misiolek \cite{Mi} showed that the equation can be realized as the geodesic flow for a certain metric on the Bott-Virasaro group. Moreover, Constantin \cite{Co} studied initial-value problems for Camassa-Holm on $\RR$ by using a group of diffeomorphisms of the form (\ref{Id+f}) with $f=o(|x|^{-3/2})$ as $|x|\to\infty$. In \cite{MT}, we used the groups of asymptotic diffeomorphims on  $\RR$  to show that the initial-value problem for the Camassa-Holm equation is well-posed with asymptotics at infinity.
For example, if the initial condition $u_0$ is in ${\AH}_{n,N}^{m,2}$ for $m\geq 3$ and $N\geq n\geq 0$, then there is a unique solution $u$ of Camassa-Holm in $C^0([0,T],{\AH}_{n,N}^{m,2})\cap C^1([0,T],{\AH}_{n,N}^{m-1,2})$. In particular, with $n=0$, this means that $f$ in  (\ref{Id+f}) is not required to decay at infinity.

With $d\geq 2$, we have also used the groups of asymptotic diffeomorphisms ${\mathcal {AD}}_{n,N}^{m,p}$ to study Euler's equation for a velocity field ${\bf u}$ and pressure ${\rm p}$ of an incompressible fluid  on $\RR^d$ with external forcing ${\bf f}$. In fact, in 
 \cite{MT2} we show that if $m>2+d/p$,  $1\leq N \leq d-1$, and  ${\bf f}\in C([0,T],\A^{m+1,p}_{1,N})$, then, for any ${\bf u}_0\in \A^{m,p}_{1,N}$ with ${\rm div}\,{\bf u}_0=0$, 
there exists $\tau\in (0,T]$ and functions ${\bf u}\in C^0([0,\tau],\A^{m,p}_{1,N})\cap C^1([0,\tau],\A^{m-1,p}_{1,N})$, and ${\rm p}\in C([0,\tau],H^{m+1,p}_{\loc})$ satisfying Euler's equation; the function {\bf u} is unique. A crucial step in this analysis is the use of the Euler projector as constructed in Section 3 of the present paper.

It is also possible to include log-terms in the asymptotics; see Appendix \ref{B} for definitions and the analog of the above Theorem. Using log-terms, one can avoid the restriction $N\leq d-1$ in the application to Euler's equation; cf.\ \cite{MT2}.

Note that asymptotic solutions of nonlinear evolution equations were studied by many authors using different methods. For example, we mention Menikoff \cite{M}, Bondareva \& Shubin \cite{BS1}, \cite{BS2}, Kenig, Ponce, and Vega \cite{KPV}, and Kappeler, Perry, Shubin, \& Topalov \cite{KPST}.

\section{ Weighted Sobolev Spaces on $\RR^d$}

Let $\langle x\rangle=\sqrt{|x|^2+1}.$
For $1<p<\infty$, $\delta\in\RR$, and a nonnegative integer $m$, we define the Banach spaces $H_\delta^{m,p}(\RR^d)$ and $W_\delta^{m,p}(\RR^d)$ to be the closures of $C_0^\infty(\RR^d)$ in the respective norms:
\begin{equation}\label{H-norm}
\|f\|_{H_\delta^{m,p}}=\sum_{|\alpha|\leq m} \| \langle x\rangle^{\delta} D^\alpha f \|_{L^p},
\end{equation}
\begin{equation}\label{W-norm}
\|f\|_{W_\delta^{m,p}}=\sum_{|\alpha|\leq m} \| \langle x\rangle^{\delta+|\alpha |} D^\alpha f \|_{L^p}.
\end{equation}
Notice that $H_\delta^{0,p}(\RR^d)=W_\delta^{0,p}(\RR^d)$ is just a weighted $L^p$-space that may be  denoted  by $L_\de^p(\RR^d)$. If we let $H^{m,p}_{\loc}(\RR^d)$ be functions whose first $m$ derivatives are  $L^p$-integrable over compact subsets of $\RR^d$,
then we can use mollifiers to show that
$
H_\delta^{m,p}(\RR^d)=\{f\in H_{\loc}^{m,p}(\RR):\|f\|_{H_\delta^{m,p}}<\infty\}
$
and
$
W_\delta^{m,p}(\RR^d)=\{f\in H_{\loc}^{m,p}(\RR):\|f\|_{W_\delta^{m,p}}<\infty\}.
$
Note also that $W_\delta^{m,p}(\RR^d)\subset H^{m,p}_\delta(\RR^d)$.

These weighted Sobolev spaces enjoy the following properties: 
\begin{lemma}\label{le:H-properties}  
Let $1<p<\infty$, $\delta\in\RR$, and $m$ be a nonnegative integer.
\begin{enumerate}
\item[(a)]\label{differentiationiscontinuous-1} For $m\geq 1$, 
$
f\mapsto \frac{\partial f}{\partial x_j}$   defines a continuous map $H_\delta^{m,p}(\RR^d)\to H^{m-1,p}_{\delta}(\RR^d).$
\item[(b)]\label{multiplicationby<x>^k-1} For  $m\geq 0$ and $\gamma\in\RR$, 
$f\mapsto \langle x\rangle^{-\gamma}f$  defines a continuous map $H_\delta^{m,p}(\RR^d)\to H^{m,p}_{\delta+\gamma}(\RR^d).$
\item[(c)]\label{SobolevInequality-1} For $mp<d$, if $f\in H^{m,p}_{\de}(\RR^d)$ then
$f\in L^q_{\de}(\RR^d)$ for all $q\in[p,dp/(d-mp)]$, and
\[
\|f\|_{L^{q}_{\de}}\leq C\, \|f\|_{H^{m,p}_{\de}}, \quad\hbox{where}\ C=C(d,m,p,q,\delta).
\]
For $mp=d$, the same conclusions hold for all $q\in[p,\infty)$.
\item[(d)]\label{Morrey-1} For $mp>d$, if $f\in H_{\delta}^{m,p}(\RR^d)$  then $f\in C^{k}(\RR^d)$ for all $k<m-(d/p)$, and
\[
\sup_{x\in \RR^d} \left( \langle x \rangle^{\delta} |D^\alpha f(x)|\right)
\leq C\,\|f\|_{H_{\de}^{m,p}} \quad\hbox{for all \ $0\leq |\alpha|\leq k$, where $C=C(d,m,p,k,\delta)$.}
\]
 In fact, for  all \ $0\leq |\alpha|\leq k$, we have
\[
 |x|^{\delta}\, |D^\alpha f(x)| \to 0\quad\hbox{as}\ |x|\to\infty.
\]
 \end{enumerate}
\end{lemma}

\begin{lemma}\label{le:W-properties}  
Let $1<p<\infty$, $\delta\in\RR$, and $m$ be a nonnegative integer.
\begin{enumerate}
\item[(a)]\label{differentiationiscontinuous-2} For $m\geq 1$, 
$f\mapsto \frac{\partial f}{\partial x_j}$  defines a continuous map $W_\delta^{m,p}(\RR^d)\to W^{m-1,p}_{\delta+1}(\RR^d).$
\item[(b)]\label{multiplicationby<x>^k-2} For  $m\geq 0$ and $\gamma\in\RR$, 
$f\mapsto \langle x\rangle^{-\gamma}f $ defines a continuous map $W_\delta^{m,p}(\RR^d)\to W^{m,p}_{\delta+\gamma}(\RR^d).$
\item[(c)]\label{SobolevInequality-2} For $mp<d$, if $f\in W^{m,p}_{\de-\frac{d}{p}}(\RR^d)$ then
$f\in L^q_{\de-\frac{d}{q}}(\RR^d)$ for all $q\in[p,dp/(d-mp)]$, and
\[
\|f\|_{L^{q}_{\de-\frac{d}{q}}}\leq C \|f\|_{W^{m,p}_{\de-\frac{d}{p}}}, \quad \hbox{where $C=C(d, m,p,q)$}.
\]
For $mp=d$, the same conclusions hold for all $q\in[p,\infty)$.
\item[(d)]\label{Morrey-2} For $mp>d$, if $f\in W_{\delta-\frac{d}{p}}^{m,p}(\RR^d)$  then $f\in C^{k}(\RR^d)$ for all $k<m-(d/p)$, and
\[
\sup_{x\in \RR^d} \left( \langle x \rangle^{\delta+|\alpha|} |D^\alpha f(x)|\right)
\leq C\,\|f\|_{W_{\de-\frac{d}{p}}^{m,p}} \quad\hbox{for all \ $0\leq |\alpha|\leq k$, where $C=C(m,p,k,\delta)$.}
\]
 In fact, for  all \ $0\leq |\alpha|\leq k$, we have
\[
 |x|^{\delta+|\alpha|}\, |D^\alpha f(x)| \to 0\quad\hbox{as}\ |x|\to\infty.
\]
 \end{enumerate}
\end{lemma}

\medskip\noindent
For both lemmas, the properties (a) and (b) are obvious; properties (c) and (d) are proved in the Appendix.

Using these lemmas, we can prove the following results about pointwise multiplication:

\begin{proposition}\label{pr:H-multiplication}  
For $(m+\ell-k)\,p>d$ where $0\leq k\leq \ell\leq m$, pointwise multiplication
\[
(f,g)\mapsto fg \ \ \hbox{defines a continuous map}  \ \ 
H_{\delta_1}^{m,p}(\RR^d) \times H_{\delta_2}^{\ell,p}(\RR^d) \to 
H^{k,p}_{\delta_1+\delta_2}(\RR^d).
\]
In fact, there is a constant $C=C(d,m,\ell,k,p,\delta_1,\delta_2)$ such that
\[
\| fg \|_{H^{k,p}_{\de_1+\de_2}} \leq C\, \|f\|_{H_{\delta_1}^{m,p}}\,\|g\|_{H_{\delta_2}^{\ell,p}}
\quad\hbox{for all $f\in H_{\delta_1}^{m,p}$ and $g\in H_{\delta_2}^{\ell,p}$.}
\]
\end{proposition}

\medskip \noindent
{\bf Proof.} To begin with, it is easy to check that the following weighted H\"older inequality holds:
\begin{equation}\label{eq:weightedHolder}
\|fg\|_{L^p_{\de}}\leq \|f\|_{L^{q_1}_{\de_1}}\,  \|g\|_{L^{q_2}_{\de_2}},
\quad\hbox{if}\ \ \frac{1}{p}=\frac{1}{q_1}+\frac{1}{q_2}\ \ \hbox{and} \ \ \de=\de_1+\de_2.
\end{equation}
Consequently, we will know that $(f,g)\mapsto fg$ defines a continuous map 
$H_{\delta_1}^{m,p}(\RR^d) \times H_{\delta_2}^{\ell,p}(\RR^d) \to 
L^{p}_{\delta_1+\delta_2}(\RR^d)$, i.e.\ the proposition holds for $k=0$,
provided we can find $1\leq q_1,q_2\leq \infty$ so that
\begin{equation}\label{eq:H-embeddings}
H_{\delta_1}^{m,p}(\RR^d)\subset L^{q_1}_{\de_1}(\RR^d), \quad
H_{\delta_2}^{\ell,p}(\RR^d)\subset L^{q_2}_{\de_2}(\RR^d), \quad 
 \frac{1}{p}=\frac{1}{q_1}+\frac{1}{q_2}.
\end{equation}

To obtain (\ref{eq:H-embeddings}), 
let us first  assume that $mp<d$, so we also have $\ell p<d$. According to Lemma \ref{le:H-properties} (c), we have
 \[
 \begin{aligned}
 H_{\delta_1}^{m,p}(\RR^d)&\subset 
 L^{q_1}_{\de_1}(\RR^d)\quad\hbox{provided}\ p\leq q_1\leq \frac{dp}{d-mp} \\
 H_{\delta_2}^{\ell,p}(\RR^d)&\subset 
 L^{q_2}_{\de_2}(\RR^d)\quad\hbox{provided}\ p\leq q_2\leq \frac{dp}{d-\ell p}.
 \end{aligned}
 \]
Now it is clear that the function $f(q_1,q_2)=q_1^{-1}+q_2^{-1}$ takes on all values between  $2/p$ and
\[
\frac{d-mp}{dp}+ \frac{d-\ell p}{dp}=\frac{2}{p}-\frac{m+\ell}{d},
\]
so whether $f(q_1,q_2)$ ever equals $p^{-1}$ is determined by whether
\begin{equation}\label{condition-on-p}
\frac{2}{p}-\frac{m+\ell}{d}\leq \frac{1}{p} \leq \frac{2}{p}.
\end{equation}
The second inequality is trivial but the first holds precisely when $(m+\ell)p\geq d$.

How does this result change when $mp\geq d$? For $mp=d$ we need to require $q_1<\infty$, which translates into a strict inequality in (\ref{condition-on-p}), so we require $(m+\ell)p> d$. For  $mp>d$, then by Lemma 1(d) we can take $q_1=\infty$ and $q_2=p$. Thus we always have (\ref{eq:H-embeddings}) under the assumption $(m+\ell)p> d$. This proves the proposition for $k=0$.

Now, to prove the proposition for all $0\leq k\leq \ell$, we must show that $D^\alpha(fg)\in L^p_{\de_1+\de_2}$ for all $ |\alpha|\leq k$. But if we use the Leibniz rule to write
\[
D^\alpha(fg)=
\sum_{\beta\leq \alpha} \left(\begin{matrix} \alpha \\ \beta \end{matrix} \right) (D^\beta f) (D^{\alpha-\beta} g),
\]
and we observe that $D^\beta f \in H^{m-|\beta|}_{\de_1}(\RR^d)$ and
$D^{\alpha-\beta} g \in H^{\ell-|\alpha-\beta|}_{\de_2}(\RR^d)$, then we can use
(\ref{eq:H-embeddings}) provided $m-|\beta|+\ell-|\alpha-\beta|>d/p$. But this is guaranteed since $|\beta|+|\alpha-\beta|=|\alpha|\leq k$. 
 \hfill$\Box$

\begin{proposition}\label{pr:W-multiplication} 
For $(m+\ell-k)p>d$ where $0\leq k\leq \ell\leq m$, pointwise multiplication
\[
(f,g)\mapsto fg \ \ \hbox{defines a continuous map}  \ \ 
W_{\delta_1-\frac{d}{p}}^{m,p}(\RR^d) \times W_{\delta_2-\frac{d}{p}}^{\ell,p}(\RR^d) \to 
W^{k,p}_{\delta_1+\delta_2-\frac{d}{p}}(\RR^d).
\]
In fact,  there is a constant $C=C(d,m,\ell,k,p)$ such that
\[
\| fg \|_{W^{k,p}_{\de_1+\de_2-\frac{d}{p}}} \leq 
C\, \|f\|_{W_{\delta_1-\frac{d}{p}}^{m,p}}\,\|g\|_{W_{\delta_2-\frac{d}{p}}^{\ell,p}}
\quad\hbox{ for all $f\in W_{\delta_1-\frac{d}{p}}^{m,p}$ and
$g\in W_{\delta_2-\frac{d}{p}}^{\ell,p}$.}
\]
\end{proposition}

\medskip \noindent
{\bf Proof.} As a special case of (\ref{eq:weightedHolder}) we have
\[
\|fg\|_{L^p_{\de-\frac{d}{p}}}\leq \|f\|_{L^{q_1}_{\de_1-\frac{d}{q_1}}}\,  \|g\|_{L^{q_2}_{\de_2-\frac{d}{q_2}}},
\quad\hbox{if}\ \ \frac{1}{p}=\frac{1}{q_1}+\frac{1}{q_2}\ \ \hbox{and} \ \ \de=\de_1+\de_2.
\]
Using this, we will know that $(f,g)\mapsto fg$ defines a continuous map
$W_{\delta_1-\frac{d}{p}}^{m,p}(\RR^d) \times W_{\delta_2-\frac{d}{p}}^{\ell,p}(\RR^d) \to 
L^{p}_{\delta_1+\delta_2-\frac{d}{p}}(\RR^d)$
provided we can find $1\leq q_1,q_2\leq \infty$ so that
\begin{equation}\label{eq:W-embeddings}
W_{\delta_1-\frac{d}{p}}^{m,p}(\RR^d)\subset L^{q_1}_{\de_1-\frac{d}{q_1}}(\RR^d), \quad
W_{\delta_2-\frac{d}{p}}^{\ell,p}(\RR^d)\subset L^{q_2}_{\de_2-\frac{d}{q_2}}(\RR^d), \quad 
 \frac{1}{p}=\frac{1}{q_1}+\frac{1}{q_2}.
\end{equation}
Let us assume first that $mp<d$, so we also have $\ell p<d$. According to Lemma \ref{le:W-properties} (c), we have
 \[
 \begin{aligned}
 W_{\delta_1-\frac{d}{p}}^{m,p}(\RR^d)&\subset 
 L^{q_1}_{\de_1-\frac{d}{q_1}}(\RR^d)\quad\hbox{provided}\ p\leq q_1\leq \frac{dp}{d-mp} \\
 W_{\delta_2-\frac{d}{p}}^{\ell,p}(\RR^d)&\subset 
 L^{q_2}_{\de_2-\frac{d}{q_2}}(\RR^d)\quad\hbox{provided}\ p\leq q_2\leq \frac{dp}{d-\ell p}.
 \end{aligned}
 \]
For the same reasons as in the proof of Proposition \ref{pr:H-multiplication}, this is possible when $(m+\ell)p\geq d$.
The case $mp\geq d$ also follows as  in the proof of Proposition \ref{pr:H-multiplication}.

Now, to prove the proposition, we must show that $\langle x\rangle^{|\alpha|}D^\alpha(fg)\in L^p_{\de_1+\de_2-(d/p)}$ for all $ |\alpha|\leq k$. But if we write
\[
\begin{aligned}
\langle x\rangle^{|\alpha|}D^\alpha(fg)=&
\langle x\rangle^{|\alpha|}\sum_{\beta\leq \alpha} \left(\begin{matrix} \alpha \\ \beta \end{matrix} \right) 
(D^\beta f) (D^{\alpha-\beta} g) \\
=& \sum_{\beta\leq \alpha} \left(\begin{matrix} \alpha \\ \beta \end{matrix} \right)
\left(\langle x\rangle^{|\beta|}D^\beta f \right)\left(\langle x\rangle^{|\alpha-\beta|}D^{\alpha-\beta} g\right),
\end{aligned}
\]
and we observe that $\langle x\rangle^{|\beta|}D^\beta f \in W^{m-|\beta|}_{\de_1-\frac{d}{p}}(\RR^d)$ and
$\langle x\rangle^{|\alpha-\beta|}D^{\alpha-\beta} g \in W^{\ell-|\alpha-\beta|}_{\de_2-\frac{d}{p}}(\RR^d)$, then we can use
(\ref{eq:W-embeddings}) provided $m-|\beta|+\ell-|\alpha-\beta|>d/p$. But this is guaranteed since $|\beta|+|\alpha-\beta|=|\alpha|\leq k$. 
 \hfill$\Box$
 
 \medskip
In the next section, we shall also need to consider  Sobolev spaces $H^{m,p}(S^{d-1})$ on the unit sphere $S^{d-1}$ in $\RR^d$. The boundedness of multiplication on these spaces can be found in the literature or  easily  derived using the H\"older inequality and Sobolev embedding on $S^{d-1}$ as in the proofs above. We record here the result.

 \begin{proposition}\label{pr:compact-multiplication} 
For $(m+\ell-k)p>d-1$ where $0\leq k\leq \ell\leq m$, pointwise multiplication
\[
(f,g)\mapsto fg \ \ \hbox{defines a continuous map}  \ \ 
H^{m,p}(S^{d-1}) \times H^{\ell,p}(S^{d-1}) \to H^{k,p}(S^{d-1}).
\]
In fact,  there is a constant $C=C(d,m,\ell,k,p)$ such that
\[
\| fg \|_{H^{k,p}} \leq 
C\, \|f\|_{H^{m,p}}\,\|g\|_{H^{\ell,p}}
\quad\hbox{ for all $f\in H^{m,p}$ and
$g\in H^{\ell,p}$.}
\]
\end{proposition}
 
\section{Asymptotic Spaces of Functions on $\RR^d$}

We want to consider functions $u\in H^{m,p}_{\loc}(\RR^d)$ which are bounded on $\RR^n$ and admit a partial asymptotic expansion as $|x|\to\infty$. To describe this partial asymptotic expansion,
let $\chi(t)$ be a smooth function satisfying
$\chi(t)= 0 $ for $ t \leq 1$, $\chi(t)= 1$ for $t\geq 2$, and $|\chi^{(k)}(t)|\leq M$ for $0\leq k\leq m$ and all $t$.
For a nonegative integer $N$, the functions that we consider are of the following form:  
\begin{subequations} \label{asymptoticexpansion}
\begin{equation}  \label{asymptoticexpansion1} u(x)=a(x)+f(x),
\quad\hbox{where }\end{equation}
\begin{equation}  \label{asymptoticexpansion2} a(x)=\chi(r)\left(a_0(\theta)+\frac{a_1(\theta)}{r}+\cdots + \frac{a_N(\theta)}{r^N}\right) \ \hbox{with}\ 
a_k\in H^{m+1+N-k,p}(S^{d-1}), \ \hbox{and}
\end{equation}
\begin{equation}  \label{asymptoticexpansion3} f(x)=o\left(|x|^{-N}\right) \quad \hbox{as}\ |x|\to\infty. \end{equation}
\end{subequations}
In (\ref{asymptoticexpansion2}) and throughout this paper, we use $r=|x|$ and $\theta=x/|x|\in S^{d-1}$.

We refer to $a$ in (\ref{asymptoticexpansion2}) as the {\it asymptotic function},  the function $a_k$ on $S^{d-1}$ as the {\it asymptotic of order $k$}, and $f$ as the {\it remainder function} for $u$.
We want to achieve (\ref{asymptoticexpansion3}) by requiring the remainder function $f$ to belong to one of the weighted Sobolev spaces discussed in the previous section. Let us begin with $W^{m,p}_\delta$. 
From Lemma \ref{le:W-properties}(d) we see that $f\in W_{\delta}^{m,p}(\RR^d)$ satisfies (\ref{asymptoticexpansion3}) provided $mp>d$ and $\delta+\frac{d}{p}\geq N$. However, for reasons that will become clear in the next section, we want to avoid  values of $\delta$ for which $\delta+\frac{d}{p}$ is an integer.  Consequently,  let us define
\begin{equation}\label{def:gamma}
\gamma_N=N+\gamma_0, \ \hbox{where $\gamma_0$ has been chosen to satisfy}\ 0<\gamma_0+\frac{d}{p}<1\, ,
\end{equation}
so that $f\in W_{\gamma_N}^{m,p}(\RR^d)$ with $mp>d$ satisfies $f(x)=o(|x|^{-N-\e})$ where $\e\in (0,1)$. Now we define
\begin{equation}\label{def:AW}
{\mathcal A}_{N}^{m,p}(\RR^d):=\left\{u \hbox{ is in the form (\ref{asymptoticexpansion}) where } 
f\in W_{\gamma_N}^{m,p}(\RR^d)\right\}.
\end{equation}
When the domain $\RR^d$ is understood, we simply write  ${\mathcal A}_{N}^{m,p}$ instead of ${\mathcal A}_{N}^{m,p}(\RR^d)$; when $N$ is fixed or understood, we may simply write $\gamma$ instead of $\gamma_N$.
The norm on ${\mathcal A}_{N}^{m,p}$ is given by
\begin{equation}\label{def:AW-norm}
\|u\|_{{\mathcal A}_{N}^{m,p}}=\|a_0\|_{H^{m+1+N,p}(S^{d-1})}+\cdots +\|a_N\|_{H^{m+1,p}(S^{d-1})}+\|f\|_{W_{\gamma}^{m,p}(\RR^d)}.
\end{equation}
This norm is complete, so ${\mathcal A}_{N}^{m,p}$ is a Banach space.  For an integer $n$ with $0\leq n\leq N$, we define closed subspaces 
\begin{equation}\label{def:AW-subspaces}
{\mathcal A}_{n,N}^{m,p}={\mathcal A}_{n,N}^{m,p}(\RR^d)=\left\{u\in {\mathcal A}_{N}^{m,p}: a_{k}=0 \ \hbox{ for $k<n$} \right\}.
\end{equation}

\begin{remark}\label{re:asymptotic_regularity}  
That the regularity of the asymptotic $a_k$ depends on $k$, i.e.\ $a_k\!\in \! H^{m+1+N-k,p}(S^{d-1})$, is an important feature of (\ref{asymptoticexpansion}); it will prove essential many times in the analysis below. It is also important
that the asymptotic $a_N$ has greater regularity than the remainder $f$, so we have assumed $a_N\in H^{m+1,p}$;  it may be possible to weaken this assumption by using fractional-order Sobolev spaces on $S^{d-1}$, but we have avoided this for simplicity.
\end{remark}

\begin{remark}\label{re:R-equivalentnorms}  
In the definition (\ref{asymptoticexpansion}), the specification that $\chi(r)\equiv 1$ for $r>2$ is somewhat arbitrary. In fact, if we introduce $\chi_R(t)=\chi(R^{-1}t)$, then $\chi_R(r)\equiv 1$ for $r>2R$ and we can write
\[
u=\chi\sum_{k=0}^N \frac{a_k(\theta)}{r^k}  + f = \chi_R\sum_{k=0}^N \frac{a_k(\theta)}{r^k} + \widetilde f,
\]
where $\widetilde f$ differs from $f$ by a  function with compact support:
\[
\widetilde f = f+(\chi-\chi_R)\sum_{k=0}^N \frac{a_k(\theta)}{r^k} .
\]
But we can estimate 
\[
\left\| (\chi-\chi_R)\sum_{k=0}^N \frac{a_k(\theta)}{r^k} \right\|_{W^{m,p}_{\gamma_N}} \leq 
C\, \sum_{k=0}^N \|a_k \|_{H^{m,p}}
\leq C\,\sum_{k=0}^N \|a_k\|_{H^{m+1+N-k,p}},
\]
where $C$ depends on $R$, $\chi$, $m$, $p$, $d$, and $N$, to conclude $\|\widetilde f \|_{W^{m,p}_{\gamma_N}} \leq 
C\|u\|_{{\mathcal A}_{n,N}^{m,p}}$.
Similarly, we can estimate $f$ in terms of the $a_k$ and $\widetilde f$, so if we use $\chi_R$
 in place of $\chi$ in  (\ref{asymptoticexpansion}), we will get a norm on the  Banach space  ${{\mathcal A}_{N}^{m,p}}$ that is equivalent to (\ref{def:AW-norm}). This will be important in subsequent sections.
In fact, it is sometimes convenient to consider the restriction of $u=a+f$ to the exterior domain $B_R^c=\RR^d\backslash B_R$. This generates a Banach space
${\mathcal A}_N^{m,p}(B_R^c)$ with norm 
\[
\|u\|_{{\mathcal A}_N^{m,p}(B_R^c)}:=\sum_{k=0}^{ N} \|a_k\|_{H^{m+1+N-k,p}(S^{d-1})}+ \|f\|_{W_{\gamma_N}^{m,p}(B_R^c)}.
\]
Notice that $ \|u\|_{{\mathcal A}_N^{m,p}(\RR^d)}$ is equivalent to $\|u\|_{{\mathcal A}_N^{m,p}(B_R^c)}+\|f\|_{H^{m,p}(B_R)}$.
\end{remark}

Now let us use  $H_\delta^{m,p}$ as the remainder space. From Lemma \ref{le:H-properties}(d) we see that $f\in H_{\delta}^{m,p}(\RR^d)$ satisfies (\ref{asymptoticexpansion3}) provided $mp>d$ and $\delta\geq N$. This suggests that we use $H_N^{m,p}$ as the remainder space. However, if we use (\ref{asymptoticexpansion2})  then we would want $\chi(r)\,r^{-N-1}\in H_N^{m,p}(B_1^c)$, but this is only true if $d<p$. Consequently, for given nonegative integer $N$, let $N^*$ be the positive integer satisfying
\begin{subequations} \label{H-asymptoticexpansion}
\begin{equation} \label{H-asymptoticexpansion1}
N-1+\frac{d}{p}< N^*\leq N+\frac{d}{p}
\end{equation}
and replace (\ref{asymptoticexpansion2}) with
\begin{equation} \label{H-asymptoticexpansion2}
a(x)=\chi(|x|)\sum_{k=0}^{ N^*}\frac{a_k(\theta)}{r^k},
\quad\hbox{with}\ a_k\in H^{m+1+N^*-k,p}(S^{d-1}).
\end{equation}
\end{subequations}
We always have $N^*\geq N$, but we have $N^*=N$ when $d=1$, or more generally if $d<p$. In any case,  let us define
\begin{equation}\label{def:AH}
{\AH}_{N}^{m,p}(\RR^d):=\left\{u \hbox{ is in the form (\ref{asymptoticexpansion1}) where $a$ satisfies  (\ref{H-asymptoticexpansion2}) and } 
f\in H_{N}^{m,p}(\RR^d)\right\}.
\end{equation}
When the domain $\RR^d$ is understood, we simply write  ${\AH}_{N}^{m,p}$.  We replace  (\ref{def:AW-norm}) by
\begin{equation}\label{def:AH-norm}
\|u\|_{{\AH}_{N}^{m,p}}=
\sum_{k=0}^{ N^*} \|a_k\|_{H^{m+1+N^*-k}(S^{d-1})} + \| f\|_{H_N^{m,p}(\RR^d)}.
\end{equation}
Under this norm, ${\AH}_{N}^{m,p}$ is a Banach space, and we define closed subspaces ${\AH}_{n,N}^{m,p}$ by requiring $a_k=0$ for $k<n$. Of course, Remarks \ref{re:asymptotic_regularity} and
\ref{re:R-equivalentnorms}  apply as well to (\ref{def:AH}) and (\ref{def:AH-norm}).

We next investigate some of the  properties of these asymptotic spaces. We begin with an elementary result.
\begin{lemma}\label{le:asymptoticsinW} 
If $a\in H^{m,p}(S^{d-1})$  then $a(\theta)\,r^{-k} \in W^{m,p}_{\delta}(B_1^c)\subset  H^{m,p}_{\delta}(B_1^c)$ for all $\delta<k-\frac{d}{p}$.
\end{lemma}
\noindent
The lemma is  easy to prove using integration in spherical coordinates and the simple computation:
\begin{equation} \label{D^k(a(theta))}
\begin{aligned}
a\in H^{m,p}(S^{d-1}), \ 0\leq |\beta|\leq m \ & \Rightarrow \ 
D^\beta (a(\theta)r^{-k}) = b_\beta(\theta)\,r^{-k-|\beta|} \ \hbox{for } |x|=r>1, \\
&\quad\hbox{where }  b_\beta\in H^{m-|\beta|,p}(S^{d-1}).
\end{aligned}
\end{equation}
\noindent
We will use Lemma \ref{le:asymptoticsinW} in confirming that our asymptotic spaces have the following properties:
\begin{proposition}\label{pr:A-properties}  
\begin{enumerate}
\item[(a)] If $n_1\geq n$ and $N_1\geq N$, then ${\mathcal A}_{n_1,N_1}^{m,p}\subset{\mathcal A}_{n,N}^{m,p}$ and  ${\AH}_{n_1,N_1}^{m,p}\subset {\AH}_{n,N}^{m,p}$.
\item[(b)]  If $m\geq 1$,  then $u\mapsto \partial u/\partial x_j$ is continuous  ${\mathcal A}_{n,N}^{m,p}\to{\mathcal A}_{n+1,N+1}^{m-1,p}$ and
 ${\AH}_{n,N}^{m,p}\to {\AH}_{n+1,N}^{m-1,p}$.
 \item[(c)] Multiplication by $\chi(r)\,r^{-k}$ is bounded ${\mathcal A}_{n,N}^{m,p}\to{\mathcal A}_{n+k,N+k}^{m,p}$ and
 ${\AH}_{n,N}^{m,p}\to {\AH}_{n+k,N+k}^{m,p}$.
 \item[(d)] Assume $m>d/p$. If $u\in {\mathcal A}_{n,N}^{m,p}$, then
 \begin{equation}\label{AW-derivatives}
 \sup_{x\in\RR^d} \x^{n+|\alpha|}\,|D^\alpha u(x)|\leq C\, \|u\|_{{\mathcal A}_{n,N}^{m,p}} \quad\hbox{for all}\ |\alpha|<m-d/p.
 \end{equation}
 If $u\in {\AH}_{n,N}^{m,p}$, then
 \begin{equation}\label{AH-derivatives}
 \sup_{x\in\RR^d} \x^n \, |D^\alpha u(x)|\leq C\,\|u\|_{{\AH}_{n,N}^{m,p}} \quad\hbox{for all}\ |\alpha|<m-d/p.
 \end{equation}
\end{enumerate}
\end{proposition}
\noindent

\smallskip\noindent{\bf Proof.} (a) Write $u\in {\mathcal A}_{n_1,N_1}^{m,p}$ as
\[
\begin{aligned}
u(x)&=\chi\,\left(\frac{a_{n_1}(\theta)}{r^{n_1}}+\cdots \frac{a_{N_1}(\theta)}{r^{N_1}}\right)+f_{N_1}
\quad\hbox{with}\ a_k\in H^{m+1+N_1-k}(S^{d-1}),\ f_{N_1}\in W^{m,p}_{\gamma_{N_1}}\\
&= \chi\,\left(\frac{a_{n}(\theta)}{r^{n}}+\cdots \frac{a_{N}(\theta)}{r^{N}}\right)+g_{N}
\quad\hbox{where}\ a_k=0 \quad\hbox{for $n\leq k <n_1$ and} \\
& \qquad\qquad g_N=
\begin{cases}
f_{N_1} & {\rm if } \ N_1=N \\
\chi\,\left(\frac{a_{N+1}}{r^{N+1}}+\cdots+\frac{a_{N_1}}{r^{N_1}}\right)+f_{N_1} & {\rm if } \ N_1>N.
\end{cases}
\end{aligned}
\] 
We clearly have $a_k\in H^{m+1+N-k}(S^{d-1})$ and (using Lemma \ref{le:asymptoticsinW})  $g_N\in W^{m,p}_{\gamma_N}$, so $u\in  {\mathcal A}_{n,N}^{m,p}$. Similarly for ${\AH}^{m,p}_{n_1,N_1}\subset {\AH}_{n,N}^{m,p}$.
 
To prove (b) let us first consider $u\in {\mathcal A}_{n,N}^{m,p}$ with asymptotics  $a_k\in H^{m+1+N-k,p}(S^{d-1})$
for $k=n,\dots,N$. We
 use (\ref{D^k(a(theta))}) with $|\beta|=1$  to compute
\begin{equation}\label{eq:d_j(a_k/r^k)}
\frac{\partial}{\partial x_j}\left( \chi(r)\frac{a_k(\theta)}{r^k}\right)=\chi(r)\, \frac{b_{k,j}(\theta)}{r^{k+1}}+\chi'(r)\, \theta_j \,\frac{a_k(\theta)}{r^k},
\end{equation}
where $b_{k,j}\in H^{m+N-k,p}(S^{d-1})=H^{(m-1)+1+(N+1)-(k+1),p}(S^{d-1})$. The term $b_{k,j}\,r^{-k-1}$ for $k=n,\dots,N$
is of the form of an asymptotic function in $\A^{m-1,p}_{n+1,N+1}$ while $\chi'(r)\, \theta_j \, a_k(\theta)\,r^{-k}$ has compact support so certainly belongs to the remainder space 
$W_{N+1}^{m-1,p}$.
Since  we also know that $\partial_j:W_{\gamma}^{m,p}\to W_{\gamma+1}^{m-1,p}$, we have
$\partial_j:{\mathcal A}_{n,N}^{m,p}\to{\mathcal A}_{n+1,N+1}^{m-1,p}$ is bounded. 

Next consider  $u\in{\AH}_{n,N}^{m,p}$, with asymptotics $a_k\in H^{m+1+N^*-k,p}(S^{d-1})$ for $k=n,\dots,N^*$.
Using \eqref{eq:d_j(a_k/r^k)} where $b_{k,j}\in H^{m+N^*-k}(S^{d-1})\subset H^{m+N^*-k-1}(S^{d-1})$ for $k=n,\dots,N^*-1$ and 
$\chi(r) \,r^{-N^*-1}\,b_{N^*,j}\in H_{N}^{m-1,p}(\RR^d)$ by Lemma \ref{le:asymptoticsinW}, together with the fact that 
$\partial_j:{H}_{n,N}^{m,p}(\RR^d)\to{H}_{n+1,N}^{m-1,p}(\RR^d)$ is bounded, we see that
$\partial_j:{\AH}_{n,N}^{m,p}\to{\AH}_{n+1,N}^{m-1,p}$ is bounded.

The proof of (c) is immediate.

To prove (d), first consider $u=\chi(r)\, a_k(\theta)\, r^{-k}$ with $n\leq k\leq N$ and $a_k\in H^{m+1+N-k}(S^{d-1})$. We generalize \eqref{eq:d_j(a_k/r^k)}  to conclude
\begin{subequations}
\begin{equation}
D^\alpha\left(\chi(r)\frac{a_k(\theta)}{r^k}\right)
=\chi(r)\, D^\alpha\left(\frac{a_k(\theta)} {r^k}\right) + g(x)
=\chi(r)\,  \frac{b_{k,\alpha}(\theta)}{ r^{k+|\alpha|}}+g(x),
\end{equation}
where $b_{k,\alpha}\in H^{m+1+N-k-|\alpha|}(S^{d-1})$ with 
\begin{equation}
\|b_{k,\alpha}\|_{H^{m+1+N-k-|\alpha|}(S^{d-1})}\leq c\,\|a_k\|_{H^{m+1+N-k}(S^{d-1})}
\end{equation}
and $g\in H^{m+2+N-k-|\alpha|}(\RR^d)$ has support in the annulus $A=\{x:1<|x|<2\}$ with
\begin{equation}
\|g\|_{H^{m+2+N-k-|\alpha|}(A)}\leq c\,\|a_k\|_{H^{m+1+N-k}(S^{d-1})}.
\end{equation}
\end{subequations}
Now $|\alpha|<m-d/p$ certainly implies $m+1+N-k-|\alpha|>(d-1)/p$, so we can use the Sobolev embedding theorem on $S^{d-1}$ to conclude
\begin{subequations}
\begin{equation}
\sup_{\theta\in S^{d-1}}|b_{k,\alpha}(\theta)|\leq \, C\,\|b_{k,\alpha}\|_{H^{m+1+N-k-|\alpha|}(S^{d-1})}
\end{equation}
and we can apply the Sobolev embedding theorem on $A$ to conclude
\begin{equation}
\sup_{1<|x|<2}|g(x)|\leq \, C\,\|g\|_{H^{m+2+N-k-|\alpha|}(A)}.
\end{equation}
\end{subequations}
Combining these inequalities, we have
\begin{equation}
\sup_{|x|>1} |x|^{n+|\alpha|} \,\left|D^\alpha\left(\chi(r)\frac{a_k(\theta)}{r^k}\right)\right| \leq C\,\|a_k\|_{H^{m+1+N-k}(S^{d-1})}.
\end{equation}
Thus, for an asymptotic function $a$ as in (\ref{asymptoticexpansion2}), we have
\begin{equation}
\sup_{|x|>1} |x|^{n+|\alpha|} \left| D^\alpha\,a(x)\,\right| \leq 
C\,\left(\|a_n\|_{H^{m+1+N-n}(S^{d-1})}+\cdots +\|a_N\|_{H^{m+1}(S^{d-1})}\right),
\end{equation}
and the same holds for the asymptotic function in \eqref{H-asymptoticexpansion2} provided we replace $N$ by $N^*$.
Now let us consider the remainder function $f$.
If $f\in W^{m,p}_{\gamma_N}(\RR^d)$, then we use Lemma \ref{le:W-properties} (d) with $\delta=N$ to conclude
\[
\sup_{x\in \RR^d} \x^{N+|\alpha|}\,\left|D^\alpha f(x)\right|  \leq C\, \|f\|_{W^{m,p}_{\gamma_N}}
\quad\hbox{for all}\ |\alpha|<m-d/p.
\]
On the other hand, if $f\in H^{m,p}_{\gamma_N}(\RR^d)$,
then we use Lemma \ref{le:H-properties} (d) with $\delta=N$ to conclude
\[
\sup_{x\in \RR^d} \x^{N}\,\left|D^\alpha f(x)\right| \leq C\,\|f\|_{H_{n}^{m,p}} \leq C\, \|f\|_{H_{N}^{m,p}}
\quad\hbox{for all}\ |\alpha|<m-d/p.
\]
Since $n\leq N$, these estimates imply  (\ref{AW-derivatives}) and (\ref{AH-derivatives}). $\hfill \Box$

\medskip
What about multiplication? The product of two partial asymptotic expansions involves a number of terms. The product of the remainder functions is covered by Propositions \ref{pr:H-multiplication} and \ref{pr:W-multiplication}; for convenience, we record here the following special case of those results:

\begin{lemma}\label{le:fg-product} 
 Assume $m>d/p$ and $k=0,\dots,m$.
\begin{enumerate}
\item[(a)] $\|fg\|_{H^{k,p}_{\delta_1+\delta_2}}\leq \|f\|_{H^{k,p}_{\delta_1}}\,\|g\|_{H^{m,p}_{\delta_2}}$
for $f\in H^{k,p}_{\delta_1}$ and $g\in H^{m,p}_{\delta_2}$.
\item[(b)]  $\|fg\|_{W^{k,p}_{\delta_1+\delta_2-d/p}}\leq \|f\|_{W^{k,p}_{\delta_1-d/p}}\,\|g\|_{W^{m,p}_{\delta_2-d/p}}$
for $f\in W^{k,p}_{\delta_1-d/p}$ and $g\in W^{m,p}_{\delta_2-d/p}$.
\end{enumerate}
\end{lemma}

\noindent
The product of an asymptotic term like $a_k(\theta)/|x|^k$ and a remainder function is covered by the following (in which we  use Lemma \ref{le:asymptoticsinW}(c) to assume $k=0$).

\begin{lemma}\label{le:af-product} 
Assume $a\in H^{s,p}(S^{d-1})$ for an integer $s>(d-1)/p$, $m$ is an integer  $0\leq m\leq s$, and $\delta\in \RR$. 
\begin{enumerate}
\item[(a)]  $ f\in W_\delta^{m,p}(B_1^c) \Rightarrow af\in W_\delta^{m,p}(B_1^c)$ and
$
\|af\|_{W_\delta^{m,p}(B_1^c)}\leq C\, \|a\|_{H^{s,p}(S^{d-1})}\|f\|_{W_\delta^{m,p}(B_1^c)}.
$
\item[(b)]  $f\in H^{m,p}_{\delta}(B_1^c) \Rightarrow a\,f \in  H^{m,p}_{\delta}(B_1^c)$ and
$
\|af\|_{H_\delta^{m,p}(B_1^c)}\leq C\, \|a\|_{H^{s,p}(S^{d-1})}\|f\|_{H_\delta^{m,p}(B_1^c)}.
$
\end{enumerate}
\end{lemma}

\noindent{\bf Proof.} For a nonnegative integer $\ell$, we simply denote by $D^\ell f$ a partial derivative of $f$ of order $\ell$. To show (a), we  want to estimate
 $\int_{|x|>1} |x|^{(\delta+\ell)p}|D^\ell(af)|^p\,dx$ for $\ell=0,\dots,m$. 
But $D^\ell (af)$ is a sum of products of the form
$D^i a\, D^j f$ where $i+j=\ell$. For $i=0$, we have
\[
\int_{|x|>1} \! |x|^{(\delta+\ell)p}|a\,D^\ell f|^p dx\leq \sup_{\theta\in S^{d-1}} |a(\theta)|^p \int_{|x|>1}\! |x|^{(\delta+\ell)p}|D^\ell f|^p dx
\leq C\,\|a\|^p_{H^{s,p}(S^{d-1})}\, \|f\|^p_{W_\delta^{m,p}(B_1^c)}.
\]
For $i>0$, $D^i a$ is a sum of products of the form $r^{-i}c_k(\theta) D_\theta^k a$ where $c_k$ is a polynomial in $\theta$ and $k=1,\dots,i$. Thus we want to estimate
\[
\int_{|x|>1} |x|^{(\delta+j)p}|D_\theta^k a|^p |D^j f|^p \,dx \quad\hbox{for}\ k=1,\dots,i; \ i+j=\ell.
\]
For fixed $r>1$, let us denote by $f_{(r)}$ the function on $S^{d-1}$ defined by $f_{(r)}(\theta)=f(r\theta)$. Now let us
use Proposition \ref{pr:compact-multiplication} (and $s> (d-1)/p$) to estimate
\[
\begin{aligned}
\int_{S^{d-1}} |D_\theta^k a(\theta)|^p\, |D^j f(r\theta)|^p\, ds_\theta & = \| D_\theta^k a\, (D^j f)_{(r)} \|^p_{L^p(S^{d-1})} \\
& \leq C \, \| D^k_\theta a \|^p_{H^{s-k,p}(S^{d-1})} \, \| (D^j f)_{(r)} \|^p_{H^{k,p}(S^{d-1})} \\
& \leq C \, \|a\|^p_{H^{s,p}(S^{d-1})} \,  \| (D^j f)_{(r)} \|^p_{H^{k,p}(S^{d-1})}.
\end{aligned}
\]
By trace theory, $D^j f \in W_{\loc}^{m-j,p}(B_1^c)$ implies $(D^j f)_{(r)}\in H^{k}(S^{d-1})$ provided $k\leq m-j-\frac{1}{p}$;  this last condition only fails when $j=m$, which does not occur since we have assumed $i>0$. Moreover, for a function $g(x)$ we can use $\partial g/\partial\theta_i=r\partial g/\partial x_i$ to estimate any derivative $D_\theta^k g$ on $S^{d-1}$ by
$| D_\theta^k g(r\theta)| \leq r |\nabla g(x)| + \cdots + r^k | \nabla^k g(x)|$, where $|\nabla^k g(x)|$ denotes the sum of the absolute values of all $x$ derivatives of $g$ of order $k$.
Applying this to $g=D^j f$, we obtain 
\[
\int_{|x|>1} |x|^{(\delta+j)p} \|D^j f\|^p_{H^{k,p}(S^{d-1})}\,dx\leq C\,\|f\|^p_{W_\delta^{\ell,p}(B_1^c)}
\quad\hbox{for } j+k\leq \ell.
\]
Thus we have shown for $i+j=\ell$ and $k=1,\dots,i$ that
\[
\int_{|x|>1} |x|^{(\delta+j)p}|D_\theta^k a|^p |D^j f|^p \,dx  \leq C \, \|a\|^p_{H^{s,p}(S^{d-1})} \, \|f\|^p_{W_\delta^{\ell,p}(B_1^c)}
\leq C \, \|a\|^p_{H^{s,p}(S^{d-1})} \, \|f\|^p_{W_\delta^{m,p}(B_1^c)}.
\]
Combining the cases $i=0$ and $i>0$, we have shown (a).

The proof of (b) follows the same outline as for (a). Again we write $D^\ell(af)$ as a sum of products $D^i a \,D^j f$ and treat the case $i=0$ by 
\[
\int_{|x|>1} |x|^{\delta p} |a D^\ell f|^p\,dx \leq \left(\sup_{\theta\in S^{d-1}}|a(\theta)|^p \right) \int_{|x|>1} |x|^{\delta p} |D^\ell f|^p\,dx
\leq  \|a \|^p_{H^{m,p}(S^{d-1})}\, \|f|^p_{H_\delta^{m,p}(B_1^c)}.
\]
For $i>0$, we want to estimate 
\[
\int_{|x|>1} |x|^{(\delta-i)p}|D_\theta^k a|^p |D^j f|^p \,dx \quad\hbox{for}\ k=1,\dots,i; \ i+j=\ell.
\]
But arguing as above and using $(\delta-i)p < \delta p$, we can conclude
\[
\int_{|x|>1} |x|^{(\delta-i)p}|D^\ell(af)|^p\,dx \leq C\, \|a\|_{H^{m,p}(S^{d-1})} \,  \|f\|_{H_\delta^{m,p}(B_1^c)} \quad\hbox{for } 
\ell \leq m. \quad\Box
\]

\medskip\noindent
Combining   Lemmas \ref{le:fg-product} and \ref{le:af-product}, we obtain the following

\begin{corollary}\label{co:WxA} 
If $m>d/p$, $0\leq n\leq N$, and $u\in {\mathcal A}_{n,N}^{m,p}$, then for any $f\in W_{\delta}^{k,p}$ where $0\leq k\leq m$ and $\delta\in\RR$ we have
\[
\|f\,u\|_{W_\delta^{k,p}}\leq C\,\|f\|_{W_\delta^{k,p}}\|u\|_{{\mathcal A}_{n,N}^{m,p}}.
\]
The analogous statement with $W$ replaced by $H$ and ${\mathcal A}$ replaced by ${\AH}$ is also true.
\end{corollary}

\medskip
We are now able to prove the following result on products for our asymptotic spaces:
\begin{proposition}\label{pr:A-products} 
For $m> d/p$ and $0\leq n_i\leq N_i$ for $i=1,2$,  let $n_0=n_1+n_2$ and $N_0=\min(N_1+n_2,N_2+n_1)$. Then
\begin{subequations}\label{est:multiplication}
\begin{equation}\label{est:AH-multiplication}
\|u\,v\|_{{\AH}_{\bar n,N_0}^{m,p} }\leq C\,\|u\|_{{\AH}_{n_1,N_1}^{m,p} }\|v\|_{{\AH}_{n_2,N_2}^{m,p} }
\quad\hbox{for}\ u\in {\AH}_{n_1,N_1}^{m,p},\ v\in {\AH}_{n_2,N_2}^{m,p},
\end{equation}
\begin{equation}\label{est:AW-multiplication}
\|u\,v\|_{{\mathcal A}_{\bar n,N_0}^{m,p} }\leq C\,\|u\|_{{\mathcal A}_{n_1,N_1}^{m,p} }\|v\|_{{\mathcal A}_{n_2,N_2}^{m,p} }
\quad\hbox{for}\ u\in {\mathcal A}_{n_1,N_1}^{m,p},\ v\in {\mathcal A}_{n_2,N_2}^{m,p}.
\end{equation}
\end{subequations}
\end{proposition}

\noindent{\bf Proof.} We shall prove (\ref{est:AW-multiplication}); the proof of (\ref{est:AH-multiplication}) is analogous. Since $p$ is fixed, we shall drop that notation, but let us introduce $\gamma_i=\gamma_{N_i}$ for $i=1,2$ and $\bar\gamma=\gamma_{\bar N}$. For $u\in {\mathcal A}_{n_1,N_1}^{m},\ v\in {\mathcal A}_{n_2,N_2}^{m}$, let us write
\[
u=\chi\left(\frac{a_{n_1}}{r^{n_1}}+\cdots+\frac{a_{N_1}}{r^{N_1}}\right)+f \quad\hbox{and}\quad 
v=\chi\left(\frac{b_{n_2}}{r^{n_2}}+\cdots+\frac{b_{N_2}}{r^{N_2}}\right)+g
\]
where $a_k\in H^{m+1+N_1-k}(S^{d-1})$, $b_k\in H^{m+1+N_2-k}(S^{d-1})$, $ f\in W_{\gamma_{1}}^{m}$, and $g\in W_{\gamma_{2}}^{m}.$
Taking the product, we can write
\begin{equation}\label{eq:product-uv}
\begin{aligned}
u\,v&=\chi^2\sum_{i=n_1}^{N_1}\sum_{j=n_2}^{N_2}\,\frac{a_ib_j}{r^{i+j}}+\chi\left(\frac{a_{n_1}}{r^{n_1}}+\cdots+\frac{a_{N_1}}{r^{N_1}}\right)g
+\chi\left(\frac{b_{n_2}}{r^{n_2}}+\cdots+\frac{b_{N_2}}{r^{N_2}}\right)f+fg \\
&=
\chi^2 \sum_{k=n_1+n_2}^{N_1+N_2}\frac{c_k}{r^k}+h\quad\hbox{where}\ c_k=\sum_{i+j=k} a_i b_j\ \hbox{and $h$ is all terms involving $f$ or $g$}.
\end{aligned}
\end{equation}
In order to show that $u\,v\in {\mathcal A}_{\bar n,\bar N}^{m}$ we need to show (i) $c_k\in H^{m+1+\bar N-k}(S^{d-1})$ and
(ii) $h\in W^{m}_{\bar\gamma}$. Of course, we also need to show 
that we can replace $\chi^2$ in (\ref{eq:product-uv}) by $\chi$. But $(\chi^2-\chi)$ is supported in $1<r<2$, so
\[
\begin{aligned}
\left\| (\chi^2-\chi) \sum_{k=\bar n}^{N_1+N_2}\frac{c_k}{r^k}\right\|_{W^{m}_{\gamma_{\bar N}}} & \leq
\sum _{k=\bar n}^{N_1+N_2}\|c_k\|_{H^{m}(S^{d-1})} \\
 \leq 
C\, &\sum_{i=n_1}^{N_1}\|a_i\|_{H^m(S^{d-1})}\sum_{j=n_2}^{N_2}\|b_j\|_{H^{m}(S^{d-1})}
\leq C\, \|u\|_{{\mathcal A}^m_{n_1,N_1}} \|v\|_{{\mathcal A}^m_{n_2,N_2}}.
\end{aligned}
\]

To prove $c_k\in H^{m+1+\bar N-k}(S^{d-1})$, we can use Proposition \ref{pr:compact-multiplication} to conclude
\[
\| a_i b_j \|_{H^{m+1+\bar N-k}}\leq 
C\,\|a_i\|_{H^{m+1+N_1-i}}\|b_j\|_{H^{m+1+N_2-j}} \quad\hbox{for}\ k=i+j,
\]
since the condition
\[
(m+1+N_1-i)+(m+1+N_2-j)-(m+1+\tilde N-k)>(d-1)/p
\]
reduces to just $m+1>(d-1)/p$, which is guaranteed by our assumption $m>d/p$. This also shows the desired estimate for 
 (\ref{est:AW-multiplication})
\[
\|c_k\|_{H^{m+1+\bar N-k}(S^{d-1})}\leq C\, \|u\|_{{\mathcal A}^m_{n_1,N_1}} \|v\|_{{\mathcal A}^m_{n_2,N_2}}.
\]

To show $h\in W^{m}_{\tilde\gamma}$ we have several terms to consider. Let us first consider $f\, g$.  
But  $f\in W^{m,p}_{\delta_1-\frac{d}{p}}$ and $g\in W^{m,p}_{\delta_2-\frac{d}{p}}$ for $N_i<\delta_i=\gamma_i+\frac{d}{p}<N_i+1$, so we can apply Proposition \ref{pr:W-multiplication} (using $m>d/p$) to conclude $fg\in W^m_{\delta_1+\delta_2-\frac{d}{p}}$. But $\bar N\leq N_1+N_2<\delta_1+\delta_2$, so we have $W^m_{\delta_1+\delta_2-\frac{d}{p}}\subset W^m_{\bar\gamma}$,
i.e.\
\[
\|f\,g\|_{W^{m}_{\bar\gamma}}\leq C\,\|f\|_{W^m_{\gamma_1}} \|g\|_{W^m_{\gamma_2}}
\]
As for the other terms, we can use $\bar\gamma\leq n_2+\gamma_1,n_1+\gamma_2$ and Lemmas \ref{le:W-properties} and \ref{le:af-product} to conclude
\[
\sum_{k=n_1}^{N_1}\left\|\frac{a_k}{r^k}\,g\right\|_{W^m_{\bar\gamma}(B_1^c)}
\leq C\,\sum_{k=n_1}^{N_1} \left\|\frac{a_k}{r^k}\,g\right\|_{W^m_{n_1+\gamma_2}(B_1^c)}
\leq C\, \left(\sum_{k=n_1}^{N_1} \|a_k\|_{H^m(S^{d-1})}\right)\|g\|_{W^m_{\gamma_2}}
\]
and
\[
\sum_{k=n_2}^{N_2}\left\|\frac{b_k}{r^k}\,f\right\|_{W^m_{\bar\gamma}(B_1^c)}
\leq C\,\sum_{k=n_2}^{N_2} \|b_k\,f\|_{W^m_{n_2+\gamma_1}(B_1^c)}
\leq C\, \left(\sum_{k=n_2}^{N_2} \|b_k\|_{H^m(S^{d-1})}\right)\|f\|_{W^m_{\gamma_1}}.
\]
Finally we use $\bar\gamma=\bar N+\gamma_0$ where $\gamma_0<1-\frac{d}{p}$ to show the term $\chi^2\sum_{k=\bar N+1}^{N_1+N_2} \,r^{-k}\, c_k$ is in the remainder space:
\[
\left\| \chi^2\sum_{k=\bar N+1}^{N_1+N_2} \frac{c_k}{r^k}\right\|_{W^m_{\bar\gamma}(B_1^c)}
\leq  \sum_{k=\bar N+1}^{N_1+N_2} \|c_k\|_{H^{m}(S^{d-1})} \|\chi^2\|_{W^m_{\gamma_0-1}}
\leq C\, \|u\|_{{\mathcal A}^m_{n_1,N_1}} \|v\|_{{\mathcal A}^m_{n_2,N_2}}.  \quad \hfill\Box
\] 

\medskip
As a special case of Proposition \ref{pr:A-products}, we obtain

\begin{corollary}\label{co:Banachalgebras} 
If $m>d/p$, and $0\leq n\leq N$, then ${\AH}_{n,N}^{m,p}$ and ${\mathcal A}_{n,N}^{m,p}$ are Banach algebras.
\end{corollary}

Since this paper is mostly concerned with diffeomorphisms of $\RR^d$, we need to consider asymptotic spaces of vector-valued functions. Here we use bold-face ${\bf u}$ for a vector-valued function and denote its components by $u^j$. Let us define the Banach spaces
\begin{subequations}\label{def:A-vector}
\begin{equation}\label{def:AH-vector}
{\AH}_{n,N}^{m,p}(\RR^d,\RR^d)=\left\{{\bf u}:\RR^d\to\RR^d \,|\, u^j\in {\AH}_{n,N}^{m,p}(\RR^d)
\right\},\quad \|{\bf u}\|_{{\AH}_{n,N}^{m,p}}=\sum_{j=1}^d \|u^j\|_{{\AH}_{n,N}^{m,p}}.
\end{equation}
and
\begin{equation}\label{def:AW-vector}
{\A}_{n,N}^{m,p}(\RR^d,\RR^d)=\left\{{\bf u}:\RR^d\to\RR^d \,|\, u^j\in {\mathcal A}_{n,N}^{m,p}(\RR^d)
\right\}, \quad \|{\bf u}\|_{{\mathcal A}_{n,N}^{m,p}}=\sum_{j=1}^d \|u^j\|_{{\mathcal A}_{n,N}^{m,p}}.
\end{equation}
\end{subequations}
As in the scalar-valued case, we will abbreviate ${\AH}_{0,N}^{m,p}$ simply as ${\AH}_N^{m,p}$ and suppress the notation $(\RR^d,\RR^d)$ when it is clear from the context that we are considering vector fields on ${\RR^d}$.

\section{Application to the Laplacian and Helmholtz Decompositions}

 The asymptotic spaces ${\A}_{n,N}^{m,p}$ are generally preferable to ${\AH}_{n,N}^{m,p}$ in applications involving the Laplacian $\Lap=\sum_{i=1}^d \partial^2/\partial x_i^2$. In this section we will illustrate this by discussing the mapping properties of $\Lap$  and an application to the Helmholtz decomposition of vector fields.
 
 To begin with, consider the mapping
 \begin{equation}\label{Delta:W->W}
\Delta: W^{m+1,p}_\delta(\RR^d)\to W^{m-1,p}_{\delta+2}(\RR^d) \quad\hbox{for}\ m\geq 1.
\end{equation}
Clearly, (\ref{Delta:W->W}) is continuous for all $\delta\in\RR$, and in \cite{Mc} it was shown that (\ref{Delta:W->W}) is injective for $\delta>-d/p$ and an isomorphism (in particular invertible) for $0<\delta+d/p<d-2$ (when $d\geq 3$).  For $N<\delta+d/p<N+1$, where $N$ is an integer $\geq d-2$, it was also shown in \cite{Mc} that (\ref{Delta:W->W}) is Fredholm with explicitly specified cokernel.  We  now observe that for arbitrary $g\in W^{m-1,p}_{\delta+2}(\RR^d)$ for $N<\delta+d/p<N+1$, we can find $u\in {\mathcal A}_{d-2,N}^{m+1,p}$ such that $\Delta u=g$; here we have used $\gamma_N=\delta$ in defining ${\mathcal A}_{d-2,N}^{m+1,p}$, so our hypothesis on $g$ can instead be written $g\in  W^{m-1,p}_{\gamma_N+2}(\RR^d)$.

\begin{lemma} \label{le:inverseLaplacian}
Suppose $d\geq 2$ and $m\geq 1$.
\begin{enumerate}
\item[(a)] For $d\geq 3$,  there is a bounded operator
\begin{subequations} \label{inverseLaplacian}
\begin{equation}  \label{inverseLaplacian1}
K:W_{\gamma_N+2}^{m-1,p}(\RR^d)\to {\mathcal A}_{d-2,N}^{m+1,p}(\RR^d),
\end{equation}
 such that $\Delta Kg=g$.
In other words,  $u=Kg$ is of the form
\begin{equation}  \label{inverseLaplacian2}
u(x)=\chi(|x|)\left( \frac{a_{d-2}}{r^{d-2}}+\cdots+\frac{a_N(\theta)}{r^N}\right)+f(x)
\end{equation}
\end{subequations}
where each $a_k(\theta)/r^k$ is harmonic for $x\not=0$ (so $a_k\in C^\infty(S^{d-1})$) and $f\in W_{\gamma_N}^{m,p}(\RR^d)$.
\item[(b)] For $d=2$, the result also holds, except $1/r^{2-d}$ in (\ref{inverseLaplacian2}) is replaced by $\log r$. Of course, this means that the asymptotic space $ {\mathcal A}_{0,N}^{m+1,p}$ in (\ref{inverseLaplacian1}) must be replaced by 
\begin{equation}
\begin{aligned}
\A^{m+1,p}_{0^*,N}(\RR^d)&=\left\{u=\chi\left(a_0^*\log r+a_0(\theta)+\cdots +\frac{a_N(\theta)}{r^N}\right)+f: \right. \\
&\left. a_0^*=const, \ a_k\in H^{m+2+N-k}(S^{d-1}),\ f\in W_{\gamma_N}^{m+1,p}(\RR^d)\right\}.
\end{aligned}
\end{equation}
\end{enumerate}
\end{lemma}

\noindent
{\bf Proof.} Let $\Gamma(|x|)$ denote the fundamental solution for the Laplace operator in $\RR^d$ and $K=\Gamma\star$ denote the convolution operator. As shown in \cite{Mc},  $K:W_{\delta+2}^{m-1,p}(\RR^d)\to W_{\delta}^{m+1,p}(\RR^d)$ is an isomorphism for $0<\delta+d/p<d-2$ (when $d\geq 3$); since we know $\gamma_N$ satisfies $N<\gamma_N+d/p<N+1$, we conclude that  (\ref{inverseLaplacian1}) is an isomorphism for $N\leq d-3$ (and ${\mathcal A}_{d-2,N}^{m+1,p}=W^{m+1,p}_{\gamma_N}$). 
For $N\geq d-2$, $K:W_{\gamma_N+2}^{m-1,p}(\RR^d)\to W_{\gamma_N}^{m+1,p}(\RR^d)$ is no longer bounded, and we either need to restrict the domain space or expand the range space. Let us first describe what happens for $N=d-2$ and then consider the general case.

For $d-2<\delta+d/p<d-1$, it was shown in \cite{Mc} that (\ref{Delta:W->W}) is injective with constants as cokernel: if we let
$\widetilde W^{m-1,p}_{\delta+2}(\RR^d)=\{g\in W^{m-1,p}_{\delta+2}(\RR^d): \int_{\RR^d}g(x)dx=0\}$, then 
$K: \widetilde W^{m-1,p}_{\delta+2}(\RR^d) \to W^{m+1,p}_\delta(\RR^d)$ is bounded. Taking $\gamma_N=\delta$,
we have $K: \widetilde W^{m-1,p}_{\gamma_N+2}(\RR^d) \to W^{m+1,p}_\gamma(\RR^d)$ is bounded. 
To extend $K$ to general 
$g\in W^{m-1,p}_{\gamma_N+2}(\RR^d)$, let us observe that $\Delta\left(\chi(r)\,\Gamma(r)\right)$ has compact support and we use Green's first identity to calculate
\[
\int_{\RR^d}\Delta(\chi(r)\,\Gamma(r))\,dx=\int_{|x|<2}\Delta(\chi(r)\,\Gamma(r))\,dx=\int_{|x|=2} \frac{\del}{\del r} \Gamma(r)\,ds=1.
\]
Now we define 
\begin{equation}
\widetilde g(x) = g(x)-c_0\, \Delta\left(\chi(r)\,\Gamma(r)\right) \quad\hbox{where}\ c_0=\int_{\RR^d} g\,dx.
\end{equation}
Notice that $c_0$ is finite; in fact, using $N=d-2$ and H\"older's inequality, we  easily confirm that
\begin{equation}
|c_0|\leq \|g\|_{L^1}\leq C\,\|g\|_{L^p_{\gamma_N+2}}.
\end{equation}
Then $\int \widetilde g\,dx=0$, so $\widetilde g\in \widetilde W^{m-1,p}_{\gamma_N+2}(\RR^d)$ and we can let $f=K\widetilde g$ to find $f\in W^{m+1,p}_{\gamma_N}(\RR^d)$. Finally, we define $Kg$ by
\begin{equation}
Kg=f+c_0\,\chi(r)\Gamma(r).
\end{equation}
For $d\geq 3$, $u=Kg$ is of the form (\ref{inverseLaplacian2}) for $N=d-2$ and satisfies $\Delta u=g$ as well as the estimate
\[
\|u\|_{{\mathcal A}_{d-2,d-2}^{m+1,p}}=|c_0|+\|f\|_{W_{\gamma_N}^{m+1,p}}\leq C\,\|g\|_{W_{\gamma_N+2}^{m+1,p}}.
\]
For $d=2$ we have $u=c_0\chi(r)\log r+f$, so it is clear how to treat this case as well. 
This proves the result for $N=d-2$.

More generally, for $k+d-2<\delta+d/p<k+d-1$ where $k=N-d+2>0$, it was shown in \cite{Mc} that (\ref{Delta:W->W}) is injective with cokernel equal to the harmonic polynomials of degree less than or equal to $k$. If we let ${\mathcal H}_k$ denote the spherical harmonics of degree $k$, let $N(k)=\hbox{dim}\,{\mathcal H}_k$, and choose an orthonormal basis $\{\phi_{k,j}:j=1,\dots,N(k)\}$ for ${\mathcal H}_k$, then
a basis for the space of harmonic polynomials that are homogeneous of degree $k$ is $\{ \phi_{k,j}(\theta)\,r^k: j=1,\dots,N(k)\}$.
Consequently, if we define
\[
\widetilde W^{m-1,p}_{\delta+2}(\RR^d)=\left\{ g\in W^{m-1,p}_{\delta+2}(\RR^d): \int_{\RR^d} g(x) \phi_{\ell,j}(\theta)\,r^\ell dx=0,
\ j=1,\dots,N(\ell),\ \ell=0,\dots,k\right\},
\]
then \cite{Mc} showed that $K: \widetilde W^{m-1,p}_{\delta+2}(\RR^d) \to W^{m+1,p}_\delta(\RR^d)$ is bounded.
Taking $\gamma_N=\delta$ and considering a
 general $g\in W^{m-1,p}_{\gamma_N+2}(\RR^d)$, we define 
\[
c_{\ell,j}=\int_{\RR^d} g(x) \phi_{\ell,j}(\theta)\,r^\ell dx \quad\hbox{for}\ j=1,\dots,N(k),\ \ell=1,\dots,k.
\]
Using H\"older's inequality, we can confirm
\[
|c_{\ell,j}|\leq C\,\|g\|_{L^p_{\gamma_N+2}},
\]
and in particular that $c_{\ell,j}$ is finite. Recall that 
\[
\frac{\phi_{\ell,j}(\theta)}{r^{d-2+\ell}} \quad\hbox{is harmonic for $r>0$,}
\]
so we can use Green's second identity to calculate
\[
\begin{aligned}
\int_{\RR^d} \phi_{\ell',j'}(\theta) r^{\ell'} \,&\Delta\left(\chi(r)\frac{\phi_{\ell,j}(\theta)}{r^{d-2+\ell}}\right)\,dx
= \int_{|x|\leq 2} \phi_{\ell',j'}(\theta) r^{\ell'} \,\Delta\left(\chi(r)\frac{\phi_{\ell,j}(\theta)}{r^{d-2+\ell}}\right)\,dx \\
&=\int_{|x|=2} \left(\phi_{\ell',j'}(\theta)\,r^{\ell'}\,\frac{\del}{\del r}\frac{\phi_{\ell,j}(\theta)}{r^{d-2+\ell}}
-\frac{\phi_{\ell,j}(\theta)}{r^{d-2+\ell}} \frac{\del}{\del r}\phi_{\ell',j'}(\theta)\,r^{\ell'} \right)\,ds \\
&=\begin{cases}
2-d-2\ell & \hbox{if}\ \ell'=\ell \hbox{ and } j'=j, \\ 0 & \hbox{otherwise.}
\end{cases}
\end{aligned}
\]
Define
\begin{equation}
\widetilde g(x)=g(x)-c_0\,\Delta(\chi(r)\Gamma(r)) - 
\sum_{\ell=1}^k \sum_{j=1}^{N(\ell)} \frac{c_{\ell,j}}{(2-d-2\ell)}\,\Delta\left(\chi(r)\frac{\phi_{\ell,j}(\theta)}{r^{d-2+\ell}}\right).
\end{equation}
Then $\widetilde g\in \widetilde W^{m-1,p}_{\gamma_N+2}(\RR^d)$ and we can let $f=K\widetilde g\in W_{\gamma_N}^{m+1,p}$. Finally, we
define 
\[
Kg=f+\chi(r)\left[c_0\,\Gamma(r)+\sum_{\ell=1}^k \sum_{j=1}^{N(\ell)} \frac{c_{\ell,j}}{(2-d-2\ell)}\frac{\phi_{\ell,j}}{r^{d-2+j}}\right].
\]
We see that $u=Kg$ is of the form (\ref{inverseLaplacian2}) and satisfies $\Delta u=g$ as well as the estimate
\[
\|u\|_{{\mathcal A}_{d-2,N}^{m+1,p}}\leq C\,\left(|c_0|+\sum_{\ell=1}^k \sum_{j=1}^{N(\ell)}|c_{\ell,j}|\right)+\|f\|_{W_\delta^{m+1,p}}\leq C\,\|g\|_{W_{\delta+2}^{m+1,p}},
\]
where $C$ depends on the Sobolev norms of $\phi_{\ell,j}$ on $S^{d-1}$, but not on $g$. 
This completes the proof.

 \hfill  $\Box$
 
 \medskip
 Notice that  \eqref{Delta:W->W} generalizes to
 \begin{equation}\label{Delta:A->A}
 \Lap: \A^{m+1,p}_N(\RR^d)\to \A^{m-1,p}_{2,N+2}(\RR^d),
 \end{equation}
 and we want to consider its invertibility.
For $v\in \A^{m-1,p}_{2,N+2}(\RR^d)$,  we write $v=b+g$ where $b=\chi(r^{-2}b_2+\cdots r^{-N-2}b_{N+2})$
with $b_{k+2}\in H^{m+N-k}(S^{d-1})$ and $g\in W^{m-1}_{\gamma_N+2}$. To define $\Lap^{-1}v$, we first try to find an asymptotic function $a=\chi(a_0+\cdots r^{-N}a_N)\in \A^{m+1,p}_N(\RR^d)$ satisfying
 \begin{equation}\label{eq:Delta(a)=b}
 \Lap\left(\frac{a_k(\theta)}{r^k}\right)=\frac{b_{k+2}(\theta)}{r^{k+2}} \quad\hbox{for $k=0,\dots,N$}.
 \end{equation}
 Then we will use Lemma \ref{le:inverseLaplacian} to find a remainder function $f$ so that
 $u=a+f\in\A^{m+1,p}_N(\RR^d)$ is an exact solution of $\Lap u=v$. 
 To solve \eqref{eq:Delta(a)=b},  we can use separation of variables. In fact, 
 using 
$
 \Lap=\partial_r^2+(d-1)\,r^{-1}\partial_r+r^{-2}\Lap_h,
$
 where $h$ is the induced metric on $S^{d-1}$, we find that $a_k$ must satisfy
 \begin{equation}\label{eq:Delta_h(ak)=}
 \Lap_h a_k-k(d-2-k)a_k=b_{k+2}\quad\hbox{on $S^{d-1}$}.
 \end{equation}
 If $k(d-2-k)>0$, then we can uniquely solve \eqref{eq:Delta_h(ak)=} to find $a_k$. However, 
 for $k=0$ or $k=d-2$, we have a simple solvability condition, namely $\int_{S^{d-1}} b_{k+2}\,ds=0$, and the solution $a_k$ is only unique up to an additive constant; this is expected since $c_0$ and $c_{d-2}r^{d-2}$ are harmonic for $r>0$.
Let us consider two closed subspaces of ${\A}^{m-1,p}_{2,N+2}(\RR^d)$:
  \begin{equation}\label{def:tilde(A)}
\widetilde{\A}^{m-1,p}_{2,N+2}(\RR^d)=
 \left\{u=\chi\left(\frac{b_2(\theta)}{r^2}+\cdots+\frac{b_{N+2}(\theta)}{r^{N+2}}\right)+f \in {\A}^{m-1,p}_{2.N+2}(\RR^d): 
\int_{S^{d-1}}b_{d}(\theta)\,ds=0\right\}.
 \end{equation}
 \begin{equation}\label{def:doubletilde(A)}
 \begin{aligned}
\accentset{\approx}{\A}^{m-1,p}_{2,N+2}(\RR^d)=
 \left\{u=\chi\left(\frac{b_2(\theta)}{r^2}+\cdots+\frac{b_{N+2}(\theta)}{r^{N+2}}\right)\right. &+f \in {\A}^{m-1,p}_{2.N+2}(\RR^d):  \\
& \left.
  \int_{S^{d-1}}b_2(\theta)\,ds=\int_{S^{d-1}}b_{d}(\theta)\,ds=0\right\}.
 \end{aligned}
 \end{equation}
 Of course, if $d>N+2$, then the solvability condition $\int b_d\,ds=0$ is vacuous.
 
 \begin{proposition} \label{pr:inverseLaplacian2}
 \begin{enumerate}
\item[(a)]  For $d\geq 3$, $m\geq 1$, and $0\leq N\leq d-2$, there is a bounded operator
\begin{subequations} \label{inverseLaplacian-A}
 \begin{equation}\label{inverseLaplacian-Aa}
 K:  \widetilde{\A}^{m-1,p}_{2,N+2}(\RR^d)\to  {\A}^{m+1,p}_{0^*,N}(\RR^d)
 \end{equation}
 satisfying $\Lap Kv=v$ for all $v\in  \widetilde{\A}^{m-1,p}_{2,N+2}(\RR^d)$. 
 This operator is also bounded
  \begin{equation}\label{inverseLaplacian-Ab}
 K:  \accentset{\approx}{\A}^{m-1,p}_{2,N+2}(\RR^d)\to  {\A}^{m+1,p}_{0,N}(\RR^d).
 \end{equation}
 \item[(b)] For $d=2$ and $m\geq 1$ the operator is bounded
  \begin{equation}\label{inverseLaplacian-Ac}
 K:  \widetilde{\A}^{m-1,p}_{2,2}(\RR^2)= \accentset{\approx}{\A}^{m-1,p}_{2,2}(\RR^2)\to  {\A}^{m+1,p}_{0^*,0}(\RR^2).
 \end{equation}
 \end{subequations}
 \end{enumerate}
 \end{proposition}
 
 \noindent
{\bf Proof. } As indicated above, for  $v=b+g=\chi(r^{-2}b_2+\cdots r^{-N-2}b_{N+2})+g\in \accentset{\approx}\A^{m-1,p}_{2,N+2}$ with $b_{k+2}\in H^{m+N-k,p}(S^{d-1})$, we have the necessary solvability conditions so that we can find $a_k\in H^{m+2+N-k,p}(S^{d-1})$ solving \eqref{eq:Delta_h(ak)=} with $\| a_k\|_{H^{m+2+N-k,p}(S^{d-1})}\leq C\, \| b_{k+2}\|_{H^{m+N-k,p}(S^{d-1})}$ for $k=0,\dots,N$; in fact, the $a_k$ are unique except for $k=0,d-2$. 
Of course, the same analysis applies for $d=2$. However, now let us assume $d>2$ and $v\in \widetilde{\A}^{m-1,p}_{2,N+2}(\RR^d)$ with $\int b_2(\theta)\,ds\not=0$. Then the necessary solvability condition does not hold in order to be able to solve
\eqref{eq:Delta_h(ak)=} for $k=0$. Instead, let us replace $a_0(\theta)$ by $a_0^*\log r+a_0(\theta)$ and instead of \eqref{eq:Delta(a)=b} try to solve 
\[
 \Lap (a_0^*\log r+a_0(\theta))=\frac{b_2(\theta)}{r^2}\, ,
\]
with $a_0^*$ being a constant. In place of \eqref{eq:Delta_h(ak)=} we have
$
\Lap_h a_0 + (d-2)a_0^* = b_2.
$
If we choose
\[
a_0^* = \frac{1}{d-2}\,\meanint_{S^{d-1}}b_2(\theta)\,ds,
\]
then $|a_0|\leq C\,\|b_2\|_{L^p(S^{d-1})}\leq C\|b_2\|_{H^{m+N,p}(S^{d-1})}$ 
and we can find $a_0\in H^{m+2+N,p}(S^{d-1})$ with $\| a_0\|_{H^{m+2+N,p}(S^{d-1})}\leq C\, \| b_{2}\|_{H^{m+N,p}(S^{d-1})}$. To summarize, we have defined $a_0^*,a_0(\om),\dots,a_N(\om)$ (where $a_0^*=0$ unless $d>2$ and $\int b_2\,ds\not= 0$) so that
\[
\Lap\left(a_0^*\log r+a_0(\theta)+\cdots +r^{-N}a_N(\theta)\right)=r^{-2}b_2(\theta)+\cdots r^{-N-2}b_{N+2}(\theta).
\]

Now let $u=a+f$ where
\[
a=\chi\left( a_0^*\,\log r+a_0(\theta)+\cdots+r^{-N}a_N(\theta)\right).
\]
We need to use Lemma \ref{le:inverseLaplacian} to find the remainder function $f$ so that $\Lap u=v$.
We compute
\[
\Lap a=\chi\,b+\Lap\chi\,(a_0^*\,\log r+a_0+\cdots+r^{-N}a_N)+\nabla\chi\cdot\nabla(a_0^*\,\log r+a_0+\cdots+r^{-N}a_N).
\]
So we want $f$ to satisfy
\[
\Lap f=h:= g-\Lap\chi(a_0^*\,\log r+a_0+\cdots+r^{-N}a_N)-\nabla\chi\cdot\nabla(a_0^*\,\log r+a_0+\cdots+r^{-N}a_N).
\]
But $h$ and $g$ differ by a function in $H^{m+1,p}(\RR^d)$ with compact support, and $g\in W^{m-1,p}_{\gamma_N+2}$, so $h\in W^{m-1,p}_{\gamma_N+2}$. Consequently, we can apply Lemma \ref{le:inverseLaplacian} to find $f=Kh\in {\A}^{m+1,p}_{d-2,N}$ (or 
$u\in {\A}^{m+1,p}_{0^*,N}$ if $d=2$) satisfying $\Lap w=h$. We see that $u=a+f$ satisfies $\Lap u=v$ and the mapping $K:v\mapsto u$ is continuous between the appropriate spaces. \hfill $\Box$

\begin{remark} 
The problem with extending this result to $N>d-2$ is that $\log r$ terms arise in the solution of \eqref{eq:Delta(a)=b} for large values of $k$. Cf.\ Example \ref{ex:d=N=2} below.
\end{remark}
 
 \medskip
Now we turn to the application to Helmholtz decompositions. It is well-known that a $C^1$-vector field ${\bf u}$ in ${\RR^3}$ satisfying $D^k{\bf u}=O(|x|^{-1-k-\e})$ as $|x|\to\infty$ for $k=0,1$ and some $\e>0$ can be decomposed into the sum of a  unique divergence-free vector field with the same decay property and a gradient field:
\begin{equation}\label{eq:Helmholtz}
{\bf u}={\bf v}+\nabla w, \quad \hbox{div}\,{\bf v}=0.
\end{equation}
Moreover, ${\bf v}$ and $\nabla w$ are orthogonal in that $\int {\bf v}\cdot\nabla w\,dx=0$. This is called the {\it Helmholtz decomposition in $\RR^3$}.
We now show that (\ref{eq:Helmholtz}) can be achieved when $d\geq 2$ and ${\bf u}\in {\mathcal A}_{1,N}^{m,p}$; this allows some vector fields ${\bf u}$ satisfying $O(|x|^{-1})$ as $|x|\to\infty$ instead of requiring  $O(|x|^{-1-\e})$. (While ${\bf v},\nabla w\in {\mathcal A}_{1,N}^{m,p}$, the orthogonality $\int {\bf v}\cdot\nabla w\,dx=0$ need not hold for $d\geq 2$ and ${\bf u}\in {\mathcal A}_{1,N}^{m,p}$.)
\begin{theorem}\label{th:Helmholtz} 
If $d\geq 2$, $m\geq 1$, and $1\leq N\leq d-1$, then every vector field ${\bf u}\in  {\mathcal A}_{1,N}^{m,p}$ can be written 
 in the form (\ref{eq:Helmholtz}) where ${\bf v}\in  {\mathcal A}_{1,N}^{m,p}$ is divergence-free and
$w\in H^{m+1,p}_\loc$ with $|w(x)|=o(|x|)$ as $|x|\to\infty$; in fact, {\bf v} is uniquely determined and $w$ is unique up to an additive constant. The map $P_0:{\bf u}\mapsto {\bf v}$ defines a bounded linear map ${\mathcal A}_{1,N}^{m,p}\to{\mathcal A}_{1,N}^{m,p} $ that is a projection: $P_0^2{\bf u}=P_0{\bf u}$.
\end{theorem}

\noindent
The operator $P_0$ is called the {\it Euler projector}. We can reformulate  Theorem \ref{th:Helmholtz} as a statement about closed subspaces.
\begin{corollary} Under the hypotheses of Theorem \ref{th:Helmholtz}, $ {\mathcal A}_{1,N}^{m,p}$ can be decomposed into a direct sum of closed subspaces
$
{\mathcal A}_{1,N}^{m,p}=\accentset{\circ}{\mathcal A} \oplus {\mathcal G},
$
where $\accentset{\circ}{\mathcal A}=\{{\bf u}\in {\mathcal A}_{1,N}^{m,p}: {\rm div}\,{\bf u}=0\}$ and ${\mathcal G}$ is the nullspace of the Euler projector $P_0$.
\end{corollary}

\noindent
{\bf Proof of Theorem \ref{th:Helmholtz}.}  When the 1st-order derivatives of {\bf u} are $O(|x|^{-2-\e})$ as $|x|\to\infty$, the decomposition (\ref{eq:Helmholtz}) can be found by letting $w=K \,\hbox{div}\,{\bf u}$, where $K$ is defined by convolution with the fundamental solution. For a vector field ${\bf u}\in {\mathcal A}_{1,N}^{m,p}$, we will replace $K$ by the operator discussed in Proposition \ref{pr:inverseLaplacian2}; however, we first need to use separation of variables to study $\hbox{div}\,{\bf u}$.

Let $\omega^1,\dots,\omega^{d-1}$ be local coordinates on $S^{d-1}$, considered as functions of Euclidean coordinates which, in this
proof, we index by superscripts, i.e.\ 
$x^1,\dots,x^d$. Let $h=h_{\alpha\beta}\,d\omega^\alpha\,d\omega^\beta$ denote the Riemannian metric on $S^{d-1}$ induced by the Euclidean metric $g=dx^2=dr^2+r^2 \,h$. Recall that the divergence of a vector field {\bf v} may be computed in general coordinates $\bar x^1,\dots,\bar x^d$  by 
\[
\hbox{div}\,{\bf v}=\frac{1}{\sqrt{{\rm det}\,\bar g_{ij}}} \,\frac{\partial}{\partial \bar x^j} \left(\sqrt{{\rm det}\,\bar g_{ij}}\, \bar{ v}^j \right),
\quad\hbox{where}\ g=\bar g_{ij}\,d\bar{x}^i\,d\bar{x}^j\ \hbox{and}\  \bar v^j=\frac{\partial \bar x^j}{\partial x^i}\,v^i.
\]
To compute  the divergence of {\bf v} in the coordinates $(\bar x^0,\dots,{\bar x}^{d-1})=(r,\omega^1,\dots,\omega^{d-1})$, we first compute its components in these coordinates by
\[
\bar v^0=\frac{\partial\, r}{\partial x^j} \,v^j=\theta_j\,v^j, \quad\hbox{where}\ \theta_j=\frac{x^j}{r}\in S^{d-1},
\]
\[
\bar v^\alpha = \frac{\partial\, \omega^\alpha}{\partial x^j}\,v^j=\omega_j^\alpha\,r^{-1}\,v^j, 
\quad\hbox{where}\ \omega_j^\alpha = r\,\frac{\partial\, \omega^\alpha}{\partial x^j}\in S^{d-1},
\]
and then compute the divergence to find
\begin{equation}\label{eq:div-spherical}
\begin{aligned}
\hbox{div}\,{\bf v}&=\frac{1}{r^{d-1} \sqrt{\hbox{det}\,h_{\alpha\beta}} } 
\frac{\partial}{\partial \bar x^j} \left( r^{d-1}\sqrt{\hbox{det}\, h_{\alpha\beta}}\,\bar v^j\right) \\
=\frac{1}{r^{d-1} \sqrt{\hbox{det}\,h_{\alpha\beta}} } &\,
 \frac{\partial}{\partial r}\left(r^{d-1}\sqrt{\hbox{det}\, h_{\alpha\beta}}\,\theta_j\,v^j\right) 
 + \frac{1}{r\,\sqrt{\hbox{det}\,h_{\alpha\beta}} } \,
  \frac{\partial}{\partial \omega^\alpha}\left(\sqrt{\hbox{det}\,h_{\alpha\beta}}\,\omega_j^\alpha v^j \right).
\end{aligned}
\end{equation}

We can use (\ref{eq:div-spherical}) to compute the divergence of a vector field of the form ${\bf v}=r^{-k}\,{\bf a}_k$ where
$k\geq 1$ and ${\bf a}_k$ is a vector field with components $a_k^j(\omega)\in H^{m+N+1-k,p}( S^{d-1})$
for $j=1,\dots,d$. We conclude
\begin{equation}\label{eq:div(r^(-k)a_k}
\hbox{div} (r^{-k}{\bf a}_k)=\frac{d-1-k}{r^{k+1}}\,\theta_j \,{a_k^j}(\omega)
+\frac{1}{r^{k+1}\sqrt{h}}\,\frac{\partial}{\partial\omega^\alpha} \left(\sqrt{h}\, \omega_j^\alpha\,a_k^j(\omega)\right),
\end{equation}
where we have used the abbreviation $\sqrt{h}$ for $\sqrt{{\rm det}\,h_{\alpha\beta}}$. 

Now we consider ${\bf u}\in \A^{m,p}_{1,N}$ and claim that $\hbox{div}\,{\bf u}\in \widetilde{\A}^{m-1,p}_{2,N+1}$. In fact, 
using \eqref{eq:div(r^(-k)a_k} we see that $\hbox{div}\,{\bf u}=\chi(r^{-2}{c_2}+\cdots+r^{-N-1}c_{N+1})+g$ where $c_d$ is given by 
\[
c_d(\om)=\frac{1}{\sqrt{h}}\,\frac{\partial}{\partial\omega^\alpha} \left(\sqrt{h}\, \omega_j^\alpha\,a_{d-1}^j(\omega)\right).
\]
Since $c_d$ is a divergence on $S^{d-1}$, we conclude $\int_{S^{d-1}} b_d\,ds=0$ and so $\hbox{div}\,{\bf u}\in\widetilde{\A}^{m-1,p}_{2,N+1}$. 
By Proposition \ref{pr:inverseLaplacian2} we have
$w=K(\hbox{div}\,{\bf u})\in \A^{m+1,p}_{0^*,N-1}$, and hence $\nabla w\in \A^{m,p}_{1,N}$. In fact, the map
${\bf u}\mapsto \nabla K(\hbox{div}\,{\bf u})$ is bounded $\A^{m,p}_{1,N}\to\A^{m,p}_{1,N}$.

Finally, we let ${\bf v}={\bf u}-\nabla w$. Then ${\bf v}\in{\mathcal A}_{1,N}^{m,p}$ and if we compute the divergence, we get
\[
{\rm div}\, {\bf v}={\rm div}\, {\bf u}-\Delta w =0.
\]
Thus ${\bf u}={\bf v}+\nabla w$ satisfies \eqref{eq:Helmholtz}. Now let us confirm uniqueness. If we had 
\[
{\bf u}={\bf v}_1+\nabla w_1={\bf v}_2+\nabla w_2, \quad{\rm div}\,{\bf v}_1={\rm div}\,{\bf v}_2=0,
\]
then $\nabla (w_1-w_2)={\bf v}_2-{\bf v}_1$, and we take divergence to conclude $\Lap(w_1-w_2)=0$.
But our assumption $(w_1-w_2)(x)=o(|x|)$ as $|x|\to\infty$ then implies $w_1-w_2=const$, and we see that ${\bf v}_2-{\bf v}_1=0$,
i.e.\ {\bf v}  is unique. Thus $P_0:{\bf u}\to{\bf v}$ is well-defined and bounded 
${\mathcal A}_{1,N}^{m,p}\to{\mathcal A}_{1,N}^{m,p}$. 
If ${\rm div}\,{\bf u}=0$ then $w=\Lap^{-1}{\rm div}\,{\bf u}=0$, so $P_0\,{\bf u}={\bf u}$. In particular,
$P_0^2\,{\bf u}=P_0\,{\bf u}$, so $P_0$ is indeed a projection. \hfill $\Box$

\medskip
The restriction $N\leq d-1$ is necessary to avoid $\log r$ terms in the Helmholtz decomposition. To see this, let us consider an example.

\begin{example} \label{ex:d=N=2} 
Let us consider $d=N=2$ and try to obtain the Helmholtz decomposition (\ref{eq:Helmholtz}). In fact, let 
$x=r\,\cos\phi$, $y=r\sin\phi$, and simply consider
\[
{\bf u}=\chi(r)\,\frac{{\bf a}_2(\phi)}{r^2}, \quad \phi\in S^1.
\]
If we try to find $b_2(\phi)$ satisfying
\[
\Lap\left(\frac{b_2(\phi)}{r}\right)={\rm div}\, \left(\frac{{\bf a}_2(\phi)}{r^2}\right),
\]
a computation shows that $b_2$ must satisfy
\[
\partial_\phi^2 b_2+b_2=-\sin\phi\, \frac{\partial a_2^1}{\partial\phi}+\cos\phi\, \frac{\partial a_2^2}{\partial\phi}
-2\cos\phi\,a_2^1-2\sin\phi\,a_2^2\,.
\]
But to solve this, we must have the right hand side orthogonal to Ker$(\partial_\phi^2+1)=\{\cos\phi,\sin\phi\}$, which need not be the case. Consequently, we must modify our solution: replace $r^{-1}b_2(\phi)$ by 
$r^{-1}(b_2(\phi)+(c_1\,\cos\phi+c_2\,\sin\phi)\log r)$ where the constants $c_1$, $c_2$ are chosen so that the right hand side of
\[
\partial_\phi^2 b_2+b_2=2(c_1\cos\phi+c_2\sin\phi)-\sin\phi\, \frac{\partial a_2^1}{\partial\phi}+\cos\phi\, \frac{\partial a_2^2}{\partial\phi}
-2\cos\phi\,a_2^1-2\sin\phi\,a_2^2\,
\]
is orthogonal to $\{\cos\phi,\sin\phi\}$, and hence we can find $b_2$. This means that $w$ contains the following asymptotic:
\[
\frac{b_2(\phi)+(c_1\,\cos\phi+c_2\,\sin\phi)\log r}{r}.
\]
Consequently, the terms  $\nabla w$ and ${\bf v}$ in (\ref{eq:Helmholtz}) will both contain asymptotics of the order $r^{-2}\log r$ as $r\to\infty$, and hence will not be in our asymptotic function space ${\mathcal A}_{1,2}^{m,p}(\RR^2,\RR^2)$.
\end{example}

\section{Groups of Asymptotic Diffeomorphisms on $\RR^d$}

In this section we state the main results for diffeomorphisms of $\RR^d$ whose asymptotic behavior can be  described in terms of the asymptotic spaces $\AH_{n,N}^m$ and $\A_{n,N}^m$; proofs of all results will be given in the next two sections.  Denote by $\Diff_+^1(\RR^d,\RR^d)$ the group of orientation-preserving $C^1$-diffeomorphisms on $\RR^d$.
For simplicity of notation, we will no longer use bold face for vector-valued functions as we did in the previous two sections. 
 \begin{definition} \label{def:asym-diffeos} 
For  integers $m> 1+d/p$ and $N\ge n\ge 0$, define
\[
\AH\D_{n,N}^{m,p}(\RR^d,\RR^d):=\{\phi\in\Diff_+^1(\RR^d,\RR^d)\,|\,\phi(x)=x+u(x),\ u\in\AH_{n,N}^{m,p}(\RR^d,\RR^d)\}
\]
and
\[
\A\D_{n,N}^{m,p}(\RR^d,\RR^d):=\{\phi\in\Diff_+^1(\RR^d,\RR^d)\,|\,\phi(x)=x+u(x), \ u\in\A_{n,N}^{m,p}(\RR^d,\RR^d)\}\,.
\]
\end{definition}
\noindent
Similar to Section 2, we will 
abbreviate these collections  by $\AH\D_{n,N}^{m,p}$ and $\A\D_{n,N}^{m,p}$ when it is clear that we are considering diffeomorphisms of $\RR^d$; we also let $\AH\D_{N}^{m,p}$ and $\A\D_{N}^{m,p}$ denote $\AH\D_{0,N}^{m,p}$ and $\A\D_{0,N}^{m,p}$ respectively. 

Now we list some important properties of these spaces of asymptotic diffeomorphisms; as stated before, proofs will be given in the next section. First, we want to show that $\AH\D_{N}^{m,p}$ and $\A\D_{N}^{m,p}$ are topological groups under composition of functions. Since $\phi=Id+u$ means $\phi(\psi)=\psi+u(\psi)$, we see that continuity of $(\phi,\psi)\to \phi(\psi)$ in $\phi$ reduces to continuity of $(u,\psi)\to u(\psi)$. Consequently, we need the following:

\begin{proposition}\label{pr:compositioncontinuous} 
For integers $m> 1+d/p$ and $N\geq n\geq 0$, composition $(u,\psi)\mapsto u\circ\psi$ defines  continuous mappings
\begin{subequations}\label{eq:compositionAxAD->A}
\begin{equation}\label{eq:composition1}
 {\AH}_{n,N}^{m,p} \times {\AH \mathcal D}_{N}^{m,p} \to  {\AH}_{n,N}^{m,p} \quad\hbox{and}\quad
 {\mathcal A}_{n,N}^{m,p} \times {\mathcal A \mathcal D}_{N}^{m,p} \to  {\mathcal A}_{n,N}^{m,p},
\end{equation}
and $C^1$-mappings
\begin{equation}\label{eq:composition2}
 {\AH}_{n,N}^{m+1,p} \times {\AH \mathcal D}_{N}^{m,p} \to  {\AH}_{n,N}^{m,p} \quad\hbox{and}\quad
 {\mathcal A}_{n,N}^{m+1,p} \times {\mathcal A \mathcal D}_{N}^{m,p} \to  {\mathcal A}_{n,N}^{m,p}.
\end{equation}
\end{subequations}
\end{proposition}
 
 \noindent
 Next we  need to know that inverses of asymptotic diffeomorphisms are asymptotic diffeomorphisms. Due to the complexity of the asymptotics, this is simplest to prove for one degree of regularity greater than that required for the continuity of  composition.
 
\begin{proposition}\label{pr:inverses} 
For integers $m> 1+d/p$ and $N\geq n\geq 0$, if $\phi\in{\AH \mathcal D}_{n,N}^{m+1,p}$  then
 $\phi^{-1}\in{\AH \mathcal D}_{n,N}^{m+1,p}$, and $\phi\to\phi^{-1}$ defines a $C^1$-mapping 
${\AH \mathcal D}_{n,N}^{m+1,p}\to {\AH \mathcal D}_{n,N}^{m,p}$.
 The same result holds when 
 ${\AH \mathcal D}$ is replaced by ${\mathcal A \mathcal D}$.
\end{proposition}

These two propositions together suggest that ${\AH \mathcal D}_{n,N}^{m,p}$  is a topological group for $m> 2+d/p$, 
but we have not shown that $\phi\to\phi^{-1}$ is continuous ${\AH \mathcal D}_{n,N}^{m,p}\to {\AH \mathcal D}_{n,N}^{m,p}$. 
However, since the topology in ${\AH \mathcal D}_N^{m,p}$  is just a translation of the Banach space topology of ${\AH }_N^{m,p}$, this follows from the result of Montgomery \cite{M}. Analogous statements can be made about ${\mathcal A \mathcal D}_N^{m,p}$. Consequently, once we have proved the two propositions above, we will have shown:

\begin{theorem}\label{th:topgroup}  
For integers $m> 2+d/p$ and $N\geq n\geq 0$, ${\AH \mathcal D}_{n,N}^{m,p}$ and ${\mathcal A \mathcal D}_{n,N}^{m,p}$ are both topological groups under composition.
\end{theorem}

\section{Proof of the Continuity of Composition (Proposition \ref{pr:compositioncontinuous})}

Our first result concerns scalar functions and is useful in taking the composition of partial asymptotic expansions.
Recall from Remark \ref{re:R-equivalentnorms} the asymptotic space ${\mathcal A}_N^{m,p}(B_R^c)$ in the exterior of the ball $B_R$.

\begin{lemma}\label{le:(1+u)^{-alpha}} 
Suppose $m>d/p$, $N\geq 1$, and $\alpha>0$. If $u\in{\mathcal A}_{1,N}^{m,p}(B_R^c)$ satisfies 
 $1+u(x)\geq \e$ for some $\e>0$ and all $|x|>R$,  then $(1+u)^{-\alpha}-1\in {\mathcal A}_{1,N}^{m,p}(B_R^c)$
 and
 \[
 \|(1+u)^{-\alpha}-1\|_{{\mathcal A}_{1,N}^{m,p}}\leq C_\alpha  \|u\|_{{\mathcal A}_{1,N}^{m,p}}.
 \]
 The same result holds with ${\mathcal A}$ replaced by $\AH$.
\end{lemma}

\noindent
{\bf Proof.} 
The hypotheses imply that $u$ is continuous and bounded on $B^c_R$, so we may assume
$-1+\e\leq u(x)\leq M$ for $|x|>R$.
Now $(1+t)^{-\alpha}$ is a smooth function for $-1+\e\leq t\leq M$, so by Taylor's theorem with remainder, we have
\begin{equation}\label{eq:Taylor}
f_\alpha(t):=(1+t)^{-\alpha}=1-\alpha t+\cdots+(-1)^\ell\left[\alpha(\alpha+1)\cdots(\alpha+\ell-1)\right]\frac{t^\ell}{\ell !}+R_\ell(t),
\end{equation}
where $R_\ell$ is a smooth function of $t\in [-1+\e,M]$ satisfying
\begin{equation}\label{eq:TaylorRemainder}
\left| R_\ell^{(j)}(t)\right|\leq C\,|t|^{\ell+1-j} \quad\hbox{for}\ j=0,1,\dots,\ell+1.
\end{equation}
(The standard statement of Taylor's theorem has $j=0$ in (\ref{eq:TaylorRemainder}); but  for $j=1,\dots,\ell+1$ we can first differentiate (\ref{eq:Taylor}) with respect to $t$ and then use the Taylor estimate for $f_\alpha^{(j)}$.)
Hence we can write
\begin{equation}\label{eq:(1+u)^{-alpha}}
(1+u(x))^{-\alpha}=1-\alpha\,u(x)+\cdots+(-1)^\ell \alpha(\alpha+1)\cdots(\alpha+\ell-1)\frac{(u(x))^\ell}{\ell !}+R_\ell(u(x)).
\end{equation}

Now we have assumed  $u\in{\mathcal A}_{1,N}^{m,p}(B_R^c)$, so by Proposition \ref{pr:A-products} we know that $u^2,\dots,u^\ell\in {\mathcal A}_{1,N}^{m,p}(B_R^c)$. Consequently, we will have completed our proof provided we can show
\begin{equation}\label{Claim-W}
\|R_N(u)\|_{W_{\gamma}^{m,p}(B_R^c)}  \leq C\|u\|_{{\mathcal A}_{1,N}^{m,p}}.
\end{equation}

\smallskip\noindent
To prove (\ref{Claim-W}), we need to consider derivatives of $R_\ell(u)$ up to order $m$; for notational simplicity, at this point let us assume $d=1$. If we calculate the first few derivatives
\[
\begin{aligned}
D_x(R_\ell(u))&=R_\ell'(u)u' \\
D_x^2(R_\ell(u))&=R_\ell''(u)(u')^2+R_\ell'(u)u'' \\
D_x^3(R_\ell(u))&=R_\ell'''(u)(u')^3+3R_\ell''(u)u'u''+R_\ell'(u)u'''
\end{aligned}
\]
we see that, for each $k=1,2,\dots,m$ we have
\begin{equation}\label{eq:D(R(u))}
\left\{
\begin{aligned}
&D_x^k(R_\ell(u))=\sum_{j=1}^k R_\ell^{(j)}(u)\,P_j^k(u',u'',\dots,u^{(k)}), \ \hbox{where}\  P_j^k(t_1,\dots,t_k) \ \hbox{is a} \\
&  \hbox{ homogeneous polynomial of degree $j$ and the total number of derivatives is $k$}.
\end{aligned}
\right.
\end{equation}
In fact, we can easily prove (\ref{eq:D(R(u))}) by induction. It is certainly true for $k=1$ (in which case $P_1^1(u')=u'$).
Now assume that (\ref{eq:D(R(u))}) is true for $k$. To prove (\ref{eq:D(R(u))}) for $k+1$, we calculate
\[
\begin{aligned}
D_x^{k+1}(R_\ell(u))&=\sum_{j=1}^k D_x\left[ R_\ell^{(j)}(u)\,P_j^k(u',\dots,u^{(k)})\right] \\
=& \sum_{j=1}^k R_\ell^{(j+1)}(u)\,u'\,P_j^k(u',\dots,u^{(k)})+R^{(j)}_\ell(u)\,D_x\left[P_j^k(u',\dots,u^{(k)})\right].
\end{aligned}
\]
But $u'\,P_j^k(u',\dots,u^{(k)})$ is a homogeneous polynomial of degree $j+1$ with total number of derivatives $k+1$, and $R^{(j)}_\ell(u)\,D_x\left[P_j^k(u',\dots,u^{(k)})\right]$ is a homogeneous polynomial of degree $j$ with total number of derivatives $k+1$. Relabeling, we have (\ref{eq:D(R(u))})  for $k+1$, completing the induction step.

Now we want to use the representation (\ref{eq:D(R(u))}) to estimate $D_x^k(R_\ell(u))$. 
Using $u\in{\mathcal A}_{1,N}^{m,p}$, we have $|u(x)|\leq C\langle x\rangle^{-1}$, and so 
(\ref{eq:TaylorRemainder}) implies
\begin{equation}\label{est:TaylorRemainder}
|R_\ell^{(j)}(u)|\leq C\,\langle u\rangle^{\ell+1-j}\leq C\,\langle x\rangle^{-\ell-1+j}.
\end{equation}
We also have $|u'(x)|\leq C\, \langle x\rangle^{-2}$,\dots, $|u^{(k)}(x)|\leq C\,\langle x\rangle^{-k-1}$
for $0\leq k\leq m-1$. So to estimate $P_j^k(u',\dots,u^{(k)})$, every occurrence of $u$ contributes $\langle x\rangle^{-1}$ and each derivative of $u$ contributes an additional $\langle x\rangle^{-1}$, so we obtain
\begin{equation}\label{est:Pjk-W}
| P_j^k(u',\dots,u^{(k)})|\leq C\,\langle x\rangle^{-j-k}.
\end{equation}
Combining (\ref{est:TaylorRemainder}) and (\ref{est:Pjk-W}) with (\ref{eq:D(R(u))}), we obtain the estimate
\begin{equation}\label{eq:D^k(R_ell(u)-W}
|D_x^k(R_\ell(u))|\leq C\,\langle x\rangle^{-\ell-1-k} \quad\hbox{for}\ k=0,\dots,m.
\end{equation}
Thus $\langle x\rangle^{N+k}|D_x^k(R_\ell(u))|\leq C\langle x\rangle^{N-\ell-1}$, which is in $L^p$ for $\ell\geq N$ and $k=0,\dots,m$, showing that $R_\ell(u)\in W_N^{m,p}(B_R^c)$. Since $W_N^{m,p}(\RR)\subset W_{N-\frac{1}{p}}^{m,p}(B_R^c)$, this proves (\ref{Claim-W}) and hence the Lemma when $u\in{\mathcal A}_{1,N}^{m,p}(B_R^c)$. 

If we instead assume $u\in {\AH }_{1,N}^{m,p}(B_R^c)$ and perform the same calculations, then $|u^{(k)}(x)|\leq C\x^{-1}$ for
$k=0,\dots,m-1$, so (\ref{est:TaylorRemainder}) still holds, but in place of (\ref{est:Pjk-W}) we obtain
\begin{equation}\label{est:Pjk-H}
| P_j^k(u',\dots,u^{(k)})|\leq C\,\langle x\rangle^{-j},
\end{equation}
and in place of (\ref{eq:D^k(R_ell(u)-W}) we have 
\begin{equation}\label{eq:D^k(R_ell(u)-H}
|D_x^k(R_\ell(u))|\leq C\,\langle x\rangle^{-\ell-1} \quad\hbox{for}\ k=0,\dots,m.
\end{equation}
Then $\x^N|D_x^k(R_\ell(u))|\leq C\x^{N-\ell-1}$, which is in $L^p$ for $\ell\geq N$, showing that in place of (\ref{Claim-W}) we have
\begin{equation}\label{Claim-H}
R_\ell(u)\in H_{N}^{m,p}(\RR) \quad\hbox{for}\ \ell\geq N,m.
\end{equation}
But this shows $(1+u)^{-\alpha}\in {\AH}_N^{m,p}(B_R^c)$ and proves the Lemma when $u\in {\AH }_{1,N}^{m,p}(B_R^c)$.
\hfill$\Box$

\medskip
We will need upper and lower bounds on $\langle\phi(x)\rangle$ when $\phi=Id+u$ is an asymptotic diffeomorphism. But these estimates do not require that $\phi$ be a diffeomorphism, so we formulate them simply in terms of the vector function $u$.

\begin{lemma}\label{le:est<phi>} 
If $u\in {\mathcal A}_0^{m,p}$ where $m>d/p$, then there exist positive constants $c_1$ and $c_2$ so that
\begin{equation}\label{eq:est<phi>}
c_1\,\langle x \rangle\leq\langle x+u(x)\rangle \leq c_2\, \langle x\rangle \quad \hbox{for all } x\in\RR^d
\end{equation}
holds uniformly for bounded $\|u\|_{  {\mathcal A}_0^{m,p}}$.
The same result holds with ${\mathcal A}$ replaced by $\AH$.
\end{lemma}

\noindent
{\bf Proof.} By Lemma \ref{le:W-properties} (d),  we have $\|u\|_\infty\leq C\, \|u\|_{  {\mathcal A}_0^{m,p}}=:M$.
Then we can use $2x\cdot u\geq -\frac{1}{2}|x|^2-2|u|^2$ to conclude
\[
|x|\geq 2M\ \Rightarrow\ \frac{1+|x+u|^2}{1+|x|^2}
\geq \frac{1+\frac{|x|^2}{2}-|u|^2}{1+|x|^2}
\geq \frac{1+\frac{M^2}{4}}{1+|x|^2}\geq \frac{1}{4},
\]
and
\[
|x|\leq 2M\ \Rightarrow\ \frac{1+|x+u|^2}{1+|x|^2}\geq \frac{1}{1+4M^2}.
\]
Similarly, we use $2x\cdot u\leq |x|^2+|u|^2$ to conclude
\[
 \frac{1+|x+u|^2}{1+|x|^2}\leq \frac{1+2|x|^2+2|u|^2}{1+|x|^2}\leq \frac{1+2|x|^2+2M^2}{1+|x|^2}\leq c(M).
\]
 These three estimates complete the proof for $u\in {\mathcal A}_0^{m,p}$.
Since Lemma \ref{le:W-properties} (d) also implies $\|u\|_\infty\leq C\, \|u\|_{  {\AH}_0^{m,p}}$, the same proof establishes (\ref{eq:est<phi>}) for $u\in {\AH }_0^{m,p}$.
\hfill$\Box$

\medskip
Our next estimates concern the Jacobian matrix $D\phi$ of an asymptotic diffeomorphism $\phi$. (We adopt the convention that the $i$-th row and $j$-th column element of $D\phi$ is $(D\phi)^i_{j}=\partial \phi^i/\partial x^j$, so the chain rule may be written $D(\phi\circ\psi)=(D\phi\circ\psi)\cdot D\psi$.)
Let $|D\phi|$ denote the sum of the absolute values of all elements of $D\phi$. For 
$\phi\in {\mathcal A \mathcal D}_0^{m,p}$  (or $\phi\in {\AH \mathcal D}_0^{m,p}$) we want to show that the estimate
\begin{equation}\label{eq:Dphi<C}
| D\phi(x)| \leq C \ \hbox{for all $x\in\RR^d$}
\end{equation}
holds uniformly for bounded $\|u\|_{{\mathcal A}_{0}^{m,p}}$.
Moreover, since we assumed that our diffeomorphisms are orientation-preserving, we know that $\det(D\phi(x))>0$, but we need to confirm a lower bound at infinity, so that we have
\begin{equation}\label{eq:0<det(Dphi)}
0<\e\leq \det(D\phi(x))  \ \hbox{for all $x\in\RR^d$}.
\end{equation}
In fact,  we want to show that  (\ref{eq:0<det(Dphi)})  holds {\bf locally uniformly} for $u\in {\mathcal A}_0^{m,p}$ (or $u\in {\AH}_0^{m,p}$): this means that for fixed $\phi^*=Id+u^*\in {\mathcal A \mathcal D}_0^{m,p}$ we can choose $\e$  so that  (\ref{eq:0<det(Dphi)})  holds for all $\phi=Id+u\in {\mathcal A \mathcal D}_0^{m,p}$ with $\|u-u^*\|_{{\mathcal A}_0^{m,p}}\leq\de$ and $\de$ sufficiently small. 

\begin{lemma}\label{le:est-|Dphi|} 
Suppose $m>1+d/p$.
If $\phi=Id+u\in {\mathcal A \mathcal D}_0^{m,p}$, then (\ref{eq:Dphi<C}) holds uniformly for bounded $\|u\|_{{\mathcal A}_0^{m,p}}$ and (\ref{eq:0<det(Dphi)})  holds locally uniformly for $u\in {\mathcal A}_0^{m,p}$. The analogous statement holds for $\phi\in {\AH \mathcal D}_0^{m,p}$.
\end{lemma}

\noindent
{\bf Proof.} Write $\phi(x)=x+u(x)$ where $u\in {\mathcal A}_0^{m,p}$. Then $D\phi(x)=I+Du(x)$ where $Du\in {\mathcal A}_{1,1}^{m-1,p}$ and $|Du(x)|\leq C\|Du\|_{ {\mathcal A}_{1,1}^{m-1,p}}$ using $m-1>d/p$ and Proposition \ref{pr:A-properties}. Thus  (\ref{eq:Dphi<C}) holds uniformly for bounded $\|u\|_{{\mathcal A}_0^{m,p}}$. But $Du\in {\mathcal A}_{1,1}^{m-1,p}$ with $m-1>d/p$ also implies $|Du(x)|\to 0$ as $|x|\to\infty$, so $\phi$ satisfies (\ref{eq:0<det(Dphi)}). However, we need to show that 
(\ref{eq:0<det(Dphi)}) holds locally uniformly.

Now suppose we fix $\phi^*=Id+u^*\in  {\mathcal A \mathcal D}_0^{m,p}$ satisfying (\ref{eq:0<det(Dphi)}), i.e.\ $\det(D\phi^*(x))\geq\e>0$ for all $x\in\RR$.
If we choose $\widetilde u\in {\mathcal A}_0^{m,p}$ with $\|\widetilde u\|_{ {\mathcal A}_0^{m,p}}$ sufficiently small, we can make 
$\|D\widetilde u\|_\infty$ as small as we like. But $\det:\RR^{2d}\to\RR$ is continuous, so we can arrange
that $\phi=\phi^*+\widetilde u$ satisfies
\[
\det(D\phi(x))= \det(D\phi^*(x)+D\widetilde u(x))\geq\frac{\e}{2}
\]
uniformly for small $\|\widetilde u\|_{ {\mathcal A}_0^{m,p}}$.
Thus we can choose $\e>0$ and $\delta>0$ so that (\ref{eq:0<det(Dphi)})  holds not only for $\phi^*$, but for all $\phi=\phi^*+\widetilde u$  with $\|\widetilde u\|_{{\mathcal A}_0^{m,p}}\leq\de$; i.e.\ (\ref{eq:0<det(Dphi)}) holds locally uniformly. 

The proof for $\phi\in  {\AH \mathcal D}_0^{m,p}$ leads to $Du\in H_0^{m-1,p}$ but otherwise is analogous.
\hfill$\Box$

\medskip
Given an asymptotic diffeomorphism $\phi$ in either ${\AH \mathcal D}_0^{m,p}$ or ${\mathcal {AD}}_0^{m,p}$ with $m>1+d/p$,  the inverse function $\phi^{-1}$ is a diffeomorphism, but we want to obtain estimates at infinity.
Letting $x=\phi^{-1}(y)$ in (\ref{eq:est<phi>}) yields in particular that
\begin{equation}\label{eq:est<psi>}
C_1\,\langle y\rangle \leq \langle \phi^{-1}(y)\rangle\leq C_2\,\langle y\rangle \quad\hbox{for all}\ y\in\RR^d,
\end{equation}
where $C_1=1/c_2$ and $C_2=1/c_1$. But since $c_1$ and $c_2$ were uniform  for bounded $\|u\|_{  {\mathcal A}_0^{m,p}(\RR^d)}$, we find the same is true of (\ref{eq:est<psi>}).

Similarly, for $\phi=Id+u$, we want to show
\begin{equation}\label{eq:Dpsi<C}
|D(\phi^{-1})(x)|\leq C_3 \quad\hbox{for all}\ x\in\RR^d,
\end{equation}
and
\begin{equation}\label{eq:0<det(DPsi)}
0<\e\leq \det (D(\phi^{-1})(x)) \quad\hbox{for all}\ x\in\RR^d,
\end{equation}
both holding  locally uniformly for bounded $\|u\|_{{\mathcal A}_0^{m,p}}$ (and similarly for  ${\AH}_0^{m,p}$). 

\begin{lemma}\label{le:Dpsi} 
Let $\phi=Id+u\in {\mathcal {AD}}_0^{m,p}$ where $m>1+d/p$. Then  (\ref{eq:Dpsi<C}) and (\ref{eq:0<det(DPsi)}) 
both hold  locally uniformly for $u\in {\mathcal A}_0^{m,p}$. The same is true with ${\mathcal A}$ replaced by $\AH$.
\end{lemma}

\noindent
{\bf Proof.} Let us fix $\phi^*=Id+u^*$ and consider $\phi=Id+u=\phi^*+\widetilde u$ with $\|\widetilde u\|_{{\mathcal A}^{m,p}_N}$ small. Let $\psi=\phi^{-1}$ and let $y=\phi(x)$. Since $D\phi(x)=I+Du(x)$, we can use Lemma \ref{le:est-|Dphi|} to conclude
\begin{equation}\label{eq:|I+Du.psi|<C}
|I+Du\circ\psi(y)|\leq C \ \hbox{for all $y\in\RR^d$,}
\end{equation}
\begin{equation}\label{eq:0<det(I+Du.psi)<C}
0<\e\leq \det(I+Du\circ\psi(y)) \ \hbox{for all $y\in\RR^d$,}
\end{equation}
where $C$ and $\e$ are uniform for $\|\widetilde u\|_{ {\mathcal A}_0^{m,p}}\leq\de$.
If we differentiate $\phi(\psi(y))=\psi(y)+u(\psi(y))=y$, we obtain $D\psi+(Du\circ\psi) D\psi=I$, or
$
(I+Du\circ \psi)D\psi=I.
$
However, (\ref{eq:0<det(I+Du.psi)<C}) shows that $I+Du\circ \psi$ is invertible, so we can write
\begin{equation}\label{eq:Dpsi}
D\psi=(I+Du\circ\psi)^{-1}.
\end{equation}
Consequently we have
\[
\det D\psi=\det[(I+Du\circ\psi)^{-1}]=[\det(I+Du\circ\psi)]^{-1}.
\]
Using (\ref{eq:|I+Du.psi|<C}) and (\ref{eq:0<det(I+Du.psi)<C}), we conclude
\[
0<\e_1\leq \det D\psi(y) \leq C_1 \ \hbox{for all $y\in\RR^d$,}
\]
where $\e_1=1/C_0$ and $C_1=1/\e$ are uniform for $\|\widetilde u\|_{ {\mathcal A}_0^{m,p}}\leq\de$.
In particular, this confirms (\ref{eq:0<det(DPsi)}). To prove (\ref{eq:Dpsi<C}) we want to bound $|D\psi|$ uniformly for small $\|\widetilde u\|_{ {\mathcal A}_0^{m,p}}$. But if we use the adjoint formula for the inverse of a matrix,
\[
(I+Du\circ\psi)^{-1}=\frac{1}{\det(I+Du\circ\psi)}\,{\rm Adj}(I+Du\circ\psi),
\]
we see that (\ref{eq:Dpsi<C}) follows from (\ref{eq:|I+Du.psi|<C}), (\ref{eq:0<det(I+Du.psi)<C}), and (\ref{eq:Dpsi}).

The proof for the result with ${\mathcal A}$ replaced by $\AH$ is strictly analogous.
\hfill$\Box$

\medskip
We now consider properties of the composition $f\circ\phi$ when $f$ is in the remainder space and $\phi$ is an asymptotic diffeomorphism. In our first result,  we allow $f$ to be less regular than $\phi$ since this will be useful for later application. We may assume that $f$ is scalar-valued, and we denote its gradient by $\nabla f$.
 
\begin{lemma}\label{le:bound-f(phi)}  
Suppose $m> 1+d/p$  and $\delta\in\RR$. 
For every $\phi\in{\mathcal A \mathcal D}_0^{m,p}$ and every $0\leq k\leq m$, we have
\begin{equation}\label{est:f(phi)inWmp}
\| f\circ \phi\|_{W_{\delta}^{k,p}}\leq C\,\|f\|_{W_{\delta}^{k,p}}\quad\hbox{for all $f\in W_{\delta}^{k,p}$},
\end{equation}
where $C$ may be taken locally uniformly in $\phi\in{\mathcal A \mathcal D}_0^{m,p}$. The analogous result with $\AH$ replacing $\A$ (and $H_\delta$ replacing $W_{\delta}$) is also true.
\end{lemma}

\noindent
{\bf Proof.}  We prove  (\ref{est:f(phi)inWmp}) by induction. 
For $k=0$, we simply use the change of variables $x=\psi(y)=\phi^{-1}(y)$:
\[
\begin{aligned}
\|f\circ\phi\|_{L_\delta^p}^p&=\int\left(\langle x\rangle^\delta|f\circ\phi(x)|\right)^p\,dx\\
&=\int\left(\langle \psi(y)\rangle^\delta|f(y)|\right)^p\det(D\psi(y))\,dy\\
&\leq C\int \left(\langle y\rangle^\delta |f(y)|\right)^p\,dy=C\,\|f\|_{L_\delta^p}^p,
\end{aligned}
\]
where $C$ can be taken locally uniformly  by Lemma \ref{le:Dpsi}. 

Now we assume  (\ref{est:f(phi)inWmp}) holds for $k<m$ and prove it for $k+1$. It suffices to assume 
$f\in W_{\delta}^{k+1,p}$ and show
\begin{equation}\label{eq:est-D(f(phi))}
\|\nabla(f\circ\phi)\|^p_{W_{\delta+1}^{k,p}}\leq C\,\|f\|^p_{W_{\delta}^{k+1,p}}.
\end{equation}
But $\nabla(f\circ\phi)=(\nabla f\circ\phi)\cdot D\phi$ where $\nabla f\in W_{\delta+1}^{k,p}$ and $D\phi=I+Du$ with $Du\in  {\mathcal A}_{1,1}^{m-1,p}$, so we can use Corollary \ref{co:WxA} concerning products (since $m-1>d/p$)  to conclude
\[
\begin{aligned}
\|\nabla(f\circ\phi)\|^p_{W_{\delta+1}^{k,p}} & =
\|(\nabla f\circ\phi)\cdot (I+Du)\|^p_{W_{\delta+1}^{k,p}} \\
& \leq C\,\| I+Du \|_{{\mathcal A}_0^{m-1,p}} \,\|\nabla f\circ\phi\|^p_{W_{\delta+1}^{k,p}}
\leq C\,\|\nabla f\circ\phi\|^p_{W_{\delta+1}^{k,p}},
\end{aligned}
\]
where $C$ can be chosen uniformly for $\|u\|_{{\mathcal A}_{0}^{m,p}}\leq M$. Now we can apply 
 (\ref{est:f(phi)inWmp}) to $\nabla f$ (with $\delta+1$ in place of $\delta$) to conclude
 \[
 \|\nabla f\circ\phi\|^p_{W_{\delta+1}^{k,p}}\leq C\,\|\nabla f\|^p_{W_{\delta+1}^{k,p}}\leq C\,\|f\|^p_{W_{\delta}^{k+1,p}},
 \]
 where $C$ may be taken locally uniformly for $\phi\in{\mathcal A \mathcal D}_0^{m,p}$. 
 Putting these two inequalities together yields (\ref{eq:est-D(f(phi))}), where $C$  may be taken  locally uniformly for $\phi\in{\mathcal A \mathcal D}_0^{m,p}$. 
 
The proof for $f\in H_\delta^{m,p}$ and $\phi\in \AH{\mathcal D}_0^{m,p}$ leads to $\nabla f\in H_\delta^{m-1,p}$ and
$Du\in H_0^{m-1,p}$ but is otherwise analogous. \hfill$\Box$

 \medskip
 The following result will play an important role in proving  the continuity of $f\circ\phi$ with respect to $\phi$.
 
 \begin{lemma}\label{le:f(phi_j)->f(phi)} 
 Assume $m> 1+d/p$, $\delta\in\RR$, and $f\in C_0^\infty(\RR^d)$. If $\phi_k,\phi\in {\mathcal{AD}}_0^{m,p}$ with 
 $\phi_k\to\phi$ in ${\mathcal{AD}}_0^{m,p}$ as $k\to\infty$, then  $f\circ\phi_k\to f\circ\phi$ in $W_\delta^{m,p}$.
 The same is true with ${\mathcal A}$ replaced by $\AH$ (and $W_\delta^{m,p}$ replaced by $H_\delta^{m,p}$).
 \end{lemma}
 
 \noindent
{\bf Proof.} Since $m> d/p$, we know that $\phi_k$ and $\phi$ are continuous functions with
$\phi_k\to\phi$ uniformly on compact sets in $\RR^d$. Moreover, since $\phi_k(x)=x+u_k(x)$ and $\phi(x)=x+u(x)$ where $u_k$ and $u$ are bounded functions while $f$ has compact support, there is a compact set $K$ such that
\begin{equation}\label{eq:f(phi_k(x))=f(phi(x))=0}
f(\phi_k(x))=0=f(\phi(x))\quad\hbox{for all}\ x\in K^c.
\end{equation}

Now we show $f\circ\phi_k\to f\circ\phi$ in $W_\delta^{\ell,p}$ for $0\leq \ell\leq m$ by induction. For $\ell=0$, we use the estimate
\[
|f(\phi_k(x))-f(\phi(x))|\leq \left(\max_{y\in \RR^d}|Df(y)|\right)\,|\phi_k(x)-\phi(x)|
\]
along with \eqref{eq:f(phi_k(x))=f(phi(x))=0} and the fact that $\phi_k\to\phi$ uniformly on $K$ to conclude that
\[
\int_{\RR^d} \x^{\delta p} |f\circ\phi_k(x)-f\circ\phi(x)|^p\,dx \leq C\int_K  |\phi_k(x)-\phi(x)|^p\,dx \to 0.
\]
Next we assume the result for $0\leq \ell<m$ and prove it for $\ell+1$.
Since we may assume $f$ is scalar-valued, this means showing 
$\nabla(f\circ\phi_k)\to \nabla(f\circ\phi)$ in $W^{\ell,p}_{\delta+1}$.
But we may compute
\[
\nabla(f\circ\phi_k)=(\nabla f)\circ\phi_k \cdot D\phi_k\quad\hbox{and}\quad \nabla(f\circ\phi)=
(\nabla f)\circ\phi \cdot D\phi.
\]
Moreover, we know  $D\phi_k\to D\phi$ in ${\mathcal{AD}}_{1,1}^{m-1,p}$ and (by the induction hypothesis)  $(\nabla f)\circ\phi_k\to (\nabla f)\circ\phi$ in $W^{\ell,p}_{\delta+1}$. Hence, by Corollary \ref{co:WxA} concerning products, we find that
$\nabla(f\circ\phi_k)\to \nabla(f\circ\phi)$ in $W^{\ell,p}_{\delta+1}$, as desired.

The proof for $ {\AH\mathcal D}_0^{m,p}$ and $H_\delta^{m,p}$ is strictly analogous.
\hfill$\Box$

 \medskip
 The previous two lemmas may be used to obtain the following continuity result.
\begin{corollary}\label{co:(f,phi)->f(phi)_continuous} 
Assume $m>1+d/p$ and any $\delta\in\RR$. Then composition $(f,\phi)\mapsto f\circ\phi$ is continuous as a map: 
a) $H_{\delta}^{m,p}\times {\AH \mathcal D}_0^{m,p} \to H_{\delta}^{m,p}$, and 
b) $W_{\delta}^{m,p}\times {\mathcal A \mathcal D}_0^{m,p} \to W_{\delta}^{m,p}$.
\end{corollary}

\noindent
{\bf Proof.}
a) Fix $(f^*,\phi^*)\in H_{\delta}^{m,p}\times {\AH \mathcal D}_0^{m,p}$ and consider a sequence $(f_j,\phi_j)\to(f^*,\phi^*)$
in $H_{\delta}^{m,p}\times {\AH \mathcal D}_0^{m,p}$.
By the triangle inequality
\[
\|f_j\circ\phi_j-f^*\circ\phi^*\|_{H^{m,p}_\delta} \leq \|f_j\circ\phi_j-f^*\circ\phi_j\|_{H^{m,p}_\delta}
+ \|f^*\circ\phi_j-f^*\circ\phi^*\|_{H^{m,p}_\delta}.
\]
Using \eqref{est:f(phi)inWmp},  we have
\[
 \|f_j\circ\phi_j-f^*\circ\phi_j\|_{H^{m,p}_\delta}= \|(f_j-f^*)\circ\phi_j\|_{H^{m,p}_\delta}\leq C_1  \|f_j-f^*\|_{H^{m,p}_\delta}
\]
for sufficiently large $j$. Now let us use the density of $C_0^\infty(\RR^d)$ in $H^{m,p}_\delta$
to find $\widetilde f^*$ such that $\|\widetilde f^*-f^*\|_{H^{m,p}_\delta}$ is small. Now we can write
\[
\begin{aligned}
\|f^*\circ\phi_j-f^*\circ\phi^*\|_{H^{m,p}_\delta}\leq \|(f^*-\widetilde f^*)\circ\phi_j\|_{H^{m,p}_\delta} 
&+ \|\widetilde f^*\circ\phi_j-\widetilde f^*\circ\phi^*\|_{H^{m,p}_\delta} \\
&+ \|(\widetilde f^*- f^*)\circ\phi^*\|_{H^{m,p}_\delta}.
\end{aligned}
\]
Again we can use \eqref{est:f(phi)inWmp} to make
\[
\|(f^*-\widetilde f^*)\circ\phi_j\|_{H^{m,p}_\delta}+ \|(\widetilde f^*- f^*)\circ\phi^*\|_{H^{m,p}_\delta}
\leq C_2\,\|\widetilde f^*-f^*\|_{H^{m,p}_\delta}
\]
for sufficiently large $j$.
Now, given $\e>0$, we first pick $\widetilde f^*$ so that $C_2\,\|\widetilde f^*-f^*\|_{H^{m,p}_\delta}<\e/2$.
Then we pick $J$ sufficiently large that both $C_1  \|f-f^*\|_{H^{m,p}_\delta}<\e/4$ and (by Lemma \ref{le:f(phi_j)->f(phi)})
\[
\|\widetilde f^*\circ\phi_j-\widetilde f^*\circ\phi^*\|_{H^{m,p}_\delta}<\e/4 \quad\hbox{ for $j\geq J$.}
\]
This shows $\|f_j\circ\phi_j-f^*\circ\phi^*\|_{H^{m,p}_\delta}<\e$ for $j\geq J$, i.e.\ $H_{\delta}^{m,p}\times {\AH \mathcal D}_0^{m,p} \to H_{\delta}^{m,p}$ is continuous.

b) The proof is exactly the same as for a).
\hfill$\Box$

 \medskip
The following estimates provide a stronger description of the continuity of $f\circ \phi$  when $f$ has an additional degree of regularity.

\begin{lemma}\label{le:stronger_f-continuity}  
Assume $m> 1+d/p$ and  $\delta\in\RR$. \\
a) Fix $\phi_*\in {\AH\mathcal{D}}_0^{m,p}$.  For $f\in H_\delta^{m+1,p}$ and $\phi \in {\AH\mathcal{D}}_0^{m,p}$ sufficiently close to $\phi_*$ we have
\begin{subequations}\label{est:f(v+tilde-v)-f(v)}
\begin{equation}\label{est:f(v+tilde-v)-f(v)-a}
\|f\circ\phi-f\circ\phi_*\|_{H_{\delta}^{m,p}}\leq C\, \|f\|_{H^{m+1,p}_\delta} \|\phi-\phi_*\|_{{\AH}_0^{m,p}},
\end{equation}
b) Fix $\phi_*\in {\mathcal {AD}}_0^{m,p}$. For $f\in W_\delta^{m+1,p}$ and $\phi \in {\mathcal {AD}}_0^{m,p}$ 
sufficiently close to $\phi_*$ we have
\begin{equation}\label{est:f(v+tilde-v)-f(v)-b}
\|f\circ\phi-f\circ\phi_*\|_{W_{\delta+1}^{m,p}}\leq C\, \|f\|_{W^{m+1,p}_\delta} \|\phi-\phi_*\|_{{\mathcal A}_0^{m,p}}.
\end{equation}
\end{subequations}
In both \eqref{est:f(v+tilde-v)-f(v)-a}  and \eqref{est:f(v+tilde-v)-f(v)-b}, $C$ can be taken uniformly for all $\phi$
in a fixed neighborhood of $\phi_*$.
\end{lemma}

\noindent
{\bf Proof.} For  $\phi=Id+u\in  {\AH\mathcal D}_0^{m,p}$
sufficiently close to $\phi_*=Id+u_*$, let $\widetilde u=\phi-\phi_*=u-u_*$ so that $\phi_*+t\widetilde u\in  {\AH\mathcal D}_0^{m,p}$ for $0\leq t\leq 1$. Now $m> 1+d/p$ implies $f\in C^1$ (in fact $C^2$), so we can write
\begin{equation}\label{eq:f(phi)-f(phi*)}
f\circ\phi-f\circ\phi_*=\int_0^1 \frac{d}{dt} f(\phi_*+t\widetilde u) \,dt
=\int_0^1 (\nabla f)(\phi_*+t\widetilde u)\cdot \widetilde u\,dt.
\end{equation}
By Corollary \ref{co:(f,phi)->f(phi)_continuous}, we know that $t\mapsto \nabla f(\phi_*+t\widetilde u)$ is continuous $[0,1]\to H_{\delta}^{m,p}(\RR^d)$, so this mapping is Riemann integrable. 
Thus we can conclude that
\[
\begin{aligned}
\| f\circ\phi-f\circ\phi_* \|_{H_{\delta}^{m,p}}
&\leq \int_0^1 \| (\nabla f)(\phi_*+t\widetilde u)\cdot \widetilde u\|_{H_{\delta}^{m,p}}\,dt \\
&\leq \int_0^1 \|  (\nabla f)(\phi_*+t\widetilde u)\|_{H_{\delta}^{m,p}} \,dt \,\|\widetilde u\|_{{\AH}_0^{m,p}},
\end{aligned}
\]
where we have also used Corollary \ref{co:WxA}.
But now we can apply Lemma \ref{le:bound-f(phi)} to $\nabla f\in H_{\delta}^{m,p}(\RR^d)$:
\[
\|  (\nabla f)(\phi_*+t\widetilde u) \|_{H_{\delta }^{m,p}}
\leq C\,\|\nabla f \|_{H_{\delta}^{m,p}} \leq C\,\|f\|_{H^{m+1,p}_\delta} ,
\]
where $C$ can be taken uniform in a neighborhood of $\phi_*$. Putting this together yields (\ref{est:f(v+tilde-v)-f(v)-a}).

The proof of (\ref{est:f(v+tilde-v)-f(v)-b}) is analogous, except now $\nabla f\in W^{m,p}_{\delta+1}$ so the estimates
become
\[
\begin{aligned}
\| f\circ\phi-f\circ\phi_* \|_{W_{\delta+1}^{m,p}}
&\leq \int_0^1 \|  (\nabla f)(\phi_*+t\widetilde u)\|_{W_{\delta+1}^{m,p}} \,dt \,\|\widetilde u\|_{{\mathcal A}_0^{m,p}}\\
& \leq  C\,\|\nabla f \|_{W_{\delta+1}^{m,p}} \,\|\widetilde u\|_{{\mathcal A}_0^{m,p}}\leq C\,\|f\|_{W^{m+1,p}_\delta} \,\|\phi-\phi_*\|_{{\mathcal A}_0^{m,p}},
\end{aligned}
\]
where $C$ can be taken uniform in a neighborhood of $\phi_*$. 
\hfill$\Box$

\medskip
Before considering the continuity of $a\circ\phi$ when $a$ is an asymptotic function
and $\phi$ is an asymptotic diffeomorphism,  we need a refinement of Lemma
\ref{le:(1+u)^{-alpha}}. For $u\in {\mathcal A}_N^{m,p}(B_1^c,\RR^d)$ with $m>d/p$, let us introduce the {\it scalar-valued} function $\rho(u)$ defined by
\begin{equation}\label{def:rho}
\rho(u)(x):=\frac{2x\cdot u(x)+|u(x)|^2}{|x|^2}.
\end{equation}
Using Propositions \ref{pr:A-properties} and \ref{pr:A-products}, we see that $\rho: {\mathcal A}_N^{m,p}(B_1^c,\RR^d)\to {\mathcal A}_{1,N+1}^{m,p}(B_1^c)$ is continuous and
 we can calculate the  asymptotics of $\rho(u)$ in terms of the asymptotics of $u$. 
 In fact, since we also know by Proposition \ref{pr:A-properties} that $u$ is bounded on $\RR^d$, we see that 
$\rho(u(x))\to 0$ as $|x|\to\infty$, and hence we have $1+\rho(u)\geq \e>0$ for $|x|>R$ with $R$ sufficiently large. Note that $R$ depends on $u$, but we can take it uniformly on a bounded neighborhood ${\mathcal U}$ of a 
fixed $u^*\in  {\mathcal A}_N^{m,p}(B_1^c,\RR^d)$.
Using Lemma \ref{le:(1+u)^{-alpha}},  the  scalar-valued function $\sigma(u)$ defined by
\begin{equation}\label{def:sigma}
\sigma(u)(x):=(1+\rho(u(x)))^{-1/2}-1
\end{equation}
is in ${\mathcal A}_{1,N+1}^{m,p}(B_R^c)$ and we can calculate its asymptotics in terms of the asymptotics of $u$, so we have $\sigma(u)\in {\mathcal A}_{1,N+1}^{m,p}(B_R^c)$. We now want to show that we can choose ${\mathcal U}$ so that $\sigma: {\mathcal U} \to {\mathcal A}_{1,N+1}^{m,p}(B_R^c)$ is real-analytic; in particular, this map is continuous.

\begin{lemma}\label{le:sigma-continuity} 
If $u^*\in {\mathcal A}_N^{m,p}(B_1^c,\RR^d)$ for $m>d/p$ and $N\geq 0$, then there is a neighborhood $ {\mathcal U}$ of $u^*$ and $R$ sufficiently large that $\sigma:  {\mathcal U} \to {\mathcal A}_{1,N+1}^{m,p}(B_R^c)$ is real analytic.
The same is true with ${\mathcal A}$ replaced by $\AH$.
\end{lemma}

\noindent   
{\bf Proof.} Let us fix $u^*\in  {\mathcal A}_N^{m,p}(B_1^c,\RR^d)$ and consider $u=u^*+\widetilde u$. As observed above, we know $\rho(u)\in  {\mathcal A}_{1,N+1}^{m,p}(B_1^c)$ and there is a neighborhood ${\mathcal U}$ of $u^*$ such that, for $R$ sufficiently large, we have $\sigma(u)\in {\mathcal A}_{1,N+1}^{m,p}(B_R^c)$ for all $u\in {\mathcal U}$. Now we compute
\[
1+\rho(u)=(1+\rho(u^*))
\left(1+\frac{2x\cdot\widetilde u+2u^*\cdot\widetilde u+|\widetilde u|^2}{(1+\rho(u^*))|x|^2}\right).
\]
But we can take $\eta$ sufficiently small that, for all $\|\widetilde u\|_{{\mathcal A}_N^{m,p}(B_1^c)}<\eta$, we  have
\[
\left\| \frac{2x\cdot\widetilde u+2u^*\cdot\widetilde u+|\widetilde u|^2}{(1+\rho(u^*))|x|^2} \right\|_{{\mathcal A}_{1,N+1}^{m,p}(B_R^c)}
<\frac{1}{2}.
\]
If we let ${\mathcal U}=\{u=u^*+\widetilde u: \|\widetilde u\|_{{\mathcal A}_N^{m,p}(B_1^c)}<\eta\}$ and use the power series $(1+t)^{-1/2}=1-\frac{1}{2}t+\cdots$, we find that
\[
u\mapsto \left(1+\frac{2x\cdot\widetilde u+2u^*\cdot\widetilde u+|\widetilde u|^2}{(1+\rho(u^*))|x|^2}\right)^{-1/2}
\ \hbox{is real analytic\, ${\mathcal U}\to{\mathcal A}_{N+1}^{m,p}(B_R^c)$.}
\]
Consequently, the same is true of $(1+\rho(u))^{-1/2}$, from which the result follows.

The proof for the result with ${\mathcal A}$ replaced by $\AH$ is strictly analogous.
 \hfill $\Box$
 
 \medskip
Another lemma will be useful in controlling the remainder term in $a\circ\phi$ when $a$ is an asymptotic function
and $\phi$ is an asymptotic diffeomorphism. In this lemma we consider a function $b$ on $S^{d-1}$ and extend it to 
$\RR^d\backslash \{0\}$ as a function of some degree of homogeneity; however, the specific degree of homogeneity does not matter since we will only be using the behavior of $b$ near $|x|=1$.
 
 \begin{lemma}\label{le:bound-b(phi)}  
For $m> 1+d/p$ and $N\geq 0$, suppose $b\in H^{m,p}(S^{d-1})$ is extended (by homogeneity of some degree) to $\RR^d\backslash \{0\}$, and $v\in {\mathcal A}_0^{m,p}(\RR^d,\RR^d)$. Then for $R$ sufficiently large we have
 \begin{subequations}
 \begin{equation}\label{est:b(x+v/|x|inW}
\left \| \frac{b(\frac{x+v}{|x|})}{|x|^{N+1}}\right\|_{W_{\gamma_N}^{m,p}(B_R^c)}\leq C  \left\|b\right\|_{H^{m,p}(S^{d-1})},
 \end{equation}
where $C$ is locally uniform in $v\in {\mathcal A}_0^{m,p}$. The analogous estimate for $v\in\AH_0^{m,p}(\RR^d,\RR^d)$ is
 \begin{equation}\label{est:b(x+v/|x|inH}
\left \| \frac{b(\frac{x+v}{|x|})}{|x|^{N^*+1}}\right\|_{H_{N}^{m,p}(B_R^c)}\leq C  \left\|b\right\|_{H^{m,p}(S^{d-1})}.
 \end{equation}
 \end{subequations}
\end{lemma}

\noindent   
{\bf Proof.} Writing $\gamma=\gamma_N$, we shall prove by induction that
 \begin{equation}\label{est:b(x+v)/|x|^(N+1)}
b\in H^{\ell,p}(S^{d-1}) \quad\Rightarrow\quad \left \| \frac{b(\frac{x+v}{|x|})}{|x|^{N+1}}\right\|_{W_{\gamma}^{\ell,p}(B_R^c)}\leq C  \left\|b\right\|_{H^{\ell,p}(S^{d-1})}
\quad\hbox{for } \ell=0,\dots,m.
 \end{equation}
For $\ell=0$ we easily obtain
\[
\left\| \frac{b(\frac{x+v}{|x|})}{|x|^{N+1}} \right\|^p_{L_\gamma^{p}(B_R^c)}
\leq C\,\int_{R}^\infty\,r^{(\gamma-N-1) p+d-1}\, I(r) \,dr
\]
where
\[
I(r)=\int_{S^{d-1}}  \left| \, b\left(\theta+\frac{v(r\theta)}{r}\right)\right|^p ds .
\]
We first want to show that, for $R$ sufficiently large depending locally uniformly on $v$,  $I(r)$ can be estimated by $C\,\|b\|^p_{L^p(S^{d-1})}$.
To do this we consider the surface $\Xi$ in $\RR^{d}\backslash\{0\}$ parameterized by $\theta\in S^{d-1}$:
\[
\xi(\theta)=\theta + \frac{v(r\theta)}{r}=\frac{x+v(x)}{|x|}.
\]
We compute the Jacobian:
\[
\frac{\partial \xi^i}{\partial \theta^j}=\delta_{ij} +\sum_{k=1}^d \frac{\partial v^i}{\partial x^k}\delta_{jk}
=\delta_{ij}+\frac{\partial v^i}{\partial x^j}.
\]
But $\nabla v\in {\mathcal A}_{1,N+1}^{m-1,p}$ with $m>1+d/p$ implies by Proposition \ref{pr:A-properties} (d) that
$|\nabla v(x)|\leq C\,\|v\|_{ {\mathcal A}_{1,N+1}^{m-1,p}}/|x|$. So,
for $R$ sufficiently large, we conclude that for $r>R$ the Jacobian is nonsingular and we have 
\[
I(r)=\int_{S^{d-1}}  \left| b\left(\theta+\frac{v(r\theta)}{r}\right)\right|^p ds 
\leq C\,\int_{S^{d-1}} |b(\theta)|^p\,ds,
\]
where $C$ is locally uniform in $v\in {\mathcal A}_0^{m,p}$. Finally, using $\gamma-N-1<-d/p$, we conclude
\[
\int^\infty_{R}\,r^{(\gamma-N-1) p+d-1}\, I(r) \,dr
\leq C\, \|b\|^p_{L^p(S^{d-1})},
\]
which gives us (\ref{est:b(x+v)/|x|^(N+1)}) for $\ell=0$.

Now we assume (\ref{est:b(x+v)/|x|^(N+1)}) for $\ell=m-1$ and prove it for $\ell=m$. It suffices to show
\begin{equation}\label{est:grad-(b/|x|^(N+1)}
\left \| \nabla \left( \frac{b(\frac{x+v}{|x|})}{|x|^{N+1}}\right) \right\|_{W_{\gamma+1}^{m-1,p}(B_R^c)}
\leq C  \left\|b\right\|_{H^{m,p}(S^{d-1})}.
\end{equation}
But
 \[
\frac{\partial}{\partial x^j} \left( \frac{b(\frac{x+v}{|x|})}{|x|^{N+1}}\right)
=(\partial_i b) \left(\frac{x+v}{|x|}\right) \left( \frac{\delta_{ij}-\theta_i\theta_j+\partial_j v^i}{|x|} \right) |x|^{-N-1}
-(N+1)b\left(\frac{x+v}{|x|}\right)|x|^{-N-2}\theta_j.
\]
We can use Lemma \ref{le:af-product} to estimate
\[
\begin{aligned}
\left\| (\partial_i b)  \left(\frac{x+v}{|x|}\right) |x|^{-N-2} \right\|_{W_{\gamma+1}^{m-1,p}(B_R^c)}
\leq C\,\left\| (\partial_i b)  \left(\frac{x+v}{|x|}\right) |x|^{-N-1} \right\|_{W_{\gamma}^{m-1,p}(B_R^c)} \\
\leq C\, \| \partial_j b\|_{H^{m-1,p}(S^{d-1})} \leq C\, \| b \|_{H^{m,p}(S^{d-1})},
\end{aligned}
\]
where we have used the induction hypothesis for $\ell=m-1$ applied to $\partial_j b\in H^{m-1,p}(S^{d-1})$.
We can also apply the induction hypothesis to estimate
\[
\begin{aligned}
\left\| \,b\left(\frac{x+v}{|x|}\right)|x|^{-N-2} \right\|_{W_{\gamma+1}^{m-1,p}}
=\left\| \, b\left(\frac{x+v}{|x|}\right)|x|^{-N-1} \right\|_{W_{\gamma}^{m-1,p}} \\
\leq C\,\|b\|_{H^{m-1}(S^{d-1})} \leq C\, \|b\|_{H^{m}(S^{d-1})}.
\end{aligned}
\]
Putting these together proves (\ref{est:grad-(b/|x|^(N+1)}), which completes the induction.  

For $v\in\AH_0^{m,p}$ we follow the same outline, using $N-N^*-1<-d/p$ to conclude convergence of the radial integral. \hfill $\Box$

\medskip
We now consider compositions $u\circ\phi$ when $u=a+f$ as in (\ref{asymptoticexpansion1}). We may assume that $u$ is scalar-valued but the diffeomorphism $\phi=Id+u$ is necessarily vector-valued. We start with  generalizing Lemma \ref{le:bound-f(phi)}.

\begin{lemma}\label{le:bound-u(phi)}  
Suppose $m> 1+d/p$ and $N\geq n\geq 0$. For any $\phi\in {\mathcal {A D}}_{N}^{m,p}$ we have
 \begin{equation}\label{est:u(phi)inA}
\left \|u\circ\phi\right\|_{{\mathcal A}_{n,N}^{m,p}}\leq C  \left\|u\right\|_{{\mathcal A}_{n,N}^{m,p}} \quad\hbox{for all}\ u\in {\mathcal A}_{n,N}^{m,p},
 \end{equation}
where $C$ may be taken locally uniformly in $\phi\in {\mathcal {A D}}_{N}^{m,p}$. The analogous result with $\AH$ replacing $\A$  is also true.
\end{lemma}

\noindent   
{\bf Proof.} To simplify notation, we  assume $n=0$.
 Using the form \eqref{def:AW-norm} of the ${\mathcal A}$-norm and Lemma \ref{le:bound-f(phi)}, it suffices to consider 
\begin{equation}\label{a=ak}
u(x)=a(x)=\chi(|x|)a_k(\theta)/|x|^k \quad\hbox{where}\ a_k\in H^{m+1+N-k,p}(S^{d-1})
\ \hbox{and}\ 0\leq k\leq N.
\end{equation}
 Moreover, since $\phi=Id+v$ where $v\in {\mathcal A }_{N}^{m,p}\subset C_B(\RR^d)$, we may assume that $|v(x)|\leq M$. Since $\phi(x)=x+v(x)$, for $|x|>2M$ we have $M\leq |\phi(x)|\leq |x|+M\leq 3|x|/2$. Let us assume $M\geq 1$ and let $R=2M$; for $x\in B_R^c$ we have $\chi(|x|)=1=\chi(|\phi(x)|)$, so it suffices to estimate
  $a\circ\phi$ in $|x|>R$ in terms of $a_k$ on $S^{d-1}$. So we want to show
  \begin{equation}\label{est:a_k(phi)}
  \|a\circ\phi\|_{{\mathcal A }_{N}^{m,p}(B_R^c)}\leq C\,  \|a_k\|_{H^{m+1+N-k,p}(S^{d-1})}
  =C  \left\|a\right\|_{{\mathcal A}_N^{m,p}},
  \end{equation}
where $C$ is locally uniform in $v\in {\mathcal A }_{N}^{m,p}$. But to estimate $\|a\circ\phi\|_{{\mathcal A }_{N}^{m,p}(B_R^c)}$, we need a partial asymptotic expansion for $a\circ\phi$.

Let us consider $a_k(x)$ as
a homogeneous of degree $0$ function on $\RR^d\backslash\{0\}$. In particular, $a_k\in H^{m+1+N-k,p}_{\loc}(\RR^d\backslash\{0\})\subset C^{N-k+1}(\RR^d\backslash\{0\})$ since $m>d/p$.
 So, by Taylor's theorem with remainder at a point $y^*\in \RR^d\backslash\{0\}$, we can write
 \begin{subequations} \label{ak-Taylor}
 \begin{equation}\label{eq:ak-Taylor1}
a_k(y)=\sum_{|\alpha|\leq N-k} D^\alpha a_k(y^*)\,\frac{(y-y^*)^\alpha}{\alpha!}+R_{N,k}(y,y^*),
 \end{equation}
 where the remainder  $R_{N,k}(y,y^*)$ can be expressed in integral form as
 \begin{equation}\label{eq:a0-Taylor2}
R_{N,k}(y,y^*) =\sum_{|\alpha|=N-k+1}\frac{N-k+1}{\alpha!}\,\int_0^1 (1-t)^{N-k} D^\alpha  a_k(y^*+t(y-y^*))\,dt\,(y-y^*)^\alpha.
 \end{equation}
 \end{subequations}
This approximation holds for $y$ in a neighborhood of $y^*$, and more generally provided $0\not\in\{y^*+t(y-y^*):0\leq t\leq 1\}$. But we now want to take both $y$ and $y^*$ on $ S^{d-1}$. In fact, we shall replace $y$ by $\phi(x)/|\phi(x)|$ and $y^*$ by $\theta=x/|x|$:
 \begin{equation}\label{eq:ak-Taylor3}
 \begin{aligned} 
a_k\left(\frac{\phi(x)}{|\phi(x)|}\right) 
  =  \sum_{|\alpha|\leq N-k} \frac{D^\alpha  a_k(\theta)}{\alpha!}\left(\frac{\phi(x)}{|\phi(x)|}-\frac{x}{|x|}\right)^\alpha+R_{N,k}\left(\frac{\phi(x)}{|\phi(x)|},\frac{x}{|x|}\right).
 \end{aligned}
 \end{equation}
Notice that $\phi(x)=x+v(x)$, where $v$ is bounded, means that $\phi(x)/|\phi(x)|\to x/|x|$ as $|x|\to \infty$, so for $|x|>R$ with $R$ sufficiently large we can arrange
\[
0\not\in \left\{\frac{x}{|x|}+t\left(\frac{\phi(x)}{|\phi(x)|}-\frac{x}{|x|}\right): 0\leq t\leq 1\right\}.
\]
But we need to investigate the difference $\phi(x)/|\phi(x)|-x/|x|$ in more detail.

Notice that we can write
 \begin{equation}\label{eq:|phi|^(-1)}
|\phi(x)|^{-1}= |x+v|^{-1}=|x|^{-1}(1+\rho(v)(x))^{-1/2}=|x|^{-1}(1+\sigma(v)(x)),
 \end{equation}
 where $\rho(v)$ and $\sigma(v)$ are defined in (\ref{def:rho}) and (\ref{def:sigma}) respectively.
But by Lemma \ref{le:sigma-continuity} we know that $\|\sigma(v)\|_{{\mathcal A}_{1,N+1}^{m,p}(B_R^c)}$ is  bounded locally uniformly for $v$ in  ${\mathcal A}_{N}^{m,p}(B_1^c)$.
It is easy to confirm that
\begin{equation}\label{eq:phi/|phi|-x/|x|}
\frac{\phi(x)}{|\phi(x)|}-\frac{x}{|x|}
=\frac{w(x)}{|x|} \quad\hbox{where}\ w(x):=\sigma(v)(x)\,[x+v(x)] + v(x).
\end{equation}
Note that $w \in {\mathcal A}_{N}^{m,p}(B_R^c)$ and that its asymptotics can be computed in terms of $v$; in particular, 
$\|w\|_{ {\mathcal A}_{N}^{m,p}(B_R^c)}$ is bounded locally uniformly in $v\in  {\mathcal A}_{N}^{m,p}(B_1^c)$.

We plug (\ref{eq:phi/|phi|-x/|x|}) and (\ref{eq:|phi|^(-1)}) into (\ref{eq:ak-Taylor3}) to conclude
\begin{equation}\label{eq:asymexp(a(phi))}
\begin{aligned}
a\circ\phi(x)&=\chi(|\phi(x)|)\,a_k\left(\frac{\phi(x)}{|\phi(x)|}\right) |\phi(x)|^{-k}\\
&=\chi(|\phi(x)|)\,\left[ \sum_{|\alpha|\leq N-k} \frac{D^\alpha  a_k(\theta)\,w^\alpha}{\alpha!\,|x|^{k+|\alpha|}}
 +R_{N,k}\left(\frac{\phi(x)}{|\phi(x)|},\frac{x}{|x|}\right)\frac{1}{|x|^k}\right]\,(1+\sigma(v))^{k}.
 \end{aligned}
\end{equation}
Although this is not quite the partial asymptotic expansion for $a\circ\phi$, it can be used to estimate 
$\|a\circ\phi\|_{{\mathcal A}^{m,p}_{N}}$. In fact, using Proposition \ref{pr:A-products} and Lemma \ref{le:sigma-continuity}, we know $\|a\,(1+\sigma(v))^k\|_{{\mathcal A}_{N}^{m,p}(B_R^c)} \leq C\,\|a\|_{{\mathcal A}_{N}^{m,p}(B_R^c)}$, where
$C$ is locally uniform in $v\in  {\mathcal A}_{N}^{m,p}(B_1^c)$; consequently, to estimate $\|a\circ\phi\|_{{\mathcal A}_{N}^{m,p}(B_R^c)}$ we need only estimate the ${\mathcal A}_{N}^{m,p}(B_R^c)$-norm of the two terms in the brackets in \eqref{eq:asymexp(a(phi))}. 

First of all, we claim
\begin{equation}\label{est:A-norm(Taylor-asymptotics)}
\left\| \sum_{|\alpha|\leq N-k}  \frac{D^\alpha  a_k(\theta)\,w^\alpha}{|x|^{k+|\alpha|}}\right\|_{{\mathcal A}_{N}^{m,p}(B_R^c)}
\leq C\,\|a\|_{{\mathcal A}_{N}^{m,p}(B_R^c)},
\end{equation}
where $C$ is locally uniform in $v\in  {\mathcal A}_{N}^{m,p}(B_1^c)$.
To see this, we use the algebra property to estimate
\[
\begin{aligned}
\left\| \sum_{|\alpha|\leq N-k}  \frac{D^\alpha  a_k(\theta)\,w^\alpha}{\alpha !\,|x|^{k+|\alpha|}}\right\|_{{\mathcal A}_{N}^{m,p}(B_R^c)} 
& \leq \sum_{|\alpha|\leq N-k} \left\| \frac{D^\alpha a_k(\theta)}{|x|^{k+|\alpha|}} \right\|_{{\mathcal A}_{N}^{m,p}(B_R^c)} 
\|w^\alpha\|_{{\mathcal A}_{N}^{m,p}(B_R^c)} 
\\
 \leq & \sup_{|\alpha|\leq N-k}\|w^\alpha\|_{{\mathcal A}_{N}^{m,p}(B_R^c)}  \sum_{|\alpha|\leq N-k} \left\| D^\alpha  a_k(\theta)\right\|_{H^{m+1+N-k-|\alpha|,p}(S^{d-1})} 
  \\
\leq & \, C\,\|a_k\|_{H^{m+1+N-k,p}(S^{d-1})}=C\,\|a\|_{{\mathcal A}_{N}^{m,p}(B_R^c)},
 \end{aligned}
\]
where $C$ is locally uniform in $w\in {\mathcal A}_{N}^{m,p}(B_R^c)$.
If we recall  that $v\mapsto w$ is continuous as a map ${\mathcal A}_{N}^{m,p}(B_1^c)\to {\mathcal A}_{N}^{m,p}(B_R^c)$, then we can consider $C$ in \eqref{est:A-norm(Taylor-asymptotics)} as being locally uniform in $v$.

Secondly, we claim
\begin{equation}\label{est:|RN|}
\left\|\frac{R_{N,k}(\phi/|\phi|,x/|x|)}{|x|^{k}}\right\|_{W_\gamma^{m,p}(B_R^c)}\leq C\,\|a_k\|_{H^{m+1+N-k,p}(S^{d-1})},
\end{equation}
where $C$ may be taken locally uniform in $v\in {\mathcal A}_{N}^{m,p}(B_1^c)$. To see this, notice from (\ref{eq:a0-Taylor2}) that
\begin{equation}\label{eq:TaylorRemainder(phi)}
\frac{R_{N,k}(\phi/|\phi|,x/|x|)}{|x|^{k}}
=\sum_{|\alpha|=N-k+1}\frac{N-k+1}{\alpha!}\int_0^1(1-t)^{N-k}D^\alpha a_k\left(\frac{x}{|x|}+t\frac{w(x)}{|x|}\right)dt \,
\frac{w^\alpha}{|x|^{N+1}}.
\end{equation}
We can apply Lemma \ref{le:bound-b(phi)} with $b=D^{\alpha}a_k\in H^{m,p}(S^{d-1})$ to conclude
\[
\left\|\frac{D^\alpha a_k(\frac{x+tw}{|x|})}{|x|^{N+1}}\right\|_{W_\gamma^{m,p}(B_R^c)}\leq C\,\|D^\alpha a_k\|_{H^{m,p}(S^{d-1})}
\leq C\, \|a_k\|_{H^{m+1+N-k,p}(S^{d-1})}.
\]
Using Corollary \ref{co:WxA} and the above remarks regarding $v\mapsto w$, we obtain \eqref{est:|RN|}.
 Putting this together with (\ref{eq:asymexp(a(phi))}) and (\ref{est:A-norm(Taylor-asymptotics)}), we obtain (\ref{est:a_k(phi)}), as desired.
 
 To prove the corresponding result for $\AH$, we replace $N$ by $N^*$ in \eqref{eq:ak-Taylor1} and \eqref{eq:a0-Taylor2} and replace \eqref{est:|RN|} by
 \[
 \left\|\frac{R_{N^*,k}(\phi/|\phi|,x/|x|)}{|x|^{k}}\right\|_{H_N^{m,p}(B_R^c)}\leq C\,\|a_k\|_{H^{m+1+N^*-k,p}(S^{d-1})}.
 \]
 The details are straight-forward.
\hfill 
 $\Box$
 
 \begin{lemma}\label{le:u(phi_j)->u(phi)}  
Suppose $m>1+d/p$ and $N\geq 0$. Let $u(x)=\chi(|x|)a_k(\theta)/|x|^k$ where $n\leq k\leq N$ and $a_k\in C^\infty(S^{d-1})$. If $\phi,\phi_j\in {\mathcal AD}^{m,p}_{N}$ with $\phi_j\to \phi$ in ${\mathcal AD}^{m,p}_{N}$, then $u\circ\phi_j\to u\circ \phi$ in ${\mathcal A}^{m,p}_{n,N}$. The same is true if $\A$ is replaced by $\AH$.
\end{lemma}

\noindent
{\bf Proof.} As in the proof of Lemma \ref{le:bound-u(phi)}, we assume  $n=0$.
Since we can write $\phi_j(x)=x+v_j(x)$ and $\phi(x)=x+v(x)$ where $v_j$ and $v$ are uniformly bounded functions, we can take $R$ large enough that $\chi(|\phi_j(x)|)=\chi(|\phi(x)|)=1$ for all $|x|>R$, so we want to estimate in $\A^{m,p}_{N}(B_R^c)$ the difference
\[
\frac{a_k\left(\frac{\phi_j(x)}{|\phi_j(x)|}\right)}{|\phi_j(x)|^k}
-\frac{a_k\left(\frac{\phi(x)}{|\phi(x)|}\right)}{|\phi(x)|^k}.
\]
Using the scalar function $\sigma$ defined in \eqref{def:sigma}, let us introduce (as we did in \eqref{eq:phi/|phi|-x/|x|}) the vector functions
\[
w_j:=\sigma(v_j)[I+v_j]+v_j \quad\hbox{and}\quad w:=\sigma(v)[I+v]+v.
\]
Since $v_j\to v$ in $\A_{N}^{m,p}$, we see by Lemma \ref{le:sigma-continuity} that $w_j\to w$ in $\A_{N}^{m,p}$. Now if we apply  \eqref{eq:asymexp(a(phi))} to both $u\circ\phi_j$ and $u\circ\phi$, we find for $|x|>R$ that
\begin{equation}\label{eq:u(phi_j)-u(phi)}
\begin{aligned}
u\circ\phi_j-u\circ\phi & =
\sum_{|\alpha|\leq N-k} \frac{D^\alpha a_k(\theta)\,[w_j^\alpha(1+\sigma(v_j))^k-w^\alpha(1+\sigma(v))^k]}{\alpha ! \, |x|^{k+|\alpha|}} \\
+ & \frac{1}{|x|^k}\left[ R_{N,k}\left(\frac{\phi_j(x)}{|\phi_j(x)|},\frac{x}{|x|}\right)(1+\sigma(v_j))^k
- R_{N,k}\left(\frac{\phi(x)}{|\phi(x)|},\frac{x}{|x|}\right)(1+\sigma(v))^k \right].
\end{aligned}
\end{equation}
Using Lemma \ref{le:sigma-continuity} again, for each fixed $\alpha$ we know
\[
\| w_j^\alpha(1+\sigma(v_j))^k-w^\alpha(1+\sigma(v))^k \|_{\A_{N}^{m,p}} \to 0.
\]
As observed in the proof of Lemma \ref{le:bound-u(phi)}, $|x|^{k+|\alpha|}D^\alpha a_k(\theta)\in \A^{m,p}_N(B^c_R)$, and so
\[
\left\| \sum_{|\alpha|\leq N-k} \frac{D^\alpha a_k(\theta)\,[w_j^\alpha(1+\sigma(v_j))^k-w^\alpha(1+\sigma(v))^k]}{\alpha ! \, |x|^{k+|\alpha|}} \right\|_{\A_{N}^{m,p}(B^c_R)} \to 0.
\]
To handle the remainder terms in \eqref{eq:u(phi_j)-u(phi)}, it suffices to show
\[
\frac{1}{|x|^k}\left[ R_{N,k}\left(\frac{\phi_j}{|\phi_j|},\frac{x}{|x|}\right)-
R_{N,k}\left(\frac{\phi}{|\phi|},\frac{x}{|x|}\right)\right]\to 0 \quad \hbox{in $\A^{m,p}_N(B^c_R)$}.
\]
But, using \eqref{eq:TaylorRemainder(phi)}, this quantity is given by
\[
\sum_{|\alpha|=N-k+1}\frac{N-k+1}{\alpha!}
\int_0^1 \frac{(1-t)^{N-k}}{|x|^{N+1}}\left[D^\alpha a_k\left(\frac{x}{|x|}+t\frac{w_j(x)}{|x|}\right)w_j^\alpha
-D^\alpha a_k\left(\frac{x}{|x|}+t\frac{w(x)}{|x|}\right)w^\alpha\right]\,dt.
\]
However, $w_j\to w$ in $\A^{m,p}_N$  implies $w_j^\alpha\to w^\alpha$  in $\A^{m,p}_N$  and also $D^\alpha a_k((x+tw_j(x))/|x|)\to D^\alpha a_k((x+tw(x))/|x|)$ in $\A^{m,p}_N(B^c_N)$, so the remainder term in \eqref{eq:u(phi_j)-u(phi)} also tends to zero in $\A^{m,p}_N(B^c_N)$ as $j\to\infty$, which is what we needed to show. The proof for $\AH$ is identical.
\hfill $\Box$
 
 \medskip
Similar to Corollary \ref{co:(f,phi)->f(phi)_continuous}, we can use Lemmas  \ref{le:bound-u(phi)} and \ref{le:u(phi_j)->u(phi)}  to show that  composition on our asymptotic spaces is continuous; we shall not repeat the argument.
\begin{corollary}\label{co:(u,phi)->u(phi)_continuous} 
Suppose $m>1+d/p$ and $N\geq n\geq 0$. Then composition $(u,\phi)\mapsto u\circ\phi$ is continuous as a map: 
a) $ {\mathcal A}_{n,N}^{m,p} \times {\mathcal A \mathcal D}_{N}^{m,p} \to  {\mathcal A}_{n,N}^{m,p}$, and
b)  $ {\AH}_{n,N}^{m,p} \times {\AH \mathcal D}_{N}^{m,p} \to  {\AH}_{n,N}^{m,p}$.
\end{corollary}

Now we extend Lemma \ref{le:stronger_f-continuity} to general asymptotic functions $u$. Again we may assume that $u$ is scalar-valued. 

\begin{lemma}\label{le:a(phi_j)->a(phi)}  
Suppose $m> 1+d/p$ and $N\geq n\geq 0$. \\
a) For $u\in {\AH}_{n,N}^{m+1,p}$, $\phi_* \in {\AH \mathcal D}_{N}^{m,p}$, and  all
$\phi \in {\AH \mathcal D}_{N}^{m,p}$ sufficiently close to $\phi_*$ we have
\begin{subequations}
\begin{equation}\label{eq:u(phi)-u(phi*)-a}
\|u\circ \phi-u\circ \phi_* \|_{{\AH}_{n,N}^{m,p}}
\leq C\, \|u\|_{{\AH}_{n,N}^{m+1,p}} \|\phi-\phi_*\|_{{\AH}_{N}^{m,p}}.
\end{equation}
b) For $u\in {\mathcal A}_{n,N}^{m+1,p}$, $\phi_* \in {\mathcal A \mathcal D}_{N}^{m,p}$, and  all
$\phi \in {\mathcal A \mathcal D}_{N}^{m,p}$ sufficiently close to $\phi_*$ we have 
\begin{equation}\label{eq:u(phi)-u(phi*)-b}
\|u\circ \phi-u\circ \phi_* \|_{{\mathcal A}_{n,N}^{m,p}}
\leq C\, \|u\|_{{\mathcal A}_{n,N}^{m+1,p}} \|\phi-\phi_*\|_{{\mathcal A}_{N}^{m,p}}.
\end{equation}
\end{subequations}
In both cases, the constant $C$ is locally uniform in $\phi$. 
\end{lemma}

\medskip\noindent
{\bf Proof. } As in the proof of Lemma \ref{le:f(phi_j)->f(phi)},  let $\phi=Id+v$ and $\phi_*=Id+v_*$, and let $\widetilde v=\phi-\phi_*=v-v_*$. Use the fact that $u\in C^1$ to write
\begin{equation}\label{est:u(phi)-u(phi*)}
u\circ \phi-u\circ \phi_*=\int_0^1 (\nabla u)(\phi_*+t\widetilde v)\cdot\widetilde v\,dt.
\end{equation}
Assuming $u\in {\AH}_{n,N}^{m+1,p}$,
we know $\nabla u\in {\AH}_{n+1,N}^{m,p}\subset {\AH}_{n,N}^{m,p}$ and $\widetilde v\in{\mathcal A}_N^{m,p}$. By Corollary \ref{co:(u,phi)->u(phi)_continuous}, the function
$F: t\mapsto  (\nabla u)(\phi_*+t\widetilde v)\cdot\widetilde v$
is continuous as a map $F:[0,1]\to {\AH}_{n,N}^{m,p}$, hence Riemann integrable. Consequently, we can apply the algebra property and then Lemma \ref{le:bound-u(phi)} to \eqref{est:u(phi)-u(phi*)} to obtain the  estimate:
\[
\begin{aligned}
\| u\circ \phi-u\circ \phi_* \|_{{\AH}_{n,N}^{m,p}}  & \leq C \sup_{0\leq t\leq 1} \|(\nabla u)(\phi_*+t\widetilde v)\|_{{\AH}_{n,N}^{m,p}}\,
\,\| \widetilde v \|_{{\AH}_{N}^{m,p}} \\
& \leq C\, \|\nabla u \|_{{\AH}_{n,N}^{m,p}}
\,\| \widetilde v \|_{{\AH}_{N}^{m,p}}.
\end{aligned}
\]
But $\|\nabla u\|_{{\AH}_{n,N}^{m,p}}\leq C\|u\|_{{\AH}_{n,N}^{m+1,p}}$ since $\nabla:\AH_{n,N}^{m+1,p}\to\AH_{n+1,N}^{m,p}\subset\AH_{n,N}^{m,p}$,  so we obtain \eqref{eq:u(phi)-u(phi*)-a}.
Now we consider \eqref{eq:u(phi)-u(phi*)-b}. For $u\in {\mathcal A}_{n,N}^{m+1,p}$, we know $\nabla u\in {\mathcal A}_{n+1,N+1}^{m,p}\subset {\mathcal A}_{n,N}^{m,p}$, so the above steps show
\[
\| u\circ \phi-u\circ \phi_* \|_{{\A}_{n,N}^{m,p}}  \leq C\,\|\nabla u\|_{{\A}_{n,N}^{m,p}}\,\|\widetilde v\|_{{\A}_{N}^{m,p}}
\leq C\,\|u\|_{{\A}_{n,N}^{m+1,p}}\,\|\widetilde v\|_{{\A}_{N}^{m,p}}. \quad \hfill\Box
\]

 \medskip\noindent 
{\bf Proof of Proposition \ref{pr:compositioncontinuous}.} The continuity of \eqref{eq:composition1} is contained in Corollary   \ref{co:(u,phi)->u(phi)_continuous}. 
For the proof that (\ref{eq:composition2}) is $C^1$, we shall abbreviate our notation.
Let $X^m$ represent   ${\AH}_{n,N}^{m,p}$ (or ${\mathcal A}_{n,N}^{m,p}$)
and $X{\mathcal D}^m$ represent   ${\AH \mathcal D}_{n,N}^{m,p}$ (or ${\mathcal {AD}}_{n,N}^{m,p}$); the norm in $X^m$ will be denoted simply by $\|\cdot\|_m$. We can also assume that $u$ is scalar-valued (although we will not distinguish notation between the $m$-norms of vector and scalar-valued functions).

We fix $u\in  X^{m+1}$ and $\phi \in{X \mathcal D}^m$, and consider nearby $u+\delta u\in  X^{m+1}$ and $\phi +\delta\phi \in {X \mathcal D}^m$; note that $\delta u\in  X^{m+1}$ and $\delta\phi \in X^{m}$. We want to show that
\begin{equation}\label{C^1-condition}
(u+\delta u)\circ (\phi+\delta\phi)=u\circ\phi + L_{u,\phi}(\delta u,\delta\phi)+R_{u,\phi}(\delta u,\delta\phi),
\end{equation}
where $L_{u,\phi}:X^{m+1}\times X^{m}\to X^{m}$ is a bounded linear map, and 
$\|R_{u,\phi}(\delta u,\delta\phi)\|_m=o(\|\delta u\|_{m+1}+\|\delta\phi\|_m)$ as $\|\delta u\|_{m+1}+\|\delta\phi\|_m\to 0$.

We shall repeatedly use the following simple identity for $u\in C^1(\RR^d)$:
 \begin{subequations}
 \begin{equation}
 u(y+\delta y)=u(y)+\nabla u(y)\cdot\delta y + R(u,y,\delta y),
\end{equation}
where
\begin{equation}
R(u,y,\delta y)=\int_0^1(\nabla u(y+t\delta y)-\nabla u(y))\cdot \delta y\,dt.
\end{equation}
\end{subequations}
Applying this identity with $y$ replaced by $\phi$, we find
\[
u\circ(\phi+\delta\phi)=u\circ\phi +(\nabla u\circ\phi) \cdot\delta\phi + R_1(u,\phi,\delta\phi),
\]
where
\[
R_1(u,\phi,\delta\phi)=\int_0^1 (\nabla u\circ(\phi+t\delta\phi)-\nabla u\circ\phi)\cdot\delta\phi\,dt.
\]
Then, replacing $u$ by $\delta u$, we find
\[
\delta u\circ(\phi+\delta\phi)=\delta u\circ\phi +(\nabla (\delta u)\circ\phi) \cdot\delta\phi + R_2(u,\delta u,\phi,\delta\phi),
\]
where
\[
R_2(u,\delta u,\phi,\delta\phi)=\int_0^1 (\nabla (\delta u)\circ(\phi+t\delta\phi)-\nabla (\delta u)\circ\phi)\cdot\delta\phi\,dt.
\]
Putting these together, we obtain (\ref{C^1-condition}) where
\begin{equation}
L_{u,\phi}(\delta u,\delta\phi)=\delta u\circ\phi+(\nabla u\circ\phi)\cdot \delta\phi,
\end{equation}
and
\begin{equation}
R_{u,\phi}(\delta u,\delta\phi)=(\nabla (\delta u)\circ\phi)\cdot\delta\phi+R_1(u,\phi,\delta\phi)+R_2(u,\phi,\delta u,\delta\phi).
\end{equation}
Clearly $L_{u,\phi}$ is linear in $\delta u$ and $\delta\phi$ and bounded as desired, so we need to show 
$\|R_{u,\phi}(\delta u,\delta\phi)\|_m=o(\|\delta u\|_{m+1}+\|\delta\phi\|_m)$ as $\|\delta u\|_{m+1}+\|\delta\phi\|_m\to 0$.
But, applying Lemma \ref{le:bound-u(phi)} and the algebra property, we can estimate the first term in $R_{u,\phi}$:
\[
\|\nabla(\delta u)\circ\phi \cdot \delta\phi\|_m\leq C\|\nabla(\delta u)\|_m\,\|\delta\phi\|_m \leq C \|\delta u\|_{m+1}\,\|\delta \phi\|_m.
\]
We can also use Lemma \ref{le:bound-u(phi)} and the algebra property to estimate the third term in $R_{u,\phi}$, namely $R_2$:
\[
\begin{aligned}
\|R_2(u,\phi,\delta u,\delta\phi)\|_m\leq \int_0^1 \|(\nabla(\delta u)\circ(\phi+t\delta\phi)-\nabla(\delta u)\circ\phi)\|_m\,dt\,\|\delta\phi\|_m \\
\leq 2 \sup_{0<t<1} \|\nabla(\delta u)\circ(\phi+t\delta\phi)\|_m\,\|\delta\phi\|_m \leq C\,\|\delta u\|_{m+1}\,\|\delta\phi\|_m.
\end{aligned}
\]
Using $\|\delta u\|_{m+1}\,\|\delta\phi\|_m\leq \|\delta u\|_{m+1}^2+\|\delta\phi\|_m^2$ in both estimates above, we see that the $X^m$-norms of the 
first and third terms
 in $R_{u,\phi}$ are actually $O( \|\delta u\|_{m+1}^2+\|\delta\phi\|_m^2)$ as $\|\delta u\|_{m+1}+\|\delta\phi\|_m\to 0$.

To estimate $R_1$ in the $X^m$-norm as $\|\delta \phi\|_m\to 0$, we use the continuity of $\phi\to\nabla u\circ \phi$ in $X^m$ from the first part of Proposition \ref{pr:compositioncontinuous} to conclude that $\|\nabla u \circ (\phi+t\delta\phi)-\nabla u\circ\phi\|_m=o(1)$ uniformly for $0<t<1$ as $\|\delta \phi\|_m\to 0$. Using this in the definition of $R_1(u,\phi,\delta\phi)$ we find that
\[
\|R_1(u,\phi,\delta\phi)\|_m=o(\|\delta\phi\|_m) \quad\hbox{as}\ \|\delta\phi\|_m\to 0.
\]
This completes the proof.
\hfill  $\Box$

\section{Proof of  Invertibility (Proposition \ref{pr:inverses})}

Before we begin the proof of  Proposition \ref{pr:inverses} we prove several lemmas that will be useful. The first shows invertibility near the identity; but we require one additional order of differentiability.

\begin{lemma}  \label{le:U-nbhd_id} 
For $m>1+d/p$ and any $\phi\in {\AH \mathcal D}^{m+1,p}_N$ with $\|\phi-Id\|_{\AH^{m+1,p}_N}<\e$ sufficiently small, then $\phi^{-1}\in {\AH \mathcal D}^{m+1,p}_N$. The same result holds when ${\AH \mathcal D}$ is replaced by  ${\mathcal {AD}}$. 
\end{lemma}

\noindent
{\bf Proof.} We proceed in two steps: we first use a fixed point argument to show $\phi^{-1}\in {\AH \mathcal D}^{m,p}_N$ and then we  show in fact that $\phi^{-1}\in {\AH \mathcal D}^{m+1,p}_N$. Having fixed $p$ and $N$, let us denote the norm in $ {\AH}^{m,p}_N$ simply by $\|\cdot\|_m$. 

To formulate the fixed point argument, let $X=\{ v\in {\AH}^{m,p}_N: \|v\|_m<\eta\}$, where we have chosen $\eta$ small enough that $Id+v\in {\AH \mathcal D}^{m,p}_N$ for all $v\in X$.
Now let us write $\phi=Id+u$ with $\|u\|_{m+1}<\e$, where $\e>0$ will be specified below, and let  $\psi:=\phi^{-1}$. Then $\psi$ is a diffeomorphism and $\phi\circ\psi=Id$
implies $\psi=Id-u\circ\psi.$ But we know that $u\in  {\AH}_{N}^{m+1,p}$ is bounded, so we can write $\psi=Id+v$, where $v:=-u\circ\psi$ is a bounded function; we want to show that $v\in {\AH}_{N}^{m,p}$. However, we know that $v$ satisfies
$Id+v=Id-u\circ(Id+v)$. Consequently, let us define
\begin{equation}\label{def:F_u}
F_u(w):= -u\circ(Id+w).
\end{equation}
If we can show that $F_u:X\to X$ and $F_u$ has a fixed point $w_*$, then we will have $\phi\circ(Id+w_*)=Id$. But we also know $\phi\circ(Id+v)=Id$. Applying $\phi^{-1}$ to both sides of $\phi\circ(Id+w_*)=\phi\circ(Id+v)$, we find $Id+w_*=Id+v$, showing $v=w_*\in X$ and hence $\phi^{-1}\in 
 {\AH \mathcal D}_{N}^{m,p}$. 

First let us confirm that $F_u:X\to X$. But we can use Lemma \ref{le:bound-u(phi)} to conclude $\|F_u(w)\|_m\leq C_1\|u\|_m$
where $C_1$ is locally uniform in $Id+w\in  {\AH \mathcal D}^{m,p}_N$, so can be taken uniform for $w\in X$. Taking $\e>0$ sufficiently small that $\e\, C_1\leq\eta$, we conclude that $\|u\|_m\leq\|u\|_{m+1}<\e$ implies $F_u:X\to X$. Next we show that $F_u:X\to X$ is a contraction. But we can use  (\ref{eq:u(phi)-u(phi*)-b}) to conclude
\[
\|F_u(v_1)-F_u(v_2)\|_m\leq C_2\,\|u\|_{m+1}\,\|v_1-v_2\|_m,
\]
where $C_2$ can be taken uniform for $v\in X$. Taking $\e>0$ small enough that $\e\, C_2<1$, we have $F_u:X\to X$ is a contraction. Consequently, $F_u$this map has a fixed point which must be $v$.

At this point we know $\phi^{-1}=Id+v\in {\AH \mathcal D}^{m,p}_N$ where $v$ satisfies 
\begin{equation}\label{eq:v(phi)=-u}
v\circ\phi=-u.
\end{equation}
We want to use this equation  and $u\in {\AH}_{N}^{m+1,p}$ to show that
$v\in {\AH}_{N}^{m+1,p}$. If we differentiate (\ref{eq:v(phi)=-u}), we obtain
\[
 Dv\circ\phi \, \cdot\,(I+Du)=-Du.
\]
Now  $\sup |Du|\leq C\,\|Du\|_{m-1}\leq C\, \|u\|_{m}<\epsilon\,C$, so  we can take $\e$ small enough that  $I+Du$ is an invertible matrix.  Using the adjoint form of the inverse matrix, we have
\begin{equation}\label{eq:cofactorformula}
(I+Du)^{-1}=\frac{1}{\det(I+Du)}\,{\rm Adj}(I+Du).
\end{equation}
Now $Du$ is a matrix, all of whose elements are in the Banach algebra $ {\AH}^{m,p}_{1,N}$; since ${\rm Adj}(I+Du)$ is comprised of products of these elements, it too is a matrix with elements in $ {\AH}^{m,p}_N$. Similarly, $\det(I+Du)$ is  a product of elements of $ {\AH}^{m,p}_N$, and $\det(I+Du)=1+w$ with $w\in {\AH}^{m,p}_{1,N}$. Using Lemma \ref{le:(1+u)^{-alpha}}, we know
$(\det(I+Du))^{-1}\in {\AH}^{m,p}_N$.
We conclude that the elements of $(I+Du)^{-1}$ are in ${\AH}^{m,p}_N$,
so the elements of  $Dv\circ\phi=- Du\,\cdot \, (I+Du)^{-1}$ are also in ${\AH}^{m,p}_N$.  But we know $\phi^{-1}\in {\AH \mathcal D}^{m,p}_N$, so we can compose  $Dv\circ\phi$ on the right by $\phi^{-1}$ to conclude
\[
Dv\circ\phi\circ\phi^{-1}=Dv\in  {\AH }^{m,p}_N.
\]
Consequently, $v\in  {\AH }^{m+1,p}_N$, which shows $\phi^{-1}\in {\AH \mathcal D}^{m+1,p}_N$.
\hfill  $\Box$

\medskip
The second lemma shows that any asymptotic diffeomorphism can be continuously deformed to one whose difference from the identity has compact support.

\begin{lemma}  \label{le:phi_t} 
For $m>1+d/p$ and any $\phi\in {\AH \mathcal D}^{m,p}_N$, there exists a continuous map $\phi_t: [0,1]\to {\AH \mathcal D}^{m,p}_N$ such that $\phi_1=\phi$ and $\phi_0-Id$ has compact support. The same result holds when ${\AH \mathcal D}$ is replaced by  ${\mathcal {AD}}$.
\end{lemma}

\noindent
{\bf Proof.} Write $\phi=Id+u$ where $u\in {\AH}^{m,p}_N$. Since $\phi$ is an orientation-preserving diffeomorphism, we know that
\[
\hbox{det}\, (D\phi)=\hbox{det}\, 
\begin{pmatrix} 
1+\del_1 u^1 & \del_2 u^1 & \cdots & \del_d u^1 \\
\del_1 u^2 & 1 + \del_2 u^2 & \cdots & \del_d u^2 \\
\vdots & \vdots & \vdots & \vdots \\
\del_1 u^d & \del_2 u^d & \cdots &1+ \del_d u^d
\end{pmatrix}
\ >\, 0 \quad\hbox{for all}\ x\in \RR^d.
\]
By continuity, we see that $\hbox{det}\,(D(\phi+v))>0$ provided $v$ is chosen so that $\sup |Dv|<\e$ with $\e>0$ sufficiently small.
But recall that we can write 
$u=\chi_R(r)\,(a_0(\theta)+\cdots r^{-N^*}a_{N^*}(\theta))+f$ where $N^*$ satisfies (\ref{H-asymptoticexpansion1}) and $f\in H_N^{m,p}$, so we can choose $R$ sufficiently large that
\[
\left| D\left[ \chi_R(r)\,(a_0(\theta)+\cdots r^{-N^*}a_{N^*}(\theta))+\chi_R(r)\,f(x)\right] \right| < \epsilon.
\]
(This can be done since $|\nabla (\chi_R(r))|=R^{-1}|\chi'(R^{-1}r)|\leq M/R$ and $|Df|=o(|x|^{-N})$ as $|x|\to\infty$.)
Thus, provided $R$ is sufficiently large, we see that
\[
\phi_t:= Id + u -(1- t)\,\left[\chi_R(a_0+\cdots \frac{a_{N^*}}{r^{N^*}})+\chi_R f\right]
\]
is an asymptotic diffeomorphism for all $t\in [0,1]$. But $\phi_1=\phi$ and $\phi_0=Id+(1-\chi_R)f$, where $(1-\chi_R)f$ has compact support, so we have proved the lemma.
\hfill  $\Box$

\medskip
The next lemma concerns  the differential of an asymptotic diffeomorphism $\phi\in {\AH \mathcal D}^{m+1,p}_N$ for $m>1+d/p$ as a $C^1$-map $ {\AH \mathcal D}^{m,p}_N\to  {\AH \mathcal D}^{m,p}_N$ defined by composition. Of course, we must take the differential of $\phi$ at a particular diffeomorphism $\psi$, which we take as the identity $\psi=Id$ and denote the resultant differential $d_\psi \phi$ by $d_0\phi$; this will be a linear map $d_0\phi: {\AH}^{m,p}_N \to {\AH}^{m,p}_N$. If we write $\phi=Id+u$, then we have  $d_0\phi=I+Du$, since for $v\in  {\AH}^{m,p}_N$ we can calculate pointwise
\[
\lim_{t\to 0} \frac{\phi(Id+tv)-\phi(Id)}{t}=v+Du\cdot v = (I+Du)\cdot v.
\]
We now show that this linear map $d_0\phi: {\AH}^{m,p}_N \to {\AH}^{m,p}_N$ is invertible. 

\begin{lemma}\label{le:dphi-invertible}  
For $m>1+d/p$ and any $\phi\in {\AH \mathcal D}^{m+1,p}_N$, the linear map $d_0\phi: {\AH}^{m,p}_N \to {\AH}^{m,p}_N$ is invertible.
If  ${\AH \mathcal D}$  is replaced by ${\mathcal {AD}}$ the conclusion 
becomes $d_0\phi: {\A}^{m,p}_{N+1} \to {\A}^{m,p}_{N+1}$ is invertible.
\end{lemma}

\noindent
{\bf Proof.} Write $\phi=Id+u$, where $u\in {\AH}^{m+1,p}_N$, so $d_0\phi= I+Du$. But since $\phi$ is a diffeomorphism, considered as a matrix, $I+Du$ is invertible and its inverse is given by the adjoint formula (\ref{eq:cofactorformula}).

Now $Du$ is a matrix, all of whose elements are in $ {\AH}^{m,p}_{1,N}$; since ${\rm Adj}(I+Du)$ is comprised of products of these elements, it too is a matrix with elements in $ {\AH}^{m,p}_N$. Since $m>d/p$, $ {\AH}^{m,p}_N$ is an algebra and ${\rm Adj}(I+Du)$ maps 
${\AH}^{m,p}_N \to {\AH}^{m,p}_N$. Also, $\det(I+Du)$ is  a product of elements of $ {\AH}^{m,p}_N$, so it too is in $ {\AH}^{m,p}_N$. Since $\det(I+Du)(x)>0$ for all $x\in \RR^d$ and $Du(x)\to 0$ as $|x|\to\infty$, we have $\det(I+Du)(x)\geq \e>0$; therefore, we can use Lemma \ref{le:(1+u)^{-alpha}} to conclude that
$\det(I+Du)^{-1}\in {\AH}^{m,p}_N$. Consequently, $(I+Du)^{-1}$ is a matrix with all elements in $ {\AH}^{m,p}_N$, so it is bounded as a map $ {\AH}^{m,p}_N \to {\AH}^{m,p}_N$, completing the proof.

The proof for ${\mathcal {AD}}$ is strictly analogous.
\hfill  $\Box$

\medskip
Next we want to show that left-translation by a fixed asymptotic diffeomorphism is an open map in a neighborhood of the identity; the next lemma shows that this is true provided we have one additional order of differentiability in the fixed diffeomorphism.

\begin{lemma}  
For $m>1+d/p$ and any fixed $\phi_*\in {\AH \mathcal D}^{m+1,p}_N$, there is an open neighborhood $ U$ of $Id$ in 
${\AH \mathcal D}^{m,p}_N$ such that $\phi_*( U)$ is an open neighborhood of $\phi_*$ in ${\AH \mathcal D}^{m,p}_N$.
The same result holds when ${\AH \mathcal D}$ is replaced by ${\mathcal {AD}}$.
\end{lemma}

\noindent
{\bf Proof.} Let $\phi_*=Id+u_*$ and for $\widetilde\phi=Id+\widetilde u$ near $Id$ in ${\AH \mathcal D}^{m,p}_N$, we can write
\[
\phi_*\circ\widetilde\phi=\phi_*+F_*(\widetilde u), \quad\hbox{where}\ F_*(\widetilde u)=u_*\circ(Id+\widetilde u)-u_*+\widetilde u.
\]
Hence there is an open neighborhood $U_0$ of $0$ in ${\AH}^{m,p}_N$ such that  $F_*:U_0\to{\AH}^{m,p}_N$ and $F_*(0)=0$; in fact, since $u_*\in  {\AH}^{m+1,p}_N$, by 
Proposition \ref{pr:compositioncontinuous}, $F_*:U_0\to{\AH}^{m,p}_N$ is $C^1$. If we compute the differential of $F_*$ at $0$, which we also denote by $d_0 F_*$, we find
\[
d_0 F_*(v)=Du_*\cdot v+v=(I+Du_*)\cdot v=d_0\phi_* (v) \quad\hbox{for any}\ v\in {\AH}^{m,p}_N.
\]
But, using Lemma \ref{le:dphi-invertible} (with $m$ in place of $m-1$), we know that $d_0\phi_*=d_0F_*$ is invertible. By the inverse function theorem, we conclude that $F_*$ admits a continuous inverse near $0\in {\AH}^{m,p}_N$, which translates as the desired conclusion for $\phi_*$ and $U=Id+U_0$.
\hfill  $\Box$

\medskip
Finally, we are ready to prove Proposition  \ref{pr:inverses}.

\medskip\noindent 
{\bf Proof of Proposition \ref{pr:inverses}.} We give the proof for ${\AH \mathcal D}$, the case of ${\mathcal A  \mathcal D}$ being  analogous. Using Lemma \ref{le:phi_t}, there exists a continuous map $\phi_t: [0,1]\to {\AH \mathcal D}^{m+1,p}_N$ such that $\phi_1=\phi$ and $\phi_0-Id$ has compact support. But then $\phi_0^{-1}$ is the identity outside a compact set, so trivially we have $\phi_0^{-1}\in {\AH \mathcal D}_{N}^{m+1,p}$. We want to use the continuity method to show that this property can gradually be extended to $\phi_1=\phi$.

Let us denote by $U^{m+1}$ the neighborhood of $Id$ in ${\AH \mathcal D}_{N}^{m+1,p}$ which Lemma \ref{le:U-nbhd_id} guarantees consists of diffeomorphisms that are invertible in ${\AH \mathcal D}_{N}^{m+1,p}$. 
Now let $U^{m+1}_t=\phi_t(U^{m+1})=\{\phi_t\circ\psi:\psi\in U^{m+1}\}$. We first want to show that every $\phi_*$ in 
$U^{m+1}_0$ is invertible in ${\AH \mathcal D}_{N}^{m+1,p}$. But this is trivial since $\phi_*=\phi_0\circ\psi$  implies
$\phi_*^{-1}=\psi^{-1}\circ\phi_0^{-1}\in{\AH \mathcal D}_{N}^{m+1,p}$. Now by compactness we can cover the path $\phi_t$ for $0\leq t \leq 1$ by a finite number of these translated neighborhoods, i.e.\ $U_0^{m+1},U_{t_1}^{m+1},\dots,U_{t_K}^{m+1},U_1^{m+1}$. Next we want to show every $\phi_*$ in 
$U^{m+1}_{t_1}$ is invertible in ${\AH \mathcal D}_{N}^{m+1,p}$; it suffices to show $\phi_{t_1}^{-1}\in {\AH \mathcal D}_{N}^{m+1,p}$. But we can pick $\tilde\phi\in U_0^{m+1}\cap U_{t_1}^{m+1}$, which is both  invertible in ${\AH \mathcal D}_{N}^{m+1,p}$ and of the form $\tilde\phi=\phi_{t_1}\circ\psi$ for some $\psi\in U^{m+1}$. However, we can compose with $\psi^{-1}$ to conclude $\phi_{t_1}=\tilde\phi\circ \psi^{-1}$ and hence $\phi_{t_1}^{-1}= \psi\circ \tilde\phi^{-1}\in {\AH \mathcal D}_{N}^{m+1,p}$. Clearly this process can be continued to show every $\phi_*\in U_1^{m+1}$ is invertible in ${\AH \mathcal D}_{N}^{m+1,p}$. In particular, $\phi_1=\phi$ is 
invertible in ${\AH \mathcal D}_{N}^{m+1,p}$, as desired.

Finally, we want to show that $\phi\mapsto\phi^{-1}$ is $C^1$ as a map ${\AH \mathcal D}_{N}^{m+1,p}\to {\AH \mathcal D}_{N}^{m,p}$. We shall do this using the implicit function theorem; this is valid since a neighborhood of a fixed $\phi_*$ in ${\AH \mathcal D}_{N}^{m+1,p}$ may be parameterized by a neighborhood of $0$ in the Banach space ${\AH}_{N}^{m+1,p}$. In fact, let $F:{\AH \mathcal D}_{N}^{m+1,p}\times {\AH \mathcal D}_{N}^{m,p} \to {\AH \mathcal D}_{N}^{m,p}$ represent composition, i.e.\ $F(\phi,\psi)=\phi\circ\psi$, which we know from Proposition \ref{pr:compositioncontinuous} is $C^1$. Let us fix $\phi_* \in {\AH \mathcal D}_{N}^{m+1,p}$ and consider the differential of the map
$\psi\mapsto F(\phi_*,\psi)$ at the point $\psi=\phi_*^{-1}\in {\AH \mathcal D}_{N}^{m,p}$,
which is given by
\[
T(h)=d\phi_*\circ\phi_*^{-1}\cdot h \quad \hbox{for}\ h\in  {\AH}_{N}^{m,p}.
\]
Since $d\phi_*\in {\AH}_{N}^{m,p}$, by Proposition \ref{pr:compositioncontinuous} we know that $d\phi_*\circ\phi_*^{-1}\in {\AH}_{N}^{m,p}$, and hence (using Proposition \ref{pr:A-products}) the linear operator $T$ is bounded ${\AH}_{N}^{m,p}\to {\AH}_{N}^{m,p}$.
In fact, $T$ is invertible on ${\AH}_{N}^{m,p}$, since its inverse is just $T^{-1}(h)=d(\phi_*)^{-1}\circ \phi_* \cdot h$, which is also bounded ${\AH}_{N}^{m,p}\to {\AH}_{N}^{m,p}$.
Consequently, the implicit function theorem implies that there is a neighborhood $U$ of $\phi_*$ in ${\AH \mathcal D}_{N}^{m+1,p}$ and a unique $C^1$ map 
$G:U\to {\AH \mathcal D}_{N}^{m,p}$ such that $F(\phi,G(\phi))=Id$ holds for all $\phi$ near $\phi_*$. But uniqueness of the inverse of $\phi$ shows $G(\phi)=\phi^{-1}$, and hence $\phi\mapsto\phi^{-1}$ is $C^1$.
\hfill  $\Box$

\appendix
\section{Proofs of Lemmas \ref{le:H-properties} and \ref{le:W-properties}}\label{A}

\noindent
{\bf Proof of  (c) \& (d) in Lemma \ref{le:H-properties}.}  Let $Q$ be a $d$-box of side length $1$. First suppose $d>mp$.
For $g\in H^{m,p}(Q)$ and $p\leq q\leq pd/(d-mp)$, the Sobolev inequality states
\[
\|g\|_{L^q(Q)}\leq C(d,m,p,q)\,\|g\|_{H^{m,p}(Q)}.
\]
Apply this to $g=\x^\de f$ and for $|\alpha|\leq m$ use
\[
\begin{aligned}
\|D^\alpha(\x^\de f)\|_{L^p(Q)}&=\|\sum_{\beta\leq\alpha} \begin{pmatrix} \alpha \\ \beta \end{pmatrix} D^\beta(\x^\delta)D^{\alpha-\beta} f\,\|_{L^p(Q)}\\
&\leq C(|\alpha|,\delta)\sum_{\beta\leq\alpha} \|\x^{\delta-|\beta|}D^{\alpha-\beta} f\,\|_{L^p(Q)}\\
&\leq C(|\alpha|,\delta)\sum_{|\gamma|\leq |\alpha|}  \|\x^{\delta}D^\gamma f\,\|_{L^p(Q)}
\end{aligned}
\]
to conclude
\begin{equation}\label{est:weightedSobolev-Q}
\|\x^\delta f\|_{L^q(Q)}\leq C(d,m,p,q,\delta)\sum_{|\alpha|\leq m} \|\x^\delta D^\alpha f\|_{L^p(Q)}.
\end{equation}
Now let $Q_0$ denote the $d$-box of side length $1$ centered at the origin in $\RR^d$, and $Q_\ell$ for $\ell=0,1,\dots$ be an enumeration of all $d$-boxes of side length 1 and centers at integral coordinates, so 
\[
\RR^d=\cup_{\ell=0}^\infty \overline{Q_\ell}.
\]
Then we use the inequality \eqref{est:weightedSobolev-Q} and then  the elementary inequality $\left(\sum a_j^q\right)^{1/q}\leq\left(\sum a_j^p\right)^{1/p}$ to estimate
\[
\begin{aligned}
\|f\|_{L^q_\delta(\RR^d)}&=
\|\x^\delta f\|_{L^q(\RR^d)}=\left(\sum_{\ell=0}^\infty \|\x^\delta f\|^q_{L^q(Q_\ell)}\right)^{1/q} \\
\leq \ & C\, \left(\sum_{\ell=0}^\infty\sum_{|\alpha|\leq m} \|\x^\delta D^\alpha f\|^q_{L^p(Q_\ell)}\right)^{1/q}\\
\leq \ & C\, \left(\sum_{\ell=0}^\infty\sum_{|\alpha|\leq m} \|\x^\delta D^\alpha f\|^p_{L^p(Q_\ell)}\right)^{1/p}
= \  C\, \left(\sum_{|\alpha|\leq m} \|\x^\delta D^\alpha f\|^p_{L^p(\RR^d)}\right)^{1/p},
\end{aligned}
\]
where $C=C(d,m,p,q,\delta)$.  But the last term is equivalent to $C \|f\|_{H^{m,p}_\delta}$, which establishes the inequality in (c).
The same argument works for $d=mp$, provided we assume $p\leq q<\infty$.

Now suppose $d<mp$ and $k<m-(d/p)$. For $g\in H^{m-k,p}(Q)$,  Morrey's inequality implies $g\in C(Q)$ and
\[
\sup_{x\in Q}| g(x)|\leq C(m,p,k) \,\|g\|_{H^{m-k,p}(Q)}.
\]
Apply this on $Q_\ell$ (as above) to $g=\x^\de D^\alpha f$ for any $\alpha$ satisfying $0\leq |\alpha|\leq k$:
\[
\begin{aligned}
\sup_{x\in Q_\ell}\, \x^\de |D^\alpha f(x)| & \leq C(m,p,k)\,\|\x^\de D^\alpha f\|_{H^{m-k,p}(Q_\ell)} \\
& = C(m,p,k) \,\sum_{|\beta|\leq m-k} \| D^\beta (\x^\de D^\alpha f) \|_{L^p(Q_\ell)} \\
& \leq C(m,p,k,\de) \,\sum_{|\alpha|\leq m} \| \x^\de D^\alpha f \|_{L^p(Q_\ell)} \\
& \leq C(m,p,k,\de)\, \|f\|_{H^{m,p}_\de(\RR^d)}.
\end{aligned}
\]
Now letting $\ell$ vary on the left, we obtain the desired inequality. Moreover, since the series
$\sum_\ell \|\x^\de D^\alpha f \|_{H^{m,p}(Q_\ell)}$ converges, we must have $ \|\x^\de D^\alpha f \|_{H^{m,p}(Q_\ell)}\to 0$
as $\ell\to\infty$. Consequently,
\[
\sup_{x\in Q_\ell}\, \x^\de |D^\alpha f(x)| \to 0 \quad\hbox{as}\ \ell\to\infty, 
\]
which implies that $|x|^\de |D^\alpha f(x)|\to 0$ as $|x|\to\infty$. \quad$\Box$

\medskip
\noindent
{\bf Proof of  (c) \& (d) in Lemma  \ref{le:W-properties}.}  We shall use scaling arguments as in \cite{B}.
Let us introduce the Sobolev norm with homogeneous weight function:
\begin{equation}\label{def:homogenousweightedSobolev}
|\!|\!| f |\!|\!|^p_{m,p,\delta}=\sum_{j=0}^m\int_{\RR^d\backslash\{0\}} \left||x|^{\delta+j}D^j f(x)\right|^p\,dx.
\end{equation}
For $R>0$ define $f_R(x)=f(Rx)$; it is easy to compute that
\begin{equation}\label{eq:scaledestimate}
|\!|\!| f_R |\!|\!|_{m,p,\de}=R^{-\de-\frac{d}{p}}|\!|\!| f |\!|\!|_{m,p,\de}.
\end{equation}
Letting $A_R=B_{2R}\backslash B_R$ where $B_R=\{x:|x|<R\}$, we can integrate over annuli instead of all of ${\RR}\backslash\{0\}$ (and adjust notation in the obvious way) to obtain the localized version of (\ref{eq:scaledestimate}):
\begin{equation}\label{eq:localscaledestimate}
 |\!|\!| f_R |\!|\!|_{m,p,\de;A_1} =R^{\,-\de-\frac{d}{p}}  |\!|\!| f |\!|\!|_{m,p,\de;A_R}.
\end{equation}
But the weighting factor $|x|^\delta$ is bounded above and below by constants on $A_1$ and by $c\,R^\de$ on $A_R$. We may  conclude the following equivalence:
\begin{equation}\label{scaledestimate-A_R}
 |\!|\!| f_R |\!|\!|_{m,p;A_1} \approx R^{\,-\de}  |\!|\!| f |\!|\!|_{m,p,\de-\frac{d}{p};A_R}.
\end{equation}

To prove (c), we apply the Sobolev inequality to $f_R$ on $A_1$ to conclude:
\[
\|f_R\|_{q;A_1}\leq C\,\|f_R\|_{m,p;A_1}, \quad C=C(d,m,p,q).
\]
Then apply (\ref{scaledestimate-A_R}) to both sides of this to conclude
\begin{equation}\label{eq:SobolevInequality-A_R}
 |\!|\!| f |\!|\!|_{q,\de-\frac{d}{q};A_R} \leq C \,  |\!|\!| f |\!|\!|_{m,p,\de-\frac{d}{p};A_R}, \quad C=C(d,m,p,q).
\end{equation}
Now let us write
\[
\RR^d=B_0\cup {\bf A}_1 \cup {\bf A}_2 \cup \cdots \quad \hbox{where}\ {\bf A}_j=A_{2^{j-1}}=B_{2^j}\backslash B_{2^{j-1}}.
\]
Then we can use the Sobolev inequality on $B_0$ and (\ref{eq:SobolevInequality-A_R}) on each ${\bf A}_j$ to obtain
\[
\begin{aligned}
\| f \|_{L^q_{\de-\frac{d}{q}}} & \approx \left( \|f\|^q_{q;B_0}+ |\!|\!| f  |\!|\!|^q_{q,\de-\frac{d}{q};{\bf A}_1} +\cdots+
 |\!|\!| f  |\!|\!|^q_{q,\de-\frac{d}{q};{\bf A}_k}+\cdots \right)^{1/q} \\
& \leq C\, \left( \|f\|^q_{m,p;B_0}+ |\!|\!| f  |\!|\!|^q_{m,p,\de-\frac{d}{q};{\bf A}_1} +\cdots+
 |\!|\!| f  |\!|\!|^q_{m,p,\de-\frac{d}{q};{\bf A}_k}+\cdots \right)^{1/q} \\
 & \leq C\, \left( \|f\|^p_{m,p;B_0}+ |\!|\!| f  |\!|\!|^p_{m,p,\de-\frac{d}{q};{\bf A}_1} +\cdots+
 |\!|\!| f  |\!|\!|^p_{m,p,\de-\frac{d}{q};{\bf A}_k}+\cdots \right)^{1/p}, \\
 \end{aligned}
\]
where in the last step we used $p\leq q$ and the elementary inequality $\left(\sum a_j^q\right)^{1/q}\leq\left(\sum a_j^p\right)^{1/p}$. But this last term is equivalent to $\|f\|_{W^p_{\de-\frac{d}{p}}}$, so we have proved (c).

To prove (d), when $mp>d$ the scaling argument yields in place of (\ref{eq:SobolevInequality-A_R})
\begin{equation}\label{eq:Morrey-A_R}
\sup_{x\in A_R} |x|^{\de+|\alpha|}|D^\alpha f(x)|\leq C\, |\!|\!| f  |\!|\!|_{m,p,\delta-\frac{d}{p};A_R}.
\end{equation}
Now, we can replace the right hand side of (\ref{eq:Morrey-A_R}) by $C\,\|f\|_{W_{\de-\frac{d}{p}}^{m,p}}$ and then allow $R$ on the left hand side to range freely to conclude
\[
\sup_{x\in\RR^d}\langle x\rangle^{\de+|\alpha|}|D^\alpha f(x)|\leq C\,\|f\|_{W_{\de-\frac{d}{p}}^{m,p}}.
\]
But since the series
\[
 |\!|\!| f  |\!|\!|^p_{m,p,\delta-\frac{d}{p};{\bf A}_1}+ |\!|\!| f  |\!|\!|^p_{m,p,\delta-\frac{d}{p};{\bf A}_2}+\cdots
\]
converges, we see that $ |\!|\!| f  |\!|\!|_{m,p,\delta-\frac{d}{p};{\bf A}_j}\to 0$ as $j\to\infty$, so from (\ref{eq:Morrey-A_R}) we conclude that 
\[
|x|^{\de+|\alpha|}|D^\alpha f(x)|\to 0 \quad\hbox{as}\ |x|\to\infty. \quad\Box
\]


\section{Asymptotic Spaces with Log terms}\label{B}

As we have seen (e.g.\ Example \ref{ex:d=N=2}), it is natural for log terms to arise when dealing with asymptotics.
In this section, we shall discuss a way to include log terms in the asymptotic function spaces in such a way as to still have Banach algebras and groups of diffeomorphisms. This will enable us to extend the Helmoltz decomposition
of Section 3 to $N\geq d$, and will be used in our application to the Euler equations with asymptotics
(cf.\ \cite{MT2}).

The key idea is to replace $a(x)$ in \eqref{asymptoticexpansion2} with 
\begin{equation}  \label{logasymptoticexpansion2} a(x)=\chi(r)\left(\frac{a_n^0+\cdots+ a_n^{n+\ell}\,(\log r)^{n+\ell}}{r^n}+\cdots + \frac{a_N^{0}+\cdots+a_N^{N+\ell}\,(\log r)^{N+\ell}}{r^N}\right),
\end{equation}
where $a^j_k\in H^{m+1+N-k,p}(S^{d-1})$ for $0\leq j\leq k+\ell$ and $0\leq n\leq k\leq N$. Note that $\ell$ is an integer, which can be negative, but we require $\ell\geq -n$. If we use as the remainder function  $f\in W_{\gamma_N}^{m,p}(\RR^d)$ with $\gamma_N$  satisfying \eqref{def:gamma}, we obtain a Banach space 
\begin{equation}\label{def:AWlog}
{\mathcal A}_{n,N;\ell}^{m,p}(\RR^d):= \{ \hbox{$u$ is in the form \eqref{asymptoticexpansion1} where $a$ is given 
by \eqref{logasymptoticexpansion2} and $f\in W_{\gamma_N}^{m,p}(\RR^d)$}\}
\end{equation}
with norm given by
\begin{equation}\label{def:AWlog-norm}
\|u\|_{{\mathcal A}_{n,N;\ell}^{m,p}}=\sum_{j=0}^{n+\ell}\|a_n^j\|_{H^{m+1+N-n}}+\cdots 
+\sum_{j=0}^{N+\ell}\|a^j_N\|_{H^{m+1}}+\|f\|_{W_{\gamma}^{m,p}(\RR^d)},
\end{equation}
where we have abbreviated $H^{m+1+N-k,p}(S^{d-1})$ by $H^{m+1+N-k}$.
We shall write ${\mathcal A}_{0,N;\ell}^{m,p}$ simply as ${\mathcal A}_{N;\ell}^{m,p}$. 
Note that ${\mathcal A}_{N}^{m,p}\subset {\mathcal A}_{N;\ell}^{m,p}$ for any $\ell\geq 0$, but that
 ${\mathcal A}_{N}^{m,p}\not={\mathcal A}_{N;0}^{m,p}$ (when $N\geq 1$). 

The asymptotic spaces ${\mathcal A}_{n,N;\ell}^{m,p}$ enjoy many properties analogous to those satisfied by the spaces ${\mathcal A}_{n,N}^{m,p}$. For example, the following is the analog of parts of Proposition \ref{pr:A-properties}:
\begin{proposition}\label{pr:Alog-properties}  
\begin{enumerate}
\item[(a)] If $n_1\geq n$, $N_1\geq N$, and $\ell_1\leq \ell$ then ${\mathcal A}_{n_1,N_1;\ell_1}^{m,p}\subset{\mathcal A}_{n,N;\ell}^{m,p}$.
\item[(b)] Multiplication by $\chi(r)\,(\log r)^{j}$ is bounded ${\mathcal A}_{n,N;\ell}^{m,p}\to{\mathcal A}_{n,N;\ell+j}^{m,p}$.
 \item[(c)] Multiplication by $\chi(r)\,r^{-k}$ is bounded ${\mathcal A}_{n,N;\ell}^{m,p}\to{\mathcal A}_{n+k,N+k;\ell-k}^{m,p}$.
  \item[(d)]     If $m\geq 1$,  then $u\mapsto \partial u/\partial x_j$ is continuous  ${\mathcal A}_{n,N;\ell}^{m,p}\to{\mathcal A}_{n+1,N+1;\ell-1}^{m-1,p}$.
 \item[(e)] Assume $m>d/p$. If $u\in {\mathcal A}_{n,N;\ell}^{m,p}$, then
 \[
 \sup_{x\in\RR^d} \frac{\x^{n+|\alpha|}}{(\log\x)^{n+\ell}}\,|D^\alpha u(x)|\leq C\, \|u\|_{{\mathcal A}_{n,N;\ell}^{m,p}} \quad\hbox{for all}\ |\alpha|<m-d/p.
 \]
\end{enumerate}
\end{proposition}
\noindent
We can also consider products of functions from these spaces. The following is the analog of parts
 of Proposition \ref{pr:A-products}:
\begin{proposition}\label{pr:Alog-products} 
For $m> d/p$,  $0\leq n_i\leq N_i$ and $\ell_i+n_i\geq 0$ for $i=1,2$,  let $n_0=n_1+n_2$, $N_0=\min(N_1+n_2,N_2+n_1)$, and $\ell_0=\ell_1+\ell_2$. Then
\begin{equation}\label{est:AWlog-multiplication}
\|u\,v\|_{{\mathcal A}_{n_0, N_0;\ell_0}^{m,p} }\leq C\,\|u\|_{{\mathcal A}_{n_1,N_1;\ell_1}^{m,p} }\|v\|_{{\mathcal A}_{n_2,N_2;\ell_2}^{m,p} }
\quad\hbox{for}\ u\in {\mathcal A}_{n_1,N_1;\ell_1}^{m,p},\ v\in {\mathcal A}_{n_2,N_2;\ell_2}^{m,p}.
\end{equation}
\end{proposition}
\noindent
As a consequence, we find that many of the ${\mathcal A}_{n,N;\ell}^{m,p}$ form Banach algebras: 
\begin{corollary} 
If $m> d/p$, $0\leq n\leq N$, and $-n\leq\ell\leq 0$, then ${\mathcal A}_{n,N;\ell}^{m,p}$ is a Banach algebra.
\end{corollary}

These function spaces are also useful for describing the mapping properties of the Laplacian (which is the reason that we have introduced them).
It is not difficult to confirm that, analogous to \eqref{Delta:A->A}, the following map is bounded:
\begin{equation}\label{eq:Laplacemap}
\Lap: {\mathcal A}_{N;\ell}^{m+1,p} \to {\mathcal A}_{2,N+2;\ell-2}^{m-1,p}.
\end{equation}
In confirming that \eqref{eq:Laplacemap} is bounded, the following computations are useful: for $\ell\geq 2$
\begin{subequations}\label{Lap(asymterm)}
\begin{equation}\label{Lap(asymterm)-ell}
\begin{aligned}
\Lap & \left[ \frac{a (\log r)^\ell}{r^k}\right]  \\
&=\frac{\ell(\ell-1)a(\log r)^{\ell-2}-\ell(2k+2-d)a(\log r)^{\ell-1}+(\Lap_\theta a+k(k+2-d) a)(\log r)^\ell }{r^{k+2}}
\end{aligned}
\end{equation}
while for $\ell=1,0$ we have
\begin{equation}\label{Lap(asymterm)-1}
\Lap  \left[ \frac{a (\log r)}{r^k}\right] =\frac{(d-2-2k)a+(\Lap_\theta a+k(k+2-d)a)\log r}{r^{k+2}}
\end{equation}
\begin{equation}\label{Lap(asymterm)-0}
\Lap  \left[ \frac{a }{r^k}\right] =\frac{\Lap_\theta a+k(k+2-d)a}{r^{k+2}}.
\end{equation}
\end{subequations}

Now we want to consider the invertibility of \eqref{eq:Laplacemap}. As in Section 3, we want to use separation of variables to solve equations of the form
\begin{equation}\label{eq:Lap(u)=}
\Lap u = \frac{b\,(\log r)^{j+\ell}}{r^{k+2}} \quad\hbox{for}\  b \in H^{m+N-k,p}(S^{d-1}),
\end{equation}
when $n\leq k\leq N$ and $-\ell\leq j\leq k$. The obvious candidate for a solution is in the form
\begin{equation}\label{u=guess1}
u=\frac{ a\,(\log r)^{j+\ell} } {r^k} \quad\hbox{for}\  a \in H^{m+2+N-k,p}(S^{d-1}),
\end{equation}
but let us see whether this is always successful. Using \eqref{Lap(asymterm)} with $\ell$ replaced by $j+\ell$, which we know is nonnegative, we have the following cases:
\begin{itemize}
\item $j+\ell=0$: Then \eqref{u=guess1} becomes
\begin{equation}\label{u=guess1.0}
u=\frac{a}{r^{k}},
\end{equation}
and from \eqref{Lap(asymterm)-0} we see that
$a$ must satisfy $\Lap_\theta a+k(k+2-d)a=b$. If $k<d-2$ (which is guaranteed if $N<d-2$), then  $\Lap_\theta +k(k+2-d)$ is invertible and we have solved \eqref{eq:Lap(u)=} with $u$ in the form \eqref{u=guess1}. But if $k\geq d-2$, then $\Lap_\theta +k(k+2-d)$ is not invertible, so we should modify our choice of solution \eqref{u=guess1.0}. Instead we try
\begin{equation}\label{u=guess1.1}
u=\frac{a_0+a_1\,\log r}{r^{k}} \quad\hbox{for}\  a_0,a_1 \in H^{m+2+N-k,p}(S^{d-1}),
\end{equation}
and from \eqref{Lap(asymterm)-1} we find that $a_0$ and $a_1$ must satisfy
\[
\begin{aligned}
\Lap_\theta a_1+k(k+2-d)a_1&=0 \quad\hbox{-- from  $\log r$ terms}\\
\Lap_\theta a_0+k(k+2-d)a_0&=b+(2k+2-d)a_1\quad\hbox{-- from no-$\log$ terms}.
\end{aligned}
\]
So we choose $a_1\in \hbox{ker}(\Lap_\theta +k(k+2-d))$ so that $b+(2k+2-d)a_1$ satisfies the solvability condition (possible since $2k+2-d>0$) enabling us to find $a_0$. Since \eqref{u=guess1.0} is a special case of \eqref{u=guess1.1}, we may take \eqref{u=guess1.1} as the general form of the solution when $j+\ell=0$.
\item $j+\ell=1$: Based upon our experience with the previous case, let us take as our general solution
\begin{equation}\label{u=guess2.1}
u=\frac{a_0+a_1\,\log r+a_2\,(\log r)^2}{r^{k}} \quad\hbox{for}\  a_0, a_1, a_2 \in H^{m+2+N-k,p}(S^{d-1}),
\end{equation}
Using \eqref{Lap(asymterm)-ell} and \eqref{Lap(asymterm)-1}, we find that $a_0$, $a_1$, and $a_2$ must satisfy
\[
\begin{aligned}
\Lap_\theta a_2+k(k+2-d)a_2&=0 \quad\hbox{-- from $(\log r)^2$ terms}\\
\Lap_\theta a_1+k(k+2-d)a_1&=b+2(2k+2-d)a_2\quad\hbox{-- from $\log r$ terms}\\
\Lap_\theta a_0+k(k+2-d)a_0&=(2k+2-d)a_1-2a_2 \quad\hbox{-- from no-$\log$ terms}.
\end{aligned}
\]
We  proceed as before: if $\Lap_\theta +k(k+2-d)$ is invertible, let $a_2=0$; if it is not invertible, then pick $a_2\in \hbox{ker}(\Lap_\theta +k(k+2-d))$ so that $b+2(2k+2-d)a_2$ satisfies the solvability condition to find $a_1$, then pick $a_1\in \hbox{ker}(\Lap_\theta +k(k+2-d))$ so that $(2k+2-d)a_1-2a_2$ satisfies the condition to find $a_0$.
\item $j+\ell\geq 2$: We now see the pattern, so we take
\begin{equation}\label{u=guess3.1}
u=\frac{a_0+a_1\,\log r+\cdots+a_{j+\ell+1}\,(\log r)^{j+\ell+1}}{r^{k}},
\end{equation}
and proceed as before to successively determine $a_{j+\ell+1}$, $a_{j+\ell}$, etc.
\end{itemize}
This analysis shows that, while \eqref{eq:Laplacemap} is not always invertible, we do have an inverse that is bounded into a larger space. The following result generalizes Proposition \ref{pr:inverseLaplacian2}:
\begin{proposition}\label{pr:inverseLaplacian3} 
For $m\geq 1$ there is a bounded map
\begin{equation}\label{eq:inverseLaplacemap}
K: {\mathcal A}_{2,N+2;\ell-2}^{m-1,p}  \to {\mathcal A}_{N;\ell+1}^{m+1,p}
\end{equation}
such that $\Lap Kv=v$ for all $v\in{\mathcal A}_{2,N+2;\ell-2}^{m-1,p}$.
\end{proposition}

The above proposition can be applied to finding a Helmholz decomposition when $N\geq d$. In fact, if we begin with 
${\bf u}\in {\mathcal A}_{1,N;\ell}^{m,p}$ for any $N\geq 1$ and $\ell\geq -1$, and follow the recipe described in Section 3, we will obtain the decomposition \eqref{eq:Helmholtz} where ${\bf  v},\nabla w\in  {\mathcal A}_{1,N;\ell+1}^{m,p}$, i.e.\ they have an additional log term.
This generalizes the phenomenon encountered in Example \ref{ex:d=N=2}.

Finally, in the full spirit of this paper, we want to consider diffeomorphisms of $\RR^d$ whose asymptotic behavior as $|x|\to\infty$ is given by ${\mathcal A}_{n,N;\ell}^{m,p}$. As in Definition \ref{def:asym-diffeos}, for $m>1+d/p$ we define
\[
\A\D_{n,N;\ell}^{m,p}(\RR^d,\RR^d):=\{\phi\in\Diff_+^1(\RR^d,\RR^d)\,|\,\phi(x)=x+u(x), \ u\in\A_{n,N;\ell}^{m,p}(\RR^d,\RR^d)\}\,.
\]
First consider the continuity of composition. 
\begin{proposition}\label{pr:compositioncontinuous-log} 
For integers $m> 1+d/p$, $0\leq n\leq N$, and $-n\leq\ell\leq 0$, composition $(u,\psi)\mapsto u\circ\psi$ defines  a continuous mapping
\begin{subequations}\label{eq:compositionAxAD->A-log}
\begin{equation}\label{eq:composition1-log}
 {\mathcal A}_{n,N;\ell}^{m,p} \times {\mathcal A \mathcal D}_{N;\ell}^{m,p} \to  {\mathcal A}_{n,N;\ell}^{m,p},
\end{equation}
and a $C^1$-mapping
\begin{equation}\label{eq:composition2-log}
 {\mathcal A}_{n,N;\ell}^{m+1,p} \times {\mathcal A \mathcal D}_{N;\ell}^{m,p} \to  {\mathcal A}_{n,N;\ell}^{m,p}.
\end{equation}
\end{subequations}
\end{proposition}
 \noindent
 The proof of Proposition \ref{pr:compositioncontinuous-log} is analogous to that of Proposition \ref{pr:compositioncontinuous}, so we will not give all the details, but let us point out some of the differences in the proof. We again use the functions $\rho$ and $\sigma$ defined in \eqref{def:rho} and \eqref{def:sigma} respectively, but we consider them as maps on our log-asymptotic spaces: we can verify that
 $\rho: {\mathcal A}^{m,p}_{N;\ell}(B_1^c,\RR^d)\to {\mathcal A}^{m,p}_{1,N+1;\ell-1}(B_1^c)$ is continuous and 
 for a fixed $u^*\in {\mathcal A}^{m,p}_{N;\ell}(B_1^c,\RR^d)$ there is a neighborhood ${\mathcal U}$ for which 
 $\sigma: {\mathcal U} \to {\mathcal A}^{m,p}_{N;\ell}(B_R^c)$ is real analytic provided $R$ is sufficiently large.
 We need to replace Lemma \ref{le:bound-u(phi)} with the following statement: For $m>1+d/p$,
 $0\leq n\leq N$, $-n\leq \ell\leq 0$, and $\phi\in {\mathcal AD}^{m,p}_{N;\ell}$, we have
 \begin{equation}\label{est:u(phi)inA-log}
 \|u\circ\phi\|_{{\mathcal AD}^{m,p}_{n,N;\ell}} \leq C\,\|u\|_{{\mathcal AD}^{m,p}_{n,N;\ell}}
 \quad\hbox{for all}\ u\in {\mathcal AD}^{m,p}_{n,N;\ell},
 \end{equation}
 where $C$ may be taken locally uniformly in $\phi\in {\mathcal AD}^{m,p}_{N;\ell}$.
 To prove this result, it  suffices to consider $u$ of the form
 \begin{equation}\label{a=ak-log}
 u(x)=\chi(r)\frac{a_k(\theta) (\log r)^j}{r^k} \quad\hbox{where}\ a_k\in H^{m+1+N-k,p}(S^{d-1})
 \ \hbox{and}\ 0\leq j\leq k.
 \end{equation}
 For $j=0$, the proof of  \eqref{est:u(phi)inA-log} is exactly the same as for \eqref{est:u(phi)inA} since we 
 still have \eqref{eq:ak-Taylor3}, \eqref{eq:|phi|^(-1)}, and \eqref{eq:phi/|phi|-x/|x|}. To handle $j>0$,
 let us first define
 \begin{equation}\label{def:tau}
 \tau(v):= \frac{1}{2}\log(1+\rho(v)).
 \end{equation}
 Since $\rho(v)\to 0$ as $|x|\to\infty$, we can use the expansion $\log(1+\rho)=\rho-\frac{1}{2}\rho^2+\frac{1}{3}\rho^3+\cdots$ near $\rho=0$ to conclude, similar to Lemma \ref{le:sigma-continuity} for $\sigma$, that  any
 $u^*\in {\mathcal A}_{N;\ell}^{m,p}(B_1^c,\RR^d)$ admits a neighborhood ${\mathcal U}$ such that the map 
 $\tau:{\mathcal U}\to {\mathcal A}_{1,N+1;\ell-1}^{m,p}(B_R^c)$ is real analytic if $R$ is sufficiently large;
 in particular, the mapping is continuous.
Now, for $\phi=Id+v \in {\mathcal AD}_{N;\ell}^{m,p}$, we can write
 \[
 \log|\phi(x)|=\log|x+v(x)|=\log(|x|(1+\rho(v))^{1/2})=\log|x|+\tau(v)
 \]
 and so for $j>0$ we have
 \begin{equation}
 ( \log|\phi(x)|)^j=(\log|x|)^j+j\,\tau(v)\,(\log|x|)^{j-1}+\cdots+(\tau(v))^j.
 \end{equation}
 This enables us to complete the arguments in the proof of Lemma \ref{le:bound-u(phi)} to show that 
 \eqref{est:u(phi)inA-log} holds. It also enables us to use the arguments in the proof of Lemma
 \ref{le:u(phi_j)->u(phi)} to show that, for $u$ as in \eqref{a=ak-log}, we have the following:
 if $\phi,\phi_j\in {\mathcal AD}^{m,p}_{N;\ell}$ with $\phi-\phi_j\to 0$ in ${\mathcal A}^{m,p}_{N;\ell}$
 and $u\in {\mathcal A}^{m,p}_{n,N;\ell}$, then $u\circ\phi_j\to u\circ\phi$ in ${\mathcal A}_{n,N;\ell}^{m,p}$.
 These results show that \eqref{eq:composition1-log} is continuous. To show that \eqref{eq:composition2-log} is
 $C^1$, we can follow the proof of \eqref{eq:u(phi)-u(phi*)-b} to obtain the following extension: for $u\in {\mathcal A}_{n,N;\ell}^{m+1,p}$,
 $\phi_*\in{\mathcal AD}_{N;\ell}^{m,p}$, and all $\phi\in {\mathcal AD}_{N;\ell}^{m,p}$ sufficiently close to $\phi_*$ we have 
 \begin{equation}\label{eq:u(phi)-u(phi*)-log}
\|u\circ \phi-u\circ \phi_* \|_{{\mathcal A}_{n,N;\ell}^{m,p}}
\leq C\, \|u\|_{{\mathcal A}_{n,N;\ell}^{m+1,p}} \|\phi-\phi_*\|_{{\mathcal A}_{N;\ell}^{m,p}}.
 \end{equation}
 From here, the proof that \eqref{eq:composition2-log} is $C^1$ follows exactly as for Proposition  
 \ref{pr:compositioncontinuous}.
 
Next consider the issue of inverses.
The following result is analogous to Proposition \ref{pr:inverses}:
\begin{proposition}\label{pr:inverses-log} 
For integers $m> 1+d/p$, $0\leq n\leq N$, and $-n\leq\ell\leq 0$, if $\phi\in{\mathcal A \mathcal D}_{n,N;\ell}^{m+1,p}$  then
 $\phi^{-1}\in{\mathcal A \mathcal D}_{n,N;\ell}^{m+1,p}$, and $\phi\to\phi^{-1}$ defines a $C^1$-mapping 
${\mathcal A \mathcal D}_{n,N;\ell}^{m+1,p}\to {\mathcal A \mathcal D}_{n,N;\ell}^{m,p}$.
 \end{proposition}
 \noindent
 The proof of Proposition \ref{pr:inverses} genearlizes immediately to establish
 Proposition \ref{pr:inverses-log}, so we shall not discuss details.
 
Finally, combining Propositions \ref{pr:compositioncontinuous-log} and \ref{pr:inverses-log} (and using \cite{M}), we obtain 
the analog of Theorem \ref{th:topgroup}:
\begin{theorem}\label{th:topgroup-log}  
For integers $m> 2+d/p$, $0\leq n\leq N$, and $-n\leq\ell\leq 0$,  ${\mathcal A \mathcal D}_{n,N;\ell}^{m,p}$ is a topological group under composition.
\end{theorem}

\begin{remark}
If we use as remainder space $H^{m,p}_N(\RR^d)$ instead of 
$W_{\gamma_N}^{m,p}(\RR^d)$ in \eqref{def:AWlog}, we will get a family of asymptotic function spaces ${\AH}_{n,N;\ell}^{m,p}(\RR^d)$
analogous to ${\AH}_{n,N}^{m,p}(\RR^d)$ but now with log terms. It enjoys most of the properties that we have discussed for the family ${\mathcal A}_{n,N;\ell}^{m,p}(\RR^d)$, but is not so useful for describing the mapping properties of $\Lap$, so we shall not discuss it further.
\end{remark}



\begin{thebibliography}{999}

\bibitem{A} V. Arnold, {\em Sur la geometri{\'e} differentielle des groupes de Lie
de dimension infinie et ses applications {\`a} l'hydrodynamique des fluids parfaits},
Ann. Inst. Fourier, $\bf 16$, 1(1966), 319-361.

\bibitem{B} R. Bartnik, \textit{The mass of an asymptotically flat manifold}, Comm. Pure Appl. Math., \textbf{39} (1986), no. 5, 661-693.

\bibitem{BS1} I.N. Bondareva, M. Shubin, {\em Uniqueness of the solution of the Cauchy problem for the Korteweg - de Vries equation in classes of increasing functions}, Moscow Univ. Math. Bulletin, \textbf{40} (1985), 53-57.

\bibitem{BS2} I.N. Bondareva, M. Shubin, {\em Equations of Korteweg-de Vries type in classes of increasing functions},  J. Soviet Math. \textbf{51} (1990), no. 3, 2323-2332.

\bibitem{BB} J.P. Bourguignon, H. Brezis, {\em Remarks on the Euler equation}, J. Func. Anal., \textbf{15} (1974), 341-363.

\bibitem{CH} R. Camassa, D. Holm, {\em An integrable shallow water equation with
peaked solitons}, Phys. Rev. Lett., \textbf{71} (1993), 1661-1664.

\bibitem{C} M. Cantor, \textit{Perfect fluid flows over $\RR^n$ with asymptotic conditions}, J. Func. Anal., \textbf{18} (1975), 73-84.

\bibitem{Co} A. Constantin, \textit{Existence of permanent and breaking waves for a shallow water equation: a geometric approach}, Ann. Inst. Four. Grenoble, \textbf{50} (2000), 321-362.

\bibitem{EM} D. Ebin, J. Marsden, {\em Groups of diffeomorphisms and the motion
of an incompressible fluid}, Ann. Math., $\bf 92$ (1970), 102-163

\bibitem{IKT} H. Inci, T. Kappeler, P. Topalov, {\em On the regularity of the composition of diffeomorphisms},
 Mem. Amer. Math. Soc., 226 (2013), no.\ 1062.

\bibitem{KPST} T. Kappeler, P. Perry, M. Shubin, P. Topalov, {\em Solutions of mKdV in classes of functions unbounded at infinity}, J. Geom. Anal. \textbf{18} (2008), no. 2, 443�477.


\bibitem{KPV} C. Kenig, G. Ponce, L. Vega, {\em Global solutions for the KdV equation with unbounded data}, J. Diff. Equations, {\bf 139} (1997), 339-364.

\bibitem{Mc} R. McOwen, \textit{The behavior of the Laplacian on weighted Sobolev spaces,} Comm. Pure Appl. Math. \textbf{32} (1979), 783-795. 

\bibitem{MT} R. McOwen, P. Topalov, {\em Asymptotics in shallow water waves}, Discrete Contin. Dyn. Syst., {\bf 35} (2015), 3103-3131.

\bibitem{MT2} R. McOwen, P. Topalov, {\em Euler's equation with asymptotics on Euclidean spaces}, in preparation.

\bibitem{Men} A. Menikoff, {\em The existence of unbounded solutions of the Korteweg-de Vries equation}, Comm. Pure Appl. Math, \textbf{25} (1972), 407-432.

\bibitem{Mi} G. Misiolek, {\em A shallow water equation as a geodesic flow on the Bott-Visaro group} J. Geom. Phys., {\bf 24} (1998), 203-208.

\bibitem{M} D. Montgomery, \textit{On continuity in topological groups}, Bull. Amer. Math. Soc., \textbf{42} (1936), 879-882.


\end{thebibliography}
\end{document}